\documentclass{article}
\usepackage{amsmath}
\usepackage{amssymb}
\begin{document}
\def\e#1\e{\begin{equation}#1\end{equation}}
\def\ea#1\ea{\begin{align}#1\end{align}}
\def\eq#1{{\rm(\ref{#1})}}
\newtheorem{thm}{Theorem}[section]
\newtheorem{lem}[thm]{Lemma}
\newtheorem{prop}[thm]{Proposition}
\newtheorem{conj}[thm]{Conjecture}
\newtheorem{cor}[thm]{Corollary}
\newenvironment{dfn}{\medskip\refstepcounter{thm}
\noindent{\bf Definition \thesection.\arabic{thm}\ }}{\medskip}
\newenvironment{ex}{\medskip\refstepcounter{thm}
\noindent{\bf Example \thesection.\arabic{thm}\ }}{\medskip}
\newenvironment{proof}[1][,]{\medskip\ifcat,#1
\noindent{\it Proof.\ }\else\noindent{\it Proof of #1.\ }\fi}
{\hfill$\square$\medskip}
\def\dim{\mathop{\rm dim}}
\def\Re{\mathop{\rm Re}}
\def\Im{\mathop{\rm Im}}
\def\ind{\mathop{\rm ind}}
\def\vol{\mathop{\rm vol}}
\def\area{\mathop{\rm area}}
\def\sind{{\ts\mathop{\text{\rm s-ind}}}}
\def\id{\mathop{\rm id}}
\def\U{\mathbin{\rm U}}
\def\SU{\mathop{\rm SU}}
\def\SO{\mathop{\rm SO}}
\def\ge{\geqslant} 
\def\le{\leqslant} 
\def\N{\mathbin{\mathbb N}}
\def\R{\mathbin{\mathbb R}}
\def\Z{\mathbin{\mathbb Z}}
\def\Q{\mathbin{\mathbb Q}}
\def\C{\mathbin{\mathbb C}}
\def\al{\alpha}
\def\be{\beta}
\def\ga{\gamma}
\def\de{\delta}
\def\ep{\epsilon}
\def\la{\lambda}
\def\ka{\kappa}
\def\th{\theta}
\def\ze{\zeta}
\def\up{\upsilon}
\def\vp{\varphi}
\def\si{\sigma}
\def\om{\omega}
\def\De{\Delta}
\def\La{\Lambda}
\def\Om{\Omega}
\def\Ga{\Gamma}
\def\Si{\Sigma}
\def\Th{\Theta}
\def\Up{\Upsilon}
\def\d{{\rm d}}
\def\pd{\partial}
\def\ts{\textstyle}
\def\sst{\scriptscriptstyle}
\def\w{\wedge}
\def\sm{\setminus}
\def\op{\oplus}
\def\ot{\otimes}
\def\iy{\infty}
\def\ra{\rightarrow}
\def\longra{\longrightarrow}
\def\t{\times}
\def\na{\nabla}
\def\ha{{\textstyle\frac{1}{2}}}
\def\ti{\tilde}
\def\bs{\boldsymbol}
\def\ov{\overline}
\def\sC{{\smash{\sst C}}}
\def\sCi{{\smash{\sst C_i}}}
\def\sE{{\smash{\sst\cal E}}}
\def\sL{{\smash{\sst L}}}
\def\sLi{{\smash{\sst L_i}}}
\def\sLon{{\smash{\sst L_1,\ldots,L_n}}}
\def\sSi{{\smash{\sst\Si}}}
\def\sSii{{\smash{\sst\Si_i}}}
\def\sN{{\smash{\sst N}}}
\def\sX{{\smash{\sst X}}}
\def\sNt{{\smash{\sst N^t}}}
\def\sNst{{\smash{\sst N^{s,t}}}}
\def\sW{{\smash{\sst W}}}
\def\sWt{{\smash{\sst W^t}}}
\def\sWst{{\smash{\sst W^{s,t}}}}
\def\sYi{{\smash{\sst{\cal Y}_i}}}
\def\sZi{{\smash{\sst{\cal Z}_i}}}
\def\sXp{{\smash{\sst X'}}}
\def\sF{{\smash{\sst\cal F}}}
\def\D{{\cal D}}
\def\F{{\cal F}}
\def\G{{\cal G}}
\def\I{{\cal I}}
\def\M{{\cal M}}
\def\oM{{\,\,\ov{\!\!{\cal M}\!}\,}}
\def\O{{\cal O}}
\def\cY{{\cal Y}}
\def\cZ{{\cal Z}}
\def\ms#1{\vert#1\vert^2}
\def\md#1{\vert #1 \vert}
\def\bmd#1{\big\vert #1 \big\vert}
\def\lnm#1#2{\Vert #1 \Vert_{L^{#2}}}
\def\an#1{\langle#1\rangle}
\def\ban#1{\bigl\langle#1\bigr\rangle}
\title{Special Lagrangian submanifolds with isolated\\
conical singularities. V. Survey and applications}
\author{Dominic Joyce \\ Lincoln College, Oxford}
\date{}
\maketitle

\section{Introduction}
\label{cs1}

{\it Special Lagrangian $m$-folds (SL\/ $m$-folds)} are a
distinguished class of real $m$-dimensional minimal submanifolds
which may be defined in $\C^m$, or in {\it Calabi--Yau $m$-folds},
or more generally in {\it almost Calabi--Yau $m$-folds} (compact
K\"ahler $m$-folds with trivial canonical bundle). We write an
almost Calabi--Yau $m$-fold as $M$ or $(M,J,\om,\Om)$, where
the manifold $M$ has complex structure $J$, K\"ahler form
$\om$ and holomorphic volume form~$\Om$.

This is the fifth in a series of five papers
\cite{Joyc9,Joyc10,Joyc11,Joyc12} studying SL $m$-folds with
{\it isolated conical singularities}. That is, we consider an SL
$m$-fold $X$ in an almost Calabi--Yau $m$-fold $M$ for $m>2$ with
singularities at $x_1,\ldots,x_n$ in $M$, such that for some
special Lagrangian cones $C_i$ in $T_{\smash{x_i}}M\cong\C^m$
with $C_i\sm\{0\}$ nonsingular, $X$ approaches $C_i$ near $x_i$,
in an asymptotic $C^1$ sense.

New readers of the series are advised to begin with this paper. We
shall survey the major results of \cite{Joyc9,Joyc10,Joyc11,Joyc12},
giving explanations, but avoiding the long, technical analytic
proofs of previous papers. We also integrate the results to give
an (incomplete) description of the {\it boundary} of a moduli
space of compact SL $m$-folds, and apply them to prove some
conjectures in \cite{Joyc1,Joyc5} on {\it connected sums} of SL
$m$-folds, and $T^2$-{\it cone singularities} of SL 3-folds.

Having a good understanding of the singularities of special
Lagrangian submanifolds will be essential in clarifying the
Strominger--Yau--Zaslow conjecture on the Mirror Symmetry
of Calabi--Yau 3-folds \cite{SYZ}, and also in resolving
conjectures made by the author \cite{Joyc1} on defining
new invariants of Calabi--Yau 3-folds by counting special
Lagrangian homology 3-spheres with weights. The series aims
to develop such an understanding for simple singularities
of SL $m$-folds.

We begin in \S\ref{cs2} with an introduction to almost
Calabi--Yau and special Lagrangian geometry, and the
{\it deformation theory} of compact SL $m$-folds. Section
\ref{cs3} defines {\it SL\/ $m$-folds with conical
singularities}, our subject, gives examples of
{\it special Lagrangian cones}, and some basics on
homology and cohomology.

Section \ref{cs4} describes the first paper \cite{Joyc9}
on the {\it regularity} of SL $m$-folds $X$ with conical
singularities $x_1,\ldots,x_n$. We study the asymptotic
behaviour of $X$ and its derivatives near $x_i$, how
quickly it converges to the cone~$C_i$.

In \S\ref{cs5} we discuss the second paper \cite{Joyc10} on
the {\it deformation theory} of compact SL $m$-folds $X$ with
conical singularities in an almost Calabi--Yau $m$-fold $M$. We
find that the {\it moduli space} ${\cal M}_\sX$ of deformations
of $X$ in $M$ is locally homeomorphic to the zeroes of a smooth map
$\Phi:{\cal I}_\sXp\ra{\cal O}_\sXp$ between finite-dimensional
vector spaces, and if the {\it obstruction space} ${\cal O}_\sXp$
is zero then ${\cal M}_\sX$ is a smooth manifold. We also study
deformations in {\it smooth families} of almost Calabi--Yau
$m$-folds $(M,J^s,\om^s,\Om^s)$ for~$s\in\F\subset\R^d$.

Section \ref{cs6} is an aside on {\it Asymptotically Conical
SL\/ $m$-folds (AC SL\/ $m$-folds)} in $\C^m$, that is,
nonsingular, noncompact SL $m$-folds $L$ in $\C^m$ which are
asymptotic at infinity to an SL cone $C$ at a prescribed
{\it rate} $\la$. Our main sources are \cite[\S 7]{Joyc9} and
Marshall \cite{Mars}. The theories of AC SL $m$-folds and SL
$m$-folds with conical singularities are similar in many respects.

Section \ref{cs7} explains the third and fourth papers
\cite{Joyc11,Joyc12} on {\it desingularizations} of a compact
SL $m$-fold $X$ with conical singularities $x_i$ with cones
$C_i$ for $i=1,\ldots,n$ in an almost Calabi--Yau $m$-fold
$M$. We take AC SL $m$-folds $L_i$ in $\C^m$ asymptotic to
$C_i$ at infinity, and glue $tL_i$ into $X$ at $x_i$ for small
$t>0$ to get a smooth family of compact, {\it nonsingular}
SL $m$-folds $\smash{\ti N^t}$ in $M$, with $\smash{\ti N^t}
\ra X$ as $t\ra 0$. We also study desingularizations in {\it
families} of almost Calabi--Yau $m$-folds $(M,J^s,\om^s,\Om^s)$
for~$s\in\F$.

The new material of the paper is \S\ref{cs8}--\S\ref{cs10}.
We study the moduli space $\M_\sN$ of compact, nonsingular
SL $m$-folds $N$ in \S\ref{cs2}, the moduli space $\M_\sX$
of compact SL $m$-folds $X$ with conical singularities in
\S\ref{cs5}, and the moduli space $\M_\sL^\la$ of AC SL $m$-folds
in $\C^m$ with rate $\la$ in \S\ref{cs6}. Section \ref{cs8}
explains how these three kinds of moduli space fit together.

The idea is that $\M_\sN$ has a {\it compactification} $\oM_\sN$
with {\it boundary} $\pd\M_\sN=\oM_\sN\sm\M_\sN$ consisting of
singular SL $m$-folds. Suppose $N$ is constructed as in
\S\ref{cs7} by desingularizing $X$ with conical singularities
$x_1,\ldots,x_n$ by gluing in AC SL $m$-folds $L_1,\ldots,L_n$.
Then in good cases we expect $\M_\sX\subseteq\pd\M_\sN$, and
$\M_\sN$ may be locally modelled on a subset of $\M_\sX\t
\M_{\smash{\sst L_1}}^0\t\cdots\t\M_{\smash{\sst L_n}}^0$ near~$X$.

Section \ref{cs9} considers {\it connected sums} of SL
$m$-folds. Suppose $X$ is a compact, {\it immersed\/} SL
$m$-fold in $M$ with {\it transverse self-intersection
points} $x_1,\ldots,x_n$. This includes the case where $X$
is a union of $q>1$ embedded SL $m$-folds $X_1,\ldots,X_q$,
and the $x_i$ are intersections between $X_j$ and $X_k$. Then
we can consider $X$ to be an {\it SL $m$-fold with conical
singularities}, with each cone $C_i$ the union of two transverse
{\it SL\/ $m$-planes} $\Pi_i^+,\Pi_i^-\cong\R^m$ in~$\C^m$.

When $\Pi_i^\pm$ satisfy an {\it angle criterion}, Lawlor
\cite{Lawl} constructed a family of AC SL $m$-folds
$L_i^{\smash{\pm,A}}$ for $A>0$ with cone $\Pi_i^+\cup\Pi_i^-$,
diffeomorphic to ${\cal S}^{m-1}\t\R$. We apply the results of
\S\ref{cs7} to construct SL $m$-folds $\smash{\ti N^t}$ by gluing
$tL_i^{\smash{\pm,A}}$ into $X$ at $x_i$. These $\smash{\ti N^t}$
are {\it multiple connected sums} of $X$ with itself. In this way
we re-prove and extend a result of Lee~\cite{Lee}.

Finally, \S\ref{cs10} studies SL 3-folds $X$ with conical
singularities with cone $C$ the $\U(1)^2$-invariant SL
$T^2$-cone due to Harvey and Lawson \cite[\S III.3.A]{HaLa}.
These have particularly nice properties. For instance, the
moduli space $\M_\sX$ is always smooth, and under topological
conditions the compactified moduli space $\oM_\sN$
of desingularizations $N$ is near $X$ a nonsingular manifold
with boundary $\M_\sX$. We prove several conjectures
from~\cite{Joyc1,Joyc5}.
\medskip

\noindent{\it Acknowledgements.} I would like to thank
Stephen Marshall, Mark Haskins, Tadashi Tokieda, Ivan
Smith, Sema Salur, and Adrian Butscher for useful
conversations. I was supported by an EPSRC Advanced
Research Fellowship whilst writing this paper.

\section{Special Lagrangian geometry}
\label{cs2}

We begin with some background from symplectic geometry.
Then special Lagrangian submanifolds (SL $m$-folds) are
introduced both in $\C^m$ and in {\it almost Calabi--Yau
$m$-folds}. We also describe the {\it deformation theory}
of compact SL $m$-folds. Some references for this section
are McDuff and Salamon \cite{McSa}, Harvey and Lawson
\cite{HaLa}, McLean \cite{McLe}, and the author~\cite{Joyc8}.

\subsection{Background from symplectic geometry}
\label{cs21}

We start by recalling some elementary symplectic geometry, which
can be found in McDuff and Salamon \cite{McSa}. Here are the basic
definitions.

\begin{dfn} Let $M$ be a smooth manifold of even dimension $2m$.
A closed $2$-form $\om$ on $M$ is called a {\it symplectic form}
if the $2m$-form $\om^m$ is nonzero at every point of $M$. Then
$(M,\om)$ is called a {\it symplectic manifold}. A submanifold
$N$ in $M$ is called {\it Lagrangian} if $\dim N=m=\ha\dim M$
and~$\om\vert_N\equiv 0$.
\label{cs2def1}
\end{dfn}

The simplest example of a symplectic manifold is~$\R^{2m}$.

\begin{dfn} Let $\R^{2m}$ have coordinates $(x_1,\ldots,x_m,
y_1,\ldots,y_m)$, and define the standard metric $g'$ and
symplectic form $\om'$ on $\R^{2m}$ by
\e
g'=\sum_{j=1}^m(\d x_j^2+\d y_j^2)\quad\text{and}\quad
\om'=\sum_{j=1}^m\d x_j\w\d y_j.
\label{cs2eq1}
\e
Then $(\R^{2m},\om')$ is a symplectic manifold. When we wish to
identify $\R^{2m}$ with $\C^m$, we take the complex coordinates
$(z_1,\ldots,z_m)$ on $\C^m$ to be $z_j=x_j+iy_j$. For $R>0$, define
$B_R$ to be the open ball of radius $R$ about 0 in~$\R^{2m}$.
\label{cs2def2}
\end{dfn}

{\it Darboux's Theorem} \cite[Th.~3.15]{McSa} says that every
symplectic manifold is locally isomorphic to $(\R^{2m},\om')$.
Our version easily follows.

\begin{thm} Let\/ $(M,\om)$ be a symplectic $2m$-manifold and\/
$x\in M$. Then there exist\/ $R>0$ and an embedding $\Up:B_R\ra
M$ with\/ $\Up(0)=x$ such that\/ $\Up^*(\om)=\om'$, where $\om'$
is the standard symplectic form on $\R^{2m}\supset B_R$. Given an
isomorphism $\up:\R^{2m}\!\ra\!T_xM$ with\/ $\up^*(\om\vert_x)\!=
\!\om'$, we can choose $\Up$ with\/~$\d\Up\vert_0\!=\!\up$.
\label{cs2thm1}
\end{thm}

Let $N$ be a real $m$-manifold. Then its tangent bundle $T^*N$ has
a {\it canonical symplectic form} $\hat\om$, defined as follows.
Let $(x_1,\ldots,x_m)$ be local coordinates on $N$. Extend them to
local coordinates $(x_1,\ldots,x_m,y_1,\ldots,y_m)$ on $T^*N$ such
that $(x_1,\ldots,y_m)$ represents the 1-form $y_1\d x_1+\cdots+y_m
\d x_m$ in $T_{(x_1,\ldots,x_m)}^*N$. Then~$\hat\om=\d x_1\w\d y_1+
\cdots+\d x_m\w\d y_m$.

Identify $N$ with the zero section in $T^*N$. Then $N$ is a
{\it Lagrangian submanifold\/} of $T^*N$. The {\it Lagrangian
Neighbourhood Theorem} \cite[Th.~3.33]{McSa} shows that any
compact Lagrangian submanifold $N$ in a symplectic manifold
looks locally like the zero section in~$T^*N$.

\begin{thm} Let\/ $(M,\om)$ be a symplectic manifold and\/
$N\subset M$ a compact Lagrangian submanifold. Then there
exists an open tubular neighbourhood\/ $U$ of the zero
section $N$ in $T^*N$, and an embedding $\Phi:U\ra M$ with\/
$\Phi\vert_N=\id:N\ra N$ and\/ $\Phi^*(\om)=\hat\om$, where
$\hat\om$ is the canonical symplectic structure on~$T^*N$.
\label{cs2thm2}
\end{thm}

We shall call $U,\Phi$ a {\it Lagrangian neighbourhood} of
$N$. Such neighbourhoods are useful for parametrizing nearby
Lagrangian submanifolds of $M$. Suppose that $\ti N$ is a
Lagrangian submanifold of $M$ which is $C^1$-close to $N$.
Then $\ti N$ lies in $\Phi(U)$, and is the image $\Phi\bigl(
\Ga(\al)\bigr)$ of the graph $\Ga(\al)$ of a unique
$C^1$-small 1-form $\al$ on~$N$.

As $\ti N$ is Lagrangian and $\Phi^*(\om)=\hat\om$ we see
that $\hat\om\vert_{\Ga(\al)}\equiv 0$. But one can easily
show that $\hat\om\vert_{\Ga(\al)}=-\pi^*(\d\al)$, where
$\pi:\Ga(\al)\ra N$ is the natural projection. Hence $\d\al=0$,
and $\al$ is a {\it closed\/ $1$-form}. This establishes a
1-1 correspondence between small closed 1-forms on $N$ and
Lagrangian submanifolds $\ti N$ close to $N$ in $M$, which
is an essential tool in proving the results of \S\ref{cs5}
and~\S\ref{cs7}.

\subsection{Special Lagrangian submanifolds in $\C^m$}
\label{cs22}

We define {\it calibrations} and {\it calibrated submanifolds},
following~\cite{HaLa}.

\begin{dfn} Let $(M,g)$ be a Riemannian manifold. An {\it oriented
tangent $k$-plane} $V$ on $M$ is a vector subspace $V$ of
some tangent space $T_xM$ to $M$ with $\dim V=k$, equipped
with an orientation. If $V$ is an oriented tangent $k$-plane
on $M$ then $g\vert_V$ is a Euclidean metric on $V$, so 
combining $g\vert_V$ with the orientation on $V$ gives a 
natural {\it volume form} $\vol_V$ on $V$, which is a 
$k$-form on~$V$.

Now let $\vp$ be a closed $k$-form on $M$. We say that
$\vp$ is a {\it calibration} on $M$ if for every oriented
$k$-plane $V$ on $M$ we have $\vp\vert_V\le \vol_V$. Here
$\vp\vert_V=\al\cdot\vol_V$ for some $\al\in\R$, and 
$\vp\vert_V\le\vol_V$ if $\al\le 1$. Let $N$ be an 
oriented submanifold of $M$ with dimension $k$. Then 
each tangent space $T_xN$ for $x\in N$ is an oriented
tangent $k$-plane. We say that $N$ is a {\it calibrated 
submanifold\/} if $\vp\vert_{T_xN}=\vol_{T_xN}$ for all~$x\in N$.
\label{cs2def3}
\end{dfn}

It is easy to show that calibrated submanifolds are automatically
{\it minimal submanifolds} \cite[Th.~II.4.2]{HaLa}. Here is the 
definition of special Lagrangian submanifolds in $\C^m$, taken
from~\cite[\S III]{HaLa}.

\begin{dfn} Let $\C^m$ have complex coordinates $(z_1,\dots,z_m)$, 
and define a metric $g'$, a real 2-form $\om'$ and a complex $m$-form 
$\Om'$ on $\C^m$ by
\e
\begin{split}
g'=\ms{\d z_1}+\cdots+\ms{\d z_m},\quad
\om'&=\ts\frac{i}{2}(\d z_1\w\d\bar z_1+\cdots+\d z_m\w\d\bar z_m),\\
\text{and}\quad\Om'&=\d z_1\w\cdots\w\d z_m.
\end{split}
\label{cs2eq2}
\e
Then $g',\om'$ are as in Definition \ref{cs2def2}, and $\Re\Om'$ and
$\Im\Om'$ are real $m$-forms on $\C^m$. Let $L$ be an oriented real
submanifold of $\C^m$ of real dimension $m$. We say that $L$ is a
{\it special Lagrangian submanifold\/} of $\C^m$, or {\it SL\/
$m$-fold}\/ for short, if $L$ is calibrated with respect to $
\Re\Om'$, in the sense of Definition~\ref{cs2def3}.
\label{cs2def4}
\end{dfn}

Harvey and Lawson \cite[Cor.~III.1.11]{HaLa} give the following
alternative characterization of special Lagrangian submanifolds:

\begin{prop} Let\/ $L$ be a real $m$-dimensional submanifold 
of\/ $\C^m$. Then $L$ admits an orientation making it into an
SL submanifold of\/ $\C^m$ if and only if\/ $\om'\vert_L\equiv 0$ 
and\/~$\Im\Om'\vert_L\equiv 0$.
\label{cs2prop1}
\end{prop}

Thus special Lagrangian submanifolds are {\it Lagrangian}
submanifolds satisfying the extra condition that
$\Im\Om'\vert_L\equiv 0$, which is how they get their name.

\subsection{Almost Calabi--Yau $m$-folds and SL $m$-folds} 
\label{cs23}

We shall define special Lagrangian submanifolds not just in
Calabi--Yau manifolds, as usual, but in the much larger
class of {\it almost Calabi--Yau manifolds}.

\begin{dfn} Let $m\ge 2$. An {\it almost Calabi--Yau $m$-fold\/}
is a quadruple $(M,J,\om,\Om)$ such that $(M,J)$ is a compact
$m$-dimensional complex manifold, $\om$ is the K\"ahler form
of a K\"ahler metric $g$ on $M$, and $\Om$ is a non-vanishing
holomorphic $(m,0)$-form on~$M$.

We call $(M,J,\om,\Om)$ a {\it Calabi--Yau $m$-fold\/} if in
addition $\om$ and $\Om$ satisfy
\e
\om^m/m!=(-1)^{m(m-1)/2}(i/2)^m\Om\w\bar\Om.
\label{cs2eq3}
\e
Then for each $x\in M$ there exists an isomorphism $T_xM\cong\C^m$
that identifies $g_x,\om_x$ and $\Om_x$ with the flat versions
$g',\om',\Om'$ on $\C^m$ in \eq{cs2eq2}. Furthermore, $g$ is
Ricci-flat and its holonomy group is a subgroup of~$\SU(m)$.
\label{cs2def5}
\end{dfn}

This is not the usual definition of a Calabi--Yau manifold, but
is essentially equivalent to it.

\begin{dfn} Let $(M,J,\om,\Om)$ be an almost Calabi--Yau $m$-fold,
and $N$ a real $m$-dimensional submanifold of $M$. We call $N$ a
{\it special Lagrangian submanifold}, or {\it SL $m$-fold\/} for
short, if $\om\vert_N\equiv\Im\Om\vert_N\equiv 0$. It easily
follows that $\Re\Om\vert_N$ is a nonvanishing $m$-form on $N$.
Thus $N$ is orientable, with a unique orientation in which
$\Re\Om\vert_N$ is positive.
\label{cs2def6}
\end{dfn}

Again, this is not the usual definition of SL $m$-fold, but is
essentially equivalent to it. In Definition \ref{cs9def3} we 
give a more general definition of SL $m$-fold involving a
{\it phase} ${\rm e}^{i\th}$. Suppose $(M,J,\om,\Om)$ is an
almost Calabi--Yau $m$-fold, with metric $g$. Let
$\psi:M\ra(0,\iy)$ be the unique smooth function such that
\e
\psi^{2m}\om^m/m!=(-1)^{m(m-1)/2}(i/2)^m\Om\w\bar\Om,
\label{cs2eq4}
\e
and define $\ti g$ to be the conformally equivalent metric $\psi^2g$
on $M$. Then $\Re\Om$ is a {\it calibration} on the Riemannian manifold
$(M,\ti g)$, and SL $m$-folds $N$ in $(M,J,\om,\Om)$ are calibrated
with respect to it, so that they are minimal with respect to~$\ti g$.

If $M$ is a Calabi--Yau $m$-fold then $\psi\equiv 1$ by \eq{cs2eq3},
so $\ti g=g$, and an $m$-submanifold $N$ in $M$ is special Lagrangian
if and only if it is calibrated w.r.t.\ $\Re\Om$ on $(M,g)$, as in
Definition \ref{cs2def4}. This recovers the usual definition of
special Lagrangian $m$-folds in Calabi--Yau $m$-folds.

\subsection{Deformations of compact SL $m$-folds} 
\label{cs24}

The {\it deformation theory} of special Lagrangian submanifolds
was studied by McLean \cite[\S 3]{McLe}, who proved the following
result in the Calabi--Yau case. The extension to the almost
Calabi--Yau case is described in~\cite[\S 9.5]{Joyc8}.

\begin{thm} Let\/ $N$ be a compact SL\/ $m$-fold in an almost
Calabi--Yau $m$-fold\/ $(M,J,\om,\Om)$. Then the moduli space
$\M_\sN$ of special Lagrangian deformations of\/ $N$ is a smooth
manifold of dimension $b^1(N)$, the first Betti number of\/~$N$.
\label{cs2thm3}
\end{thm}

Here is a sketch of the proof of Theorem \ref{cs2thm3}. Let $g$
be the K\"ahler metric on $M$, and define $\psi:M\ra(0,\iy)$ by
\eq{cs2eq4}. Applying Theorem \ref{cs2thm2} gives an open
neighbourhood $U$ of $N$ in $T^*N$ and an embedding $\Phi:U\ra M$.
Let $\pi:U\ra N$ be the natural projection. Define an $m$-form
$\be$ on $U$ by $\be=\Phi^*(\Im\Om)$. If $\al$ is a 1-form on $N$
let $\Ga(\al)$ be the graph of $\al$ in $T^*N$, and write
$C^\iy(U)\subset C^\iy(T^*N)$ for the subset of 1-forms whose
graphs lie in~$U$.

Then each submanifold $\ti N$ of $M$ which is $C^1$-close to $N$
is $\Phi\bigl(\Ga(\al)\bigr)$ for some small $\al\in C^\iy(U)$.
Here is the condition for $\smash{\ti N}$ to be special Lagrangian.

\begin{lem} In the situation above, if\/ $\al\in C^\iy(U)$
then $\ti N=\Phi\bigl(\Ga(\al)\bigr)$ is a special
Lagrangian $m$-fold in $M$ if and only if\/ $\d\al=0$
and\/~$\pi_*\bigl(\be\vert_{\Ga(\al)}\bigr)=0$.
\label{cs2lem}
\end{lem}

\begin{proof} By Definition \ref{cs2def6}, $\ti N$ is an SL $m$-fold
if and only if $\om\vert_{\ti N}\equiv\Im\Om\vert_{\ti N}\equiv 0$.
Pulling back by $\Phi$ and pushing forward by $\pi:
\Ga(\al)\ra N$, we see that $\ti N$ is special Lagrangian if
and only if $\pi_*\bigl(\hat\om\vert_{\Ga(\al)}\bigr)\equiv
\pi_*\bigl(\be\vert_{\Ga(\al)}\bigr)\equiv 0$, since
$\Phi^*(\om)=\hat\om$ and $\Phi^*(\Im\Om)=\be$. But $\pi_*\bigl(
\hat\om\vert_{\Ga(\al)}\bigr)=-\d\al$, and the lemma follows.
\end{proof}

We rewrite the condition $\pi_*\bigl(\be\vert_{\Ga(\al)}\bigr)=0$
in terms of a function~$F$.

\begin{dfn} Define $F:C^\iy(U)\ra C^\iy(N)$ by $\pi_*\bigl(
\be\vert_{\Ga(\al)}\bigr)=F(\al)\,\d V_g$, where $\d V_g$ is the
volume form of $g\vert_N$ on $N$. Then Lemma \ref{cs2lem}
shows that if $\al\in C^\iy(U)$ then $\Phi\bigl(\Ga(\al)\bigr)$
is special Lagrangian if and only if~$\d\al=F(\al)=0$.
\label{cs2def7}
\end{dfn}

In \cite[Prop.~2.10]{Joyc10} we compute the expansion of $F$
up to first order in~$\al$.

\begin{prop} This function $F$ may be written
\e
F(\al)[x]=-\d^*\bigl(\psi^m\al\bigr)+Q\bigl(x,\al(x),\na\al(x)\bigr)
\quad\text{for $x\in N$,}
\label{cs2eq5}
\e
where $Q:\bigl\{(x,y,z):x\in N$, $y\in T^*_xN\cap U$,
$z\in\ot^2T^*_xN\bigr\}\ra\R$ is smooth and\/ $Q(x,y,z)
=O(\ms{y}+\ms{z})$ for small\/~$y,z$.
\label{cs2prop2}
\end{prop}

From Definition \ref{cs2def7} and Proposition \ref{cs2prop2}
we see that the moduli space $\M_\sN$ of special Lagrangian
deformations of $N$ is locally approximately isomorphic to
the vector space of 1-forms $\al$ with $\d\al=\d^*(\psi^m\al)=0$.
But by Hodge theory, this is isomorphic to the de Rham cohomology
group $H^1(N,\R)$, and is a manifold with dimension~$b^1(N)$.

To carry out this last step rigorously requires some technical
machinery: one must work with certain {\it Banach spaces} of 
sections of $\La^kT^*N$ for $k=0,1,2$, use {\it elliptic
regularity results} to prove that the map $\al\mapsto\bigl(\d\al,
\d F\vert_0(\al)\bigr)$ is {\it surjective} upon the appropriate
Banach spaces, and then use the {\it Implicit Mapping Theorem
for Banach spaces} to show that the kernel of the map is what
we expect. This concludes our sketch of the proof of
Theorem~\ref{cs2thm3}.

Finally we extend of Theorem \ref{cs2thm3} to {\it families}
of almost Calabi--Yau $m$-folds.

\begin{dfn} Let $(M,J,\om,\Om)$ be an almost Calabi--Yau
$m$-fold. A {\it smooth family of deformations of\/}
$(M,J,\om,\Om)$ is a connected open set $\F\subset\R^d$
for $d\ge 0$ with $0\in\F$ called the {\it base space}, and
a smooth family $\bigl\{(M,J^s,\om^s,\Om^s):s\in\F\bigr\}$
of almost Calabi--Yau structures on $M$
with~$(J^0,\om^0,\Om^0)=(J,\om,\Om)$.
\label{cs2def8}
\end{dfn}

If $N$ is a compact SL $m$-fold in $(M,J,\om,\Om)$, the moduli
of deformations of $N$ in each $(M,J^s,\om^s,\Om^s)$ for $s\in\F$
make up a big moduli space~$\M_\sN^\sF$.

\begin{dfn} Let $\bigl\{(M,J^s,\om^s,\Om^s):s\in\F\bigr\}$
be a smooth family of deformations of an almost Calabi--Yau
$m$-fold $(M,J,\om,\Om)$, and $N$ be a compact SL $m$-fold
in $(M,J,\om,\Om)$. Define the {\it moduli space $\M_\sN^\sF$
of deformations of\/ $N$ in the family} $\F$ to be the set
of pairs $(s,\hat N)$ for which $s\in\F$ and $\hat N$ is a
compact SL $m$-fold in $(M,J^s,\om^s,\Om^s)$ which is
diffeomorphic to $N$ and isotopic to $N$ in $M$. Define a
{\it projection} $\pi^\sF:\M_\sN^\sF\ra\F$ by
$\pi^\sF(s,\hat N)=s$. Then $\M_\sN^\sF$ has a
natural topology, and $\pi^\sF$ is continuous.
\label{cs2def9}
\end{dfn}

The following result is proved by Marshall \cite[Th.~3.2.9]{Mars},
using similar methods to Theorem~\ref{cs2thm3}.

\begin{thm} Let\/ $\bigl\{(M,J^s,\om^s,\Om^s):s\in\F\bigr\}$ be
a smooth family of deformations of an almost Calabi--Yau $m$-fold\/
$(M,J,\om,\Om)$, with base space $\F\subset\R^d$. Suppose $N$ is a
compact SL\/ $m$-fold in $(M,J,\om,\Om)$ with\/ $[\om^s\vert_N]=0$
in $H^2(N,\R)$ and\/ $[\Im\Om^s\vert_N]=0$ in $H^m(N,\R)$ for all\/
$s\in\F$. Let\/ $\M_\sN^\sF$ be the moduli space of deformations
of\/ $N$ in $\F$, and\/ $\pi^\sF:\M_\sN^\sF\ra\F$ the natural
projection.

Then $\M_\sN^\sF$ is a smooth manifold of dimension\/ $d+b^1(N)$,
and\/ $\pi^\sF:\M_\sN^\sF\ra\F$ a smooth submersion. For small
$s\in\F$ the moduli space $\M_\sN^s=(\pi^\sF)^{-1}(s)$ of
deformations of\/ $N$ in $(M,J^s,\om^s,\Om^s)$ is a nonempty
smooth manifold of dimension $b^1(N)$, with\/~$\M_\sN^0=\M_\sN$.
\label{cs2thm4}
\end{thm}

Here a necessary condition for the existence of an SL $m$-fold
$\hat N$ isotopic to $N$ in $(M,J^s,\om^s,\Om^s)$ is that
$[\om^s\vert_N]=[\Im\Om^s\vert_N]=0$ in $H^*(N,\R)$, since
$[\om^s\vert_N]$ and $[\om^s\vert_{\smash{\hat N}}]$ are
identified under the natural isomorphism between $H^2(N,\R)$
and $H^2(\hat N,\R)$, and similarly for~$\Im\Om^s$.

The point of the theorem is that these conditions $[\om^s\vert_N]
=[\Im\Om^s\vert_N]=0$ are also {\it sufficient\/} for the
existence of $\hat N$ when $s$ is close to 0 in $\F$. That is,
the only {\it obstructions} to existence of compact SL $m$-folds
when we deform the underlying almost Calabi--Yau $m$-fold are
the obvious cohomological ones.

\section{SL cones and conical singularities}
\label{cs3}

We begin in \S\ref{cs31} with some definitions on {\it special
Lagrangian cones}. Section \ref{cs32} gives {\it examples} of
SL cones, and \S\ref{cs33} defines {\it SL\/ $m$-folds with
conical singularities}, the subject of the paper. Section
\ref{cs34} discusses {\it homology} and {\it cohomology} of
SL $m$-folds with conical singularities.

\subsection{Preliminaries on special Lagrangian cones}
\label{cs31}

We define {\it special Lagrangian cones}, and some notation.

\begin{dfn} A (singular) SL $m$-fold $C$ in $\C^m$ is called a
{\it cone} if $C=tC$ for all $t>0$, where $tC=\{t\,{\bf x}:{\bf x}
\in C\}$. Let $C$ be a closed SL cone in $\C^m$ with an isolated
singularity at 0. Then $\Si=C\cap{\cal S}^{2m-1}$ is a compact,
nonsingular $(m\!-\!1)$-submanifold of ${\cal S}^{2m-1}$, not
necessarily connected. Let $g_\sSi$ be the restriction
of $g'$ to $\Si$, where $g'$ is as in~\eq{cs2eq2}.

Set $C'=C\sm\{0\}$. Define $\iota:\Si\t(0,\iy)\ra\C^m$ by
$\iota(\si,r)=r\si$. Then $\iota$ has image $C'$. By an abuse
of notation, {\it identify} $C'$ with $\Si\t(0,\iy)$ using
$\iota$. The {\it cone metric} on $C'\cong\Si\t(0,\iy)$
is~$g'=\iota^*(g')=\d r^2+r^2g_\sSi$.

For $\al\in\R$, we say that a function $u:C'\ra\R$ is
{\it homogeneous of order} $\al$ if $u\circ t\equiv t^\al u$ for
all $t>0$. Equivalently, $u$ is homogeneous of order $\al$ if
$u(\si,r)\equiv r^\al v(\si)$ for some function~$v:\Si\ra\R$.
\label{cs3def1}
\end{dfn}

In \cite[Lem.~2.3]{Joyc9} we study {\it homogeneous harmonic
functions} on~$C'$.

\begin{lem} In the situation of Definition \ref{cs3def1},
let\/ $u(\si,r)\equiv r^\al v(\si)$ be a homogeneous function
of order $\al$ on $C'=\Si\t(0,\iy)$, for $v\in C^2(\Si)$. Then
\e
\De u(\si,r)=r^{\al-2}\bigl(\De_\sSi v-\al(\al+m-2)v\bigr),
\label{cs3eq1}
\e
where $\De$, $\De_\sSi$ are the Laplacians on $(C',g')$
and\/ $(\Si,g_\sSi)$. Hence, $u$ is harmonic on $C'$
if and only if\/ $v$ is an eigenfunction of\/ $\De_\sSi$
with eigenvalue~$\al(\al+m-2)$.
\label{cs3lem}
\end{lem}

Following \cite[Def.~2.5]{Joyc9}, we define:

\begin{dfn} In the situation of Definition \ref{cs3def1},
suppose $m>2$ and define
\e
\D_\sSi=\bigl\{\al\in\R:\text{$\al(\al+m-2)$ is
an eigenvalue of $\De_\sSi$}\bigr\}.
\label{cs3eq2}
\e
Then $\D_\sSi$ is a countable, discrete subset of
$\R$. By Lemma \ref{cs3lem}, an equivalent definition is that
$\D_\sSi$ is the set of $\al\in\R$ for which there
exists a nonzero homogeneous harmonic function $u$ of order
$\al$ on~$C'$.

Define $m_\sSi:\D_\sSi\ra\N$ by taking
$m_\sSi(\al)$ to be the multiplicity of the eigenvalue
$\al(\al+m-2)$ of $\De_\sSi$, or equivalently the
dimension of the vector space of homogeneous harmonic
functions $u$ of order $\al$ on $C'$. Define
$N_\sSi:\R\ra\Z$ by
\e
N_\sSi(\de)=
-\sum_{\!\!\!\!\al\in\D_\sSi\cap(\de,0)\!\!\!\!}m_\sSi(\al)
\;\>\text{if $\de<0$, and}\;\>
N_\sSi(\de)=
\sum_{\!\!\!\!\al\in\D_\sSi\cap[0,\de]\!\!\!\!}m_\sSi(\al)
\;\>\text{if $\de\ge 0$.}
\label{cs3eq3}
\e
Then $N_\sSi$ is monotone increasing and upper semicontinuous,
and is discontinuous exactly on $\D_\sSi$, increasing by
$m_\sSi(\al)$ at each $\al\in\D_\sSi$. As the
eigenvalues of $\De_\sSi$ are nonnegative, we see that
$\D_\sSi\cap(2-m,0)=\emptyset$ and $N_\sSi\equiv 0$
on~$(2-m,0)$.
\label{cs3def2}
\end{dfn}

We define the {\it stability index} of $C$, and {\it stable}
and {\it rigid} cones~\cite[Def.~3.6]{Joyc10}.

\begin{dfn} Let $C$ be an SL cone in $\C^m$ for $m>2$ with an
isolated singularity at 0, let $G$ be the Lie subgroup of\/
$\SU(m)$ preserving $C$, and use the notation of Definitions
\ref{cs3def1} and \ref{cs3def2}. Then \cite[eq.~(8)]{Joyc10}
shows that
\e
m_\sSi(0)=b^0(\Si),\quad
m_\sSi(1)\ge 2m \quad\text{and}\quad
m_\sSi(2)\ge m^2-1-\dim G.
\label{cs3eq4}
\e

Define the {\it stability index} $\sind(C)$ to be
\e
\sind(C)=N_\sSi(2)-b^0(\Si)-m^2-2m+1+\dim G.
\label{cs3eq5}
\e
Then $\sind(C)\ge 0$ by \eq{cs3eq4}, as $N_\sSi(2)\ge
m_\sSi(0)+m_\sSi(1)+m_\sSi(2)$ by \eq{cs3eq3}.
We call $C$ {\it stable} if~$\sind(C)=0$.

Following \cite[Def.~6.7]{Joyc9}, we call $C$ {\it rigid} if
$m_\sSi(2)=m^2-1-\dim G$. As
\begin{equation*}
\sind(C)\ge m_\sSi(2)-(m^2-1-\dim G)\ge 0,
\end{equation*}
we see that {\it if\/ $C$ is stable, then $C$ is rigid}.
\label{cs3def3}
\end{dfn}

We shall see in \S\ref{cs5} that $\sind(C)$ is the dimension
of an obstruction space to deforming an SL $m$-fold $X$ with
a conical singularity with cone $C$, and that if $C$ is
{\it stable} then the deformation theory of $X$ simplifies.
An SL cone $C$ is {\it rigid} if all infinitesimal deformations
of $C$ as an SL cone come from $\SU(m)$ rotations of $C$. This
will be useful in the Geometric Measure Theory material
of~\S\ref{cs42}.

\subsection{Examples of special Lagrangian cones}
\label{cs32}

In our first example we can compute the data of \S\ref{cs31}
very explicitly.

\begin{ex} Here is a family of special Lagrangian cones
constructed by Harvey and Lawson \cite[\S III.3.A]{HaLa}.
For $m\ge 3$, define
\e
C_{\sst\rm HL}^m=\bigl\{(z_1,\ldots,z_m)\in\C^m:i^{m+1}z_1\cdots
z_m\in[0,\iy), \quad \md{z_1}=\cdots=\md{z_m}\bigr\}.
\label{cs3eq6}
\e
Then $C_{\sst\rm HL}^m$ is a {\it special Lagrangian cone} in
$\C^m$ with an isolated singularity at 0, and $\Si_{\sst\rm HL}^m
=C_{\sst\rm HL}^m\cap{\cal S}^{2m-1}$ is an $(m\!-\!1)$-torus
$T^{m-1}$. Both $C_{\sst\rm HL}^m$ and $\Si_{\sst\rm HL}^m$ are
invariant under the $\U(1)^{m-1}$ subgroup of $\SU(m)$ acting by
\e
(z_1,\ldots,z_m)\mapsto({\rm e}^{i\th_1}z_1,\ldots,{\rm e}^{i\th_m}
z_m) \quad\text{for $\th_j\in\R$ with $\th_1+\cdots+\th_m=0$.}
\label{cs3eq7}
\e
In fact $\pm\,C_{\sst\rm HL}^m$ for $m$ odd, and $C_{\sst\rm
HL}^m,iC_{\sst\rm HL}^m$ for $m$ even, are the unique SL cones in
$\C^m$ invariant under \eq{cs3eq7}, which is how Harvey and Lawson
constructed them.

The metric on $\Si_{\sst\rm HL}^m\cong T^{m-1}$ is flat, so it is
easy to compute the eigenvalues of $\De_{\smash{\sst\Si_{\rm HL}^m}}$.
This was done by Marshall \cite[\S 6.3.4]{Mars}. There is a 1-1
correspondence between $(n_1,\ldots,n_{m-1})\in\Z^{m-1}$ and
eigenvectors of $\De_{\smash{\sst\Si_{\rm HL}^m}}$ with eigenvalue
\e
m\sum_{i=1}^{m-1}n_i^2-\sum_{i,j=1}^{m-1}n_in_j.
\label{cs3eq8}
\e

Using \eq{cs3eq8} and a computer we can find the eigenvalues
of $\De_{\smash{\sst\Si_{\rm HL}^m}}$ and their multiplicities.
The Lie subgroup $G_{\sst\rm HL}^m$ of $\SU(m)$ preserving
$C_{\sst\rm HL}^m$ has identity component the $\U(1)^{m-1}$ of
\eq{cs3eq7}, so that $\dim G_{\sst\rm HL}^m=m-1$. Thus we can
calculate $\D_{\smash{\sst\Si_{\rm HL}^m}}$, $m_{\smash{\sst
\Si_{\rm HL}^m}}$, $N_{\smash{\sst\Si_{\rm HL}^m}}$, and the
stability index $\sind(C_{\sst\rm HL}^m)$. This was done by
Marshall \cite[Table 6.1]{Mars} and the author
\cite[\S 3.2]{Joyc10}. Table \ref{cs3table} gives the data
$m,N_{\smash{\sst\Si_{\rm HL}^m}}(2),m_{\smash{\sst\Si_{\rm
HL}^m}}(2)$ and $\sind(C_{\sst\rm HL}^m)$ for~$3\le m\le 12$.

\begin{table}[htb]
\center{
\begin{tabular}{|l|r|r|r|r|r|r|r|r|r|r|}\hline
$m\vphantom{\bigr(^l_j}$
&  3 &  4 &  5 &  6 &  7 &  8 &  9 &  10 & 11 & 12 \\
\hline
$N_{\smash{\sst\Si_{\rm HL}^m}}(2)\vphantom{\bigr(^l_j}$
& 13 & 27 & 51 & 93 & 169& 311& 331& 201 & 243& 289\\
\hline
$m_{\smash{\sst\Si_{\rm HL}^m}}(2)\vphantom{\bigr(^l_j}$
&  6 & 12 & 20 & 30 & 42 & 126& 240&  90 & 110& 132\\
\hline
$\sind(C_{\sst\rm HL}^m)\vphantom{\bigr(^l_j}$
&  0 &  6 & 20 & 50 & 112& 238& 240&  90 & 110& 132\\
\hline
\end{tabular}
}
\caption{Data for $\U(1)^{m-1}$-invariant SL cones
$C_{\sst\rm HL}^m$ in $\C^m$}
\label{cs3table}
\end{table}

One can also prove that
\e
N_{\smash{\sst\Si_{\rm HL}^m}}(2)=2m^2+1\;\>\text{and}\;\>
m_{\smash{\sst\Si_{\rm HL}^m}}(2)=\sind(C_{\sst\rm HL}^m)=m^2-m
\;\>\text{for $m\ge 10$.}
\label{cs3eq9}
\e
As $C_{\sst\rm HL}^m$ is {\it stable} when $\sind(C_{\sst\rm
HL}^m)=0$ we see from Table \ref{cs3table} and \eq{cs3eq9}
that $C_{\sst\rm HL}^3$ is a {\it stable} cone in $\C^3$, but
$C_{\sst\rm HL}^m$ is {\it unstable} for $m\ge 4$. Also
$C_{\sst\rm HL}^m$ is {\it rigid\/} when $m_{\smash{\sst\Si_{
\rm HL}^m}}(2)=m^2-m$, as $\dim G_{\sst\rm HL}^m=m-1$. Thus
$C_{\sst\rm HL}^m$ is {\it rigid\/} if and only if $m\ne 8,9$,
by Table \ref{cs3table} and~\eq{cs3eq9}.
\label{cs3ex1}
\end{ex}

Here is an example taken from \cite[Ex.~9.4]{Joyc2}, chosen
as it is easy to write down.

\begin{ex} Let $a_1,\ldots,a_m\in\Z$ with $a_1+\cdots+a_m=0$
and highest common factor 1, such that $a_1,\ldots,a_k>0$
and $a_{k+1},\ldots,a_m<0$ for $0<k<m$. Define
\e
\begin{split}
L^{a_1,\ldots,a_m}_0=\bigl\{
\bigl(i{\rm e}^{ia_1\th}x_1&,{\rm e}^{ia_2\th}x_2,\ldots,
{\rm e}^{ia_m\th}x_m\bigr):
\th\in[0,2\pi),\\ 
&x_1,\ldots,x_m\in\R,\qquad 
a_1x_1^2+\cdots+a_mx_m^2=0\bigr\}.
\end{split}
\label{cs3eq10}
\e
Then $L^{a_1,\ldots,a_m}_0$ is an {\it immersed SL cone} in
$\C^m$, with an isolated singularity at~0.

Define $C^{a_1,\ldots,a_m}=\bigl\{(x_1,\ldots,x_m)\in\R^m:
a_1x_1^2+\cdots+a_mx_m^2=0\bigr\}$. Then $C^{a_1,\ldots,a_m}$
is a quadric cone on ${\cal S}^{k-1}\t S^{m-k-1}$ in $\R^m$,
and $L^{a_1,\ldots,a_m}_0$ is the image of an immersion $\Phi:
C^{a_1,\ldots,a_m}\t{\cal S}^1\ra\C^m$, which is generically
2:1. Therefore $L^{a_1,\ldots,a_m}_0$ is an immersed SL cone
on~$({\cal S}^{k-1}\t{\cal S}^{m-k-1}\t{\cal S}^1)/\Z_2$.
\label{cs3ex2}
\end{ex}

Further examples of SL cones are constructed by Harvey and
Lawson \cite[\S III.3]{HaLa}, Haskins \cite{Hask}, the
author \cite{Joyc2,Joyc3}, and others. Special Lagrangian
cones in $\C^3$ are a special case, which may be treated
using the theory of {\it integrable systems}. In principle
this should yield a {\it classification} of all SL cones on
$T^2$ in $\C^3$. For more information see McIntosh \cite{McIn}
or the author~\cite{Joyc7}.

\subsection{Special Lagrangian $m$-folds with conical singularities}
\label{cs33}

Now we can define {\it conical singularities} of SL $m$-folds,
following~\cite[Def.~3.6]{Joyc9}.

\begin{dfn} Let $(M,J,\om,\Om)$ be an almost Calabi--Yau $m$-fold
for $m>2$, and define $\psi:M\ra(0,\iy)$ as in \eq{cs2eq4}. Suppose
$X$ is a compact singular SL $m$-fold in $M$ with singularities at
distinct points $x_1,\ldots,x_n\in X$, and no other singularities.

Fix isomorphisms $\up_i:\C^m\ra T_{x_i}M$ for $i=1,\ldots,n$
such that $\up_i^*(\om)=\om'$ and $\up_i^*(\Om)=\psi(x_i)^m\Om'$,
where $\om',\Om'$ are as in \eq{cs2eq2}. Let $C_1,\ldots,C_n$ be SL
cones in $\C^m$ with isolated singularities at 0. For $i=1,\ldots,n$
let $\Si_i=C_i\cap{\cal S}^{2m-1}$, and let $\mu_i\in(2,3)$ with
\e
(2,\mu_i]\cap\D_\sSii=\emptyset,
\quad\text{where $\D_\sSii$ is defined in \eq{cs3eq2}.}
\label{cs3eq11}
\e
Then we say that
$X$ has a {\it conical singularity} or {\it conical singular
point} at $x_i$, with {\it rate} $\mu_i$ and {\it cone} $C_i$
for $i=1,\ldots,n$, if the following holds.

By Theorem \ref{cs2thm1} there exist embeddings $\Up_i:B_R\ra M$
for $i=1,\ldots,n$ satisfying $\Up_i(0)=x_i$, $\d\Up_i\vert_0=\up_i$
and $\Up_i^*(\om)=\om'$, where $B_R$ is the open ball of radius $R$
about 0 in $\C^m$ for some small $R>0$. Define $\iota_i:\Si_i\t(0,R)
\ra B_R$ by $\iota_i(\si,r)=r\si$ for~$i=1,\ldots,n$.

Define $X'=X\sm\{x_1,\ldots,x_n\}$. Then there should exist a
compact subset $K\subset X'$ such that $X'\sm K$ is a union of
open sets $S_1,\ldots,S_n$ with $S_i\subset\Up_i(B_R)$, whose
closures $\bar S_1,\ldots,\bar S_n$ are disjoint in $X$. For
$i=1,\ldots,n$ and some $R'\in(0,R]$ there should exist a smooth
$\phi_i:\Si_i\t(0,R')\ra B_R$ such that $\Up_i\circ\phi_i:\Si_i
\t(0,R')\ra M$ is a diffeomorphism $\Si_i\t(0,R')\ra S_i$, and
\e
\bmd{\na^k(\phi_i-\iota_i)}=O(r^{\mu_i-1-k})
\quad\text{as $r\ra 0$ for $k=0,1$.}
\label{cs3eq12}
\e
Here $\na$ is the Levi-Civita connection of the cone metric
$\iota_i^*(g')$ on $\Si_i\t(0,R')$, $\md{\,.\,}$ is computed
using $\iota_i^*(g')$. If the cones $C_1,\ldots,C_n$ are
{\it stable} in the sense of Definition \ref{cs3def3}, then
we say that $X$ has {\it stable conical singularities}.
\label{cs3def4}
\end{dfn}

We will see in Theorem \ref{cs4thm4} that if \eq{cs3eq12} holds
for $k=0,1$ and some $\mu_i$ satisfying \eq{cs3eq11}, then we
can choose a natural $\phi_i$ for which \eq{cs3eq12} holds for
{\it all\/} $k\ge 0$, and for {\it all\/} rates $\mu_i$
satisfying \eq{cs3eq11}. Thus the number of derivatives
required in \eq{cs3eq12} and the choice of $\mu_i$ both make
little difference. We choose $k=0,1$ in \eq{cs3eq12}, and some
$\mu_i$ in \eq{cs3eq11}, to make the definition as weak
as possible.

We suppose $m>2$ for two reasons. Firstly, the only SL cones
in $\C^2$ are finite unions of SL planes $\R^2$ in $\C^2$
intersecting only at 0. Thus any SL 2-fold with conical
singularities is actually {\it nonsingular} as an immersed
2-fold, so there is really no point in studying them.
Secondly, $m=2$ is a special case in the analysis of
\cite[\S 2]{Joyc9}, and it is simpler to exclude it.
Therefore we will suppose $m>2$ throughout the paper.

Here are the reasons for the conditions on $\mu_i$ in
Definition~\ref{cs3def4}:
\begin{itemize}
\setlength{\parsep}{0pt}
\setlength{\itemsep}{0pt}
\item We need $\mu_i>2$, or else \eq{cs3eq12} does not force
$X$ to approach $C_i$ near~$x_i$.
\item The definition involves a choice of $\Up_i:B_R\ra M$.
If we replace $\Up_i$ by a different choice $\ti\Up_i$ then
we should replace $\phi_i$ by $\ti\phi_i=(\ti\Up_i^{-1}\circ
\Up_i)\circ\phi_i$ near 0 in $B_R$. Calculation shows that
as $\Up_i,\ti\Up_i$ agree up to second order, we
have~$\bmd{\na^k(\ti\phi_i-\phi_i)}=O(r^{2-k})$.

Therefore we choose $\mu_i<3$ so that these $O(r^{2-k})$
terms are absorbed into the $O(r^{\mu_i-1-k})$ in \eq{cs3eq12}.
This makes the definition independent of the choice of
$\Up_i$, which it would not be if~$\mu_i>3$.

\item Condition \eq{cs3eq11} is needed to prove the regularity
result Theorem \ref{cs4thm4}, and also to reduce to a minimum
the obstructions to deforming compact SL $m$-folds with
conical singularities studied in~\S\ref{cs5}.
\end{itemize}

\subsection{Homology and cohomology}
\label{cs34}

Next we discuss {\it homology} and {\it cohomology} of SL $m$-folds
with conical singularities, following \cite[\S 2.4]{Joyc9}. For a
general reference, see for instance Bredon \cite{Bred}. When $Y$
is a manifold, write $H^k(Y,\R)$ for the $k^{\rm th}$ {\it de Rham
cohomology group} and $H^k_{\rm cs}(Y,\R)$ for the $k^{\rm th}$
{\it compactly-supported de Rham cohomology group} of $Y$. If $Y$
is compact then $H^k(Y,\R)=H^k_{\rm cs}(Y,\R)$. The {\it Betti
numbers} of $Y$ are $b^k(Y)=\dim H^k(Y,\R)$ and~$b^k_{\rm cs}(Y)
=\dim H^k_{\rm cs}(Y,\R)$.

Let $Y$ be a topological space, and $Z\subset Y$ a subspace.
Write $H_k(Y,\R)$ for the $k^{\rm th}$ {\it real singular
homology group} of $Y$, and $H_k(Y;Z,\R)$ for the $k^{\rm th}$
{\it real singular relative homology group} of $(Y;Z)$. When
$Y$ is a manifold and $Z$ a submanifold we define $H_k(Y,\R)$
and $H_k(Y;Z,\R)$ using {\it smooth\/} simplices, as in
\cite[\S V.5]{Bred}. Then the pairing between (singular)
homology and (de Rham) cohomology is defined at the chain
level by integrating $k$-forms over $k$-simplices.

Suppose $X$ is a compact SL $m$-fold in $M$ with conical
singularities $x_1,\ldots,x_n$ and cones $C_1,\ldots,C_n$, and
set $X'=X\sm\{x_1,\ldots,x_n\}$ and $\Si_i=C_i\cap{\cal S}^{2m-1}$,
as in \S\ref{cs33}. Then $X'$ is the interior of a compact manifold
$\bar X'$ with boundary $\coprod_{i=1}^n\Si_i$. Using this we show
in \cite[\S 2.4]{Joyc9} that there is a natural long exact sequence
\e
\cdots\ra
H^k_{\rm cs}(X',\R)\ra H^k(X',\R)\ra\bigoplus_{i=1}^n
H^k(\Si_i,\R)\ra H^{k+1}_{\rm cs}(X',\R)\ra\cdots,
\label{cs3eq13}
\e
and natural isomorphisms
\begin{gather}
H_k\bigl(X;\{x_1,\ldots,x_n\},\R\bigr)^*\!\cong\!
H^k_{\rm cs}(X',\R)\!\cong\!H_{m-k}(X',\R)\!\cong\!H^{m-k}(X',\R)^*
\label{cs3eq14}\\
\text{and}\quad
H^k_{\rm cs}(X',\R)\cong H_k(X,\R)^*
\quad\text{for all $k>1$.}
\label{cs3eq15}
\end{gather}
The inclusion $\iota:X\ra M$ induces homomorphisms~$\iota_*:
H_k(X,\R)\ra H_k(M,\R)$.

\section{The asymptotic behaviour of $X$ near $x_i$}
\label{cs4}

We now review the work of \cite{Joyc9} on the {\it regularity}
of SL $m$-folds with conical singularities. Let $M$ be an
almost Calabi--Yau $m$-fold and $X$ an SL $m$-fold in $M$
with conical singularities at $x_1,\ldots,x_n$, with
identifications $\up_i$ and cones $C_i$. We study how
quickly $X$ converges to the cone $\up(C_i)$ in
$T_{\smash{x_i}}M$ near~$x_i$.

We start in \S\ref{cs41} by writing $X$ in a special
coordinate system near $x_i$, as the graph of an exact
1-form $\eta_i=\d A_i$ on $C_i'=C_i\sm\{0\}$. The special
Lagrangian condition reduces to a {\it nonlinear elliptic
p.d.e.}\ on the function $A_i$. In \S\ref{cs42} we explain
how {\it elliptic regularity} of this p.d.e.\ implies that
$A_i$ and its derivatives decay quickly near~$x_i$.

\subsection{Lagrangian Neighbourhood Theorems}
\label{cs41}

In \cite[Th.~4.3]{Joyc9} we extend the {\it Lagrangian Neighbourhood
Theorem}, Theorem \ref{cs2thm2}, to special Lagrangian cones.

\begin{thm} Let\/ $C$ be an SL cone in $\C^m$ with isolated
singularity at\/ $0$, and set\/ $\Si=C\cap{\cal S}^{2m-1}$.
Define $\iota:\Si\t(0,\iy)\ra\C^m$ by $\iota(\si,r)=r\si$,
with image $C\sm\{0\}$. For $\si\in\Si$, $\tau\in T_\si^*\Si$,
$r\in(0,\iy)$ and\/ $u\in\R$, let\/ $(\si,r,\tau,u)$ represent
the point\/ $\tau+u\,\d r$ in $T^*_{\smash{(\si,r)}}\bigl(\Si\!
\t\!(0,\iy)\bigr)$. Identify $\Si\t(0,\iy)$ with the zero section
$\tau\!=\!u\!=\!0$ in $T^*\bigl(\Si\t(0,\iy)\bigr)$. Define an
action of\/ $(0,\iy)$ on $T^*\bigl(\Si\!\t\!(0,\iy)\bigr)$ by
\e
t:(\si,r,\tau,u)\longmapsto (\si,tr,t^2\tau,tu)
\quad\text{for $t\in(0,\iy),$}
\label{cs4eq1}
\e
so that\/ $t^*(\hat\om)\!=\!t^2\hat\om$, for $\hat\om$ the
canonical symplectic structure on~$T^*\bigl(\Si\!\t\!(0,\iy)\bigr)$.

Then there exists an open neighbourhood\/ $U_\sC$ of\/
$\Si\t(0,\iy)$ in $T^*\bigl(\Si\t(0,\iy)\bigr)$ invariant under
\eq{cs4eq1} given by
\e
U_\sC=\bigl\{(\si,r,\tau,u)\in T^*\bigl(\Si\t(0,\iy)\bigr):
\bmd{(\tau,u)}<2\ze r\bigr\}\quad\text{for some $\ze>0,$}
\label{cs4eq2}
\e
where $\md{\,.\,}$ is calculated using the cone metric $\iota^*(g')$
on $\Si\t(0,\iy)$, and an embedding $\Phi_\sC:U_\sC\ra\C^m$
with\/ $\Phi_\sC\vert_{\Si\t(0,\iy)}=\iota$, $\Phi_{\sst
C}^*(\om')=\hat\om$ and\/ $\Phi_\sC\circ t=t\,\Phi_\sC$
for all\/ $t>0$, where $t$ acts on $U_\sC$ as in \eq{cs4eq1}
and on $\C^m$ by multiplication.
\label{cs4thm1}
\end{thm}

These $U_\sC,\Phi_\sC$ are a {\it Lagrangian neighbourhood\/}
of $C'$ in $\C^m$ which is {\it equivariant under the action
of dilations}. Effectively they are a special coordinate system
on $\C^m$ near $C'$, in which $\om'$ assumes a simple form. In
\cite[Th.~4.4]{Joyc9} we use $U_\sCi,\Phi_\sCi$ to construct a
particular choice of $\phi_i$ in Definition~\ref{cs3def4}.

\begin{thm} Let\/ $(M,J,\om,\Om)$, $\psi,X,n,x_i,\up_i,C_i,\Si_i,
\mu_i,R,\Up_i$ and\/ $\iota_i$ be as in Definition \ref{cs3def4}.
Theorem \ref{cs4thm1} gives $\ze>0$, neighbourhoods $U_\sCi$
of\/ $\Si_i\t(0,\iy)$ in $T^*\bigl(\Si_i\t(0,\iy)\bigr)$ and
embeddings $\Phi_\sCi:U_\sCi\ra\C^m$ for~$i=1,\ldots,n$.

Then for sufficiently small\/ $R'\in(0,R]$ there exist unique
closed\/ $1$-forms $\eta_i$ on $\Si_i\t(0,R')$ for $i=1,\ldots,n$
written $\eta_i(\si,r)=\eta_i^1(\si,r)+\eta_i^2(\si,r)\d r$ for
$\eta_i^1(\si,r)\in T_\si^*\Si_i$ and\/ $\eta_i^2(\si,r)\in\R$,
and satisfying $\md{\eta_i(\si,r)}<\ze r$ and
\e
\bmd{\na^k\eta_i}=O(r^{\mu_i-1-k})
\quad\text{as $r\ra 0$ for $k=0,1,$}
\label{cs4eq3}
\e
computing $\na,\md{\,.\,}$ using the cone metric $\iota_i^*(g')$,
such that the following holds.

Define $\phi_i:\Si_i\t(0,R')\ra B_R$ by $\phi_i(\si,r)=\Phi_{\sst
C_i}\bigl(\si,r,\eta_i^1(\si,r),\eta_i^2(\si,r)\bigr)$. Then
$\Up_i\circ\phi_i:\Si_i\t(0,R')\ra M$ is a diffeomorphism
$\Si_i\t(0,R')\ra S_i$ for open sets $S_1,\ldots,S_n$ in $X'$
with\/ $\bar S_1,\ldots,\bar S_n$ disjoint, and\/ $K=X'\sm(S_1
\cup\cdots\cup S_n)$ is compact. Also $\phi_i$ satisfies
\eq{cs3eq12}, so that\/ $R',\phi_i,S_i,K$ satisfy
Definition~\ref{cs3def4}.
\label{cs4thm2}
\end{thm}

We explained in \S\ref{cs21} that in a Lagrangian neighbourhood
$U,\Phi$ of a Lagrangian $m$-fold $N$ gives a 1-1 correspondence
between nearby Lagrangian $m$-folds $\ti N$ and closed 1-forms
on $N$. Theorem \ref{cs4thm2} uses this correspondence for the
Lagrangian neighbourhoods $U_\sCi,\Phi_\sCi$ of Theorem
\ref{cs4thm1}. This is why the 1-forms $\eta_i$ are closed.
We can extend Theorem \ref{cs2thm2} to SL $m$-folds with
conical singularities \cite[Th.~4.6]{Joyc9}, in a way
compatible with Theorems \ref{cs4thm1} and~\ref{cs4thm2}.

\begin{thm} Suppose $(M,J,\om,\Om)$ is an almost Calabi--Yau
$m$-fold and\/ $X$ a compact SL\/ $m$-fold in $M$ with conical
singularities at\/ $x_1,\ldots,x_n$. Let the notation $\psi,\up_i,
C_i,\Si_i,\mu_i,R,\Up_i$ and\/ $\iota_i$ be as in Definition
\ref{cs3def4}, and let\/ $\ze,U_\sCi,\allowbreak
\Phi_\sCi,\allowbreak R',\allowbreak \eta_i,\allowbreak
\eta_i^1,\eta_i^2,\phi_i,S_i$ and\/ $K$ be as in Theorem~\ref{cs4thm2}.

Then making $R'$ smaller if necessary, there exists an open tubular
neighbourhood\/ $U_\sXp\subset T^*X'$ of the zero section
$X'$ in $T^*X'$, such that under $\d(\Up_i\circ\phi_i):T^*\bigl(
\Si_i\t(0,R')\bigr)\ra T^*X'$ for $i=1,\ldots,n$ we have
\e
\bigl(\d(\Up_i\circ\phi_i)\bigr)^*(U_\sXp)=\bigl\{(\si,r,\tau,u)
\in T^*\bigl(\Si_i\t(0,R')\bigr):\bmd{(\tau,u)}<\ze r\bigr\},
\label{cs4eq4}
\e
and there exists an embedding $\Phi_\sXp:U_\sXp\ra M$ with\/
$\Phi_\sXp\vert_{X'}=\id:X'\ra X'$ and\/ $\Phi_\sXp^*(\om)=\hat\om$,
where $\hat\om$ is the canonical symplectic structure on $T^*X'$,
such that
\e
\Phi_\sXp\circ\d(\Up_i\circ\phi_i)(\si,r,\tau,u)\equiv\Up_i\circ
\Phi_\sCi\bigl(\si,r,\tau+\eta_i^1(\si,r),u+\eta_i^2(\si,r)\bigr)
\label{cs4eq5}
\e
for all\/ $i=1,\ldots,n$ and\/ $(\si,r,\tau,u)\in T^*\bigl(\Si_i\t(0,R')
\bigr)$ with\/ $\bmd{(\tau,u)}<\ze r$. Here $\md{(\tau,u)}$ is computed
using the cone metric $\iota_i^*(g')$ on~$\Si_i\t(0,R')$.
\label{cs4thm3}
\end{thm}

This is an essential tool in the deformation theory of \S\ref{cs5}
and desingularization results of \S\ref{cs7}, as it gives a special
coordinate system on $M$ near $X'$ in which $\om$ assumes a simple
form. In these coordinate, deformations or desingularizations of $X$
become {\it graphs of closed\/ $1$-forms} on $X'$ away from $x_i$,
as in~\S\ref{cs21}.

\subsection{Regularity of $X$ near $x_i$}
\label{cs42}

The results of \S\ref{cs41} used only the fact that $X'$
is {\it Lagrangian} in $(M,\om)$. Our next theorems make
essential use of the {\it special\/} Lagrangian condition.
In \cite[\S 5]{Joyc9} we study the asymptotic behaviour of
the maps $\phi_i$ of Theorem \ref{cs4thm2}. Combining
\cite[Th.~5.1]{Joyc9}, \cite[Lem.~4.5]{Joyc9} and
\cite[Th.~5.5]{Joyc9} proves:

\begin{thm} In the situation of Theorem \ref{cs4thm2} we have
$\eta_i=\d A_i$ for $i=1,\ldots,n$, where $A_i:\Si_i\t(0,R')\ra\R$
is given by $A_i(\si,r)=\int_0^r\eta_i^2(\si,s)\d s$. Suppose
$\mu_i'\in(2,3)$ with\/ $(2,\mu_i']\cap\D_\sSii=\emptyset$ for
$i=1,\ldots,n$. Then
\e
\begin{gathered}
\bmd{\na^k(\phi_i-\iota_i)}=O(r^{\mu_i'-1-k}),\quad
\bmd{\na^k\eta_i}=O(r^{\mu_i'-1-k})\quad\text{and}\\
\bmd{\na^kA_i}=O(r^{\mu_i'-k})
\quad\text{as $r\ra 0$ for all\/ $k\ge 0$ and\/ $i=1,\ldots,n$.}
\end{gathered}
\label{cs4eq6}
\e

Hence $X$ has conical singularities at $x_i$ with cone $C_i$
and rate $\mu_i'$, for all possible rates $\mu_i'$ allowed by
Definition \ref{cs3def4}. Therefore, the definition of
conical singularities is essentially independent of the
choice of rate~$\mu_i$.
\label{cs4thm4}
\end{thm}

Theorem \ref{cs4thm4} in effect {\it strengthens} the definition of
SL $m$-folds with conical singularities, Definition \ref{cs3def4},
as it shows that \eq{cs3eq12} actually implies the much stronger
condition \eq{cs4eq6} on all derivatives.

The proof works by treating $X'$ near $x_i$ as a deformation
of the SL cone $C_i$ in $\C^m$. Thus we can apply Proposition
\ref{cs2prop2} with $N$ replaced by $\Si_i\t(0,R')$ and 
$\al=\eta_i=\d A_i$, and we find that $A_i$ satisfies the
second-order nonlinear~p.d.e.
\e
\d^*\bigl(\psi^m\d A_i\bigr)(\si,r)=
Q\bigl(\si,r,\d A_i(\si,r),\na^2A_i(\si,r)\bigr)
\label{cs4eq7}
\e
for $(\si,r)\in\Si_i\t(0,R')$, where $Q$ is a smooth nonlinear
function.

When $r$ is small the $Q$ term in \eq{cs4eq7} is also small
and \eq{cs4eq7} approximates $\De_iA_i=0$, where $\De_i$ is
the Laplacian on the cone $C_i$. Therefore \eq{cs4eq7} is
{\it elliptic} for small $r$. Using known results on the
regularity of solutions of nonlinear second-order elliptic
p.d.e.s, and a theory of analysis on weighted Sobolev spaces
on manifolds with cylindrical ends developed by Lockhart and
McOwen \cite{Lock}, we can then prove~\eq{cs4eq6}.

Our next result \cite[Th.~6.8]{Joyc9} is an application of
{\it Geometric Measure Theory}. For an introduction to the
subject, see Morgan \cite{Morg}. Geometric Measure Theory
studies measure-theoretic generalizations of submanifolds
called {\it integral currents}, which may be very singular,
and is particularly powerful for {\it minimal\/} submanifolds.
As from \S\ref{cs2} SL $m$-folds are minimal submanifolds
w.r.t.\ an appropriate metric, many major results of Geometric
Measure Theory immediately apply to {\it special Lagrangian
integral currents}, a very general class of singular SL
$m$-folds with strong compactness properties.

\begin{thm} Let\/ $(M,J,\om,\Om)$ be an almost Calabi--Yau
$m$-fold and define $\psi:M\ra(0,\iy)$ as in \eq{cs2eq4}.
Let\/ $x\in M$ and fix an isomorphism $\up:\C^m\ra T_xM$
with\/ $\up^*(\om)=\om'$ and\/ $\up^*(\Om)=\psi(x)^m\Om'$,
where $\om',\Om'$ are as in~\eq{cs2eq2}.

Suppose that\/ $T$ is a special Lagrangian integral current
in $M$ with\/ $x\in T^\circ$, and that\/ $\up_*(C)$ is a
multiplicity $1$ tangent cone to $T$ at\/ $x$, where $C$ is
a rigid special Lagrangian cone in $\C^m$, in the sense of
Definition \ref{cs3def3}. Then $T$ has a conical singularity
at\/ $x$, in the sense of Definition~\ref{cs3def4}.
\label{cs4thm5}
\end{thm}

This is a {\it weakening} of Definition \ref{cs3def4} for
{\it rigid\/} cones $C$. Theorem \ref{cs4thm5} also holds
for the larger class of {\it Jacobi integrable} SL cones
$C$, defined in~\cite[Def.~6.7]{Joyc9}.

Basically, Theorem \ref{cs4thm5} shows that if a singular
SL $m$-fold $T$ in $M$ is locally modelled on a rigid SL
cone $C$ in only a very weak sense, then it necessarily
satisfies Definition \ref{cs3def4}. One moral of Theorems
\ref{cs4thm4} and \ref{cs4thm5} is that, at least for
rigid SL cones $C$, more-or-less {\it any} sensible
definition of SL $m$-folds with conical singularities is
equivalent to Definition~\ref{cs3def4}.

Theorem \ref{cs4thm5} is proved by applying regularity
results of Allard and Almgren, and Adams and Simon, mildly
adapted to the special Lagrangian situation, which roughly
say that if a tangent cone $C_i$ to $X$ at $x_i$ has an
isolated singularity at 0, is multiplicity 1, and rigid,
then $X$ has a parametrization $\phi_i$ near $x_i$ which
satisfies \eq{cs3eq12} for some $\mu_i>2$. It then quickly
follows that $X$ has a conical singularity at $x_i$, in
the sense of Definition~\ref{cs3def4}.

As discussed in \cite[\S 6.3]{Joyc9}, one can use other
results from Geometric Measure Theory to argue that for
tangent cones $C$ which are not Jacobi integrable, Definition
\ref{cs3def4} may be {\it too strong}, in that there could
exist examples of singular SL $m$-folds with tangent cone
$C$ which are not covered by Definition \ref{cs3def4}, as
the decay conditions \eq{cs3eq12} are too strict. 

\section{Moduli of SL $m$-folds with conical singularities}
\label{cs5}

Next we review the work of \cite{Joyc10} on {\it deformation
theory} for compact SL $m$-folds with conical singularities.
Following \cite[Def.~5.4]{Joyc10}, we define the space $\M_\sX$
of compact SL $m$-folds $\hat X$ in $M$ with conical singularities
deforming a fixed SL $m$-fold $X$ with conical singularities.

\begin{dfn} Let $(M,J,\om,\Om)$ be an almost Calabi--Yau
$m$-fold and $X$ a compact SL $m$-fold in $M$ with conical
singularities at $x_1,\ldots,x_n$ with identifications
$\up_i:\C^m\ra T_{x_i}M$ and cones $C_1,\ldots,C_n$. Define
the {\it moduli space} $\M_\sX$ {\it of deformations of\/}
$X$ to be the set of $\hat X$ such that
\begin{itemize}
\setlength{\parsep}{0pt}
\setlength{\itemsep}{0pt}
\item[(i)] $\hat X$ is a compact SL $m$-fold in $M$ with
conical singularities at $\hat x_1,\ldots,\hat x_n$ with
cones $C_1,\ldots,C_n$, for some $\hat x_i$ and
identifications~$\hat\up_i:\C^m\ra T_{\smash{\hat x_i}}M$.
\item[(ii)] There exists a homeomorphism $\hat\iota:X\ra\hat X$
with $\hat\iota(x_i)=\hat x_i$ for $i=1,\ldots,n$ such that
$\hat\iota\vert_{X'}:X'\ra\hat X'$ is a diffeomorphism and
$\hat\iota$ and $\iota$ are isotopic as continuous maps
$X\ra M$, where $\iota:X\ra M$ is the inclusion.
\end{itemize}

In \cite[Def.~5.6]{Joyc10} we define a {\it topology} on
$\M_\sX$, and explain why it is the natural choice. We will
not repeat the complicated definition here; readers are
referred to \cite[\S 5]{Joyc10} for details.
\label{cs5def1}
\end{dfn}

In \cite[Th.~6.10]{Joyc10} we describe $\M_\sX$ near $X$, in terms
of a smooth map $\Phi$ between the {\it infinitesimal deformation
space} $\I_\sXp$ and the {\it obstruction space}~$\O_\sXp$.

\begin{thm} Suppose $(M,J,\om,\Om)$ is an almost Calabi--Yau
$m$-fold and\/ $X$ a compact SL\/ $m$-fold in $M$ with conical
singularities at\/ $x_1,\ldots,x_n$ and cones $C_1,\ldots,C_n$.
Let\/ $\M_\sX$ be the moduli space of deformations
of\/ $X$ as an SL\/ $m$-fold with conical singularities in $M$,
as in Definition \ref{cs5def1}. Set\/~$X'=X\sm\{x_1,\ldots,x_n\}$.

Then there exist natural finite-dimensional vector spaces
$\I_\sXp$, $\O_\sXp$ such that\/ $\I_\sXp$ is isomorphic to
the image of\/ $H^1_{\rm cs}(X',\R)$ in $H^1(X',\R)$ and\/ 
$\dim\O_\sXp=\sum_{i=1}^n\sind(C_i)$, where $\sind(C_i)$ is
the stability index of Definition \ref{cs3def3}. There exists
an open neighbourhood\/ $U$ of\/ $0$ in $\I_\sXp$, a smooth
map $\Phi:U\ra\O_\sXp$ with\/ $\Phi(0)=0$, and a map
$\Xi:\{u\in U:\Phi(u)=0\}\ra\M_\sX$ with\/ $\Xi(0)=X$ which is
a homeomorphism with an open neighbourhood of\/ $X$ in~$\M_\sX$.
\label{cs5thm1}
\end{thm}

Here is a sketch of the proof. For simplicity, consider
first the subset of $\hat X\in\M_\sX$ which have the
same singular points $x_1,\ldots,x_n$ and identifications
$\up_1,\ldots,\up_n$ as $X$. If $\hat X$ is $C^1$ close
to $X$ in an appropriate sense then $\hat X'=\Phi_\sXp
\bigl(\Ga(\al)\bigr)$, where $U_\sXp,\Phi_\sXp$ is the
Lagrangian neighbourhood map of Theorem \ref{cs4thm3},
and $\Ga(\al)\subset U_\sXp$ is the graph of a small
1-form $\al$ on~$X'$.

Since $\hat X'$ is Lagrangian, $\al$ is {\it closed},
as in \S\ref{cs21}. Also, if $\phi_i,\eta_i$ and
$\hat\phi_i,\hat\eta_i$ are as in Theorem \ref{cs4thm2}
for $X,\hat X$ then $(\Up_i\circ\phi_i)^*(\al)=
\hat\eta_i-\eta_i$ on $\Si_i\t(0,R')$, so applying
Theorem \ref{cs4thm4} to $X,\hat X$ shows that if
$i=1,\ldots,n$ and $\mu_i'\in(2,3)$ with $(2,\mu_i']
\cap\D_\sSii=\emptyset$ then
\e
\bmd{\na^k\al(x)}=O\bigl(d(x,x_i)^{\mu_i'-1-k}\bigr)
\quad\text{near $x_i$ for all $k\ge 0$.}
\label{cs5eq2}
\e

As $\al$ is closed it has a cohomology class $[\al]\in H^1(X',\R)$.
Since \eq{cs5eq2} implies that $\al$ decays quickly near $x_i$, it
turns out that $\al$ must be {\it exact} near $x_i$. Therefore
$[\al]$ can be represented by a compactly-supported form on $X'$,
and lies in the image of $H^1_{\rm cs}(X',\R)$ in~$H^1(X',\R)$.

Choose a vector space $\I_\sXp$ of compactly-supported 1-forms
on $X'$ representing the image of $H^1_{\rm cs}(X',\R)$ in
$H^1(X',\R)$. Then we can write $\al=\be+\d f$, where $\be\in
\I_\sXp$ with $[\al]=[\be]$ is unique, and $f\in C^\iy(X')$ is
unique up to addition of constants. As $\hat X'$ is special
Lagrangian we find that $f$ satisfies a {\it second-order
nonlinear elliptic p.d.e.} similar to \eq{cs4eq7}:
\e
\d^*\bigl(\psi^m(\be+\d f)\bigr)(x)=
Q\bigl(x,(\be+\d f)(x),(\na\be+\na^2f)(x)\bigr)
\label{cs5eq3}
\e
for $x\in X'$. The {\it linearization} of \eq{cs5eq3} at
$\be=f=0$ is~$\d^*\bigl(\psi^m(\be+\d f)\bigr)=0$.

We study the family of small solutions $\be,f$ of \eq{cs5eq3}
for which $f$ has the decay near $x_i$ required by \eq{cs5eq2}.
There is a ready-made theory of analysis on manifolds with
cylindrical ends developed by Lockhart and McOwen \cite{Lock},
which is well-suited to this task. We work on certain {\it weighted
Sobolev spaces} $L^p_{k,{\bs\mu}}(X')$ of functions on~$X'$.

By results from \cite{Lock} it turns out that the
operator $f\mapsto\d^*(\psi^m\d f)$ is a {\it Fredholm} map
$L^p_{k,{\bs\mu}}(X')\ra L^p_{k-2,{\bs\mu}-2}(X')$, with
cokernel of dimension $\sum_{i=1}^nN_\sSii(2)$. This cokernel
is in effect the {\it obstruction space} to deforming $X$ with
$x_i,\up_i$ fixed, as it is the obstruction space to solving
the linearization of \eq{cs5eq3} in $f$ at~$\be=f=0$.

By varying the $x_i$ and $\up_i$, and allowing $f$ to converge
to different constant values on the ends of $X'$ rather than
zero, we can overcome many of these obstructions. This reduces
the dimension of the obstruction space $\O_\sXp$ from
$\sum_{i=1}^nN_\sSii(2)$ to $\sum_{i=1}^n\sind(C_i)$.  The
obstruction map $\Phi$ is constructed using the Implicit
Mapping Theorem for Banach spaces. This concludes our sketch.

If the $C_i$ are {\it stable} then $\O_\sXp=\{0\}$ and we
deduce~\cite[Cor.~6.11]{Joyc10}:

\begin{cor} Suppose $(M,J,\om,\Om)$ is an almost Calabi--Yau
$m$-fold and\/ $X$ a compact SL\/ $m$-fold in $M$ with stable
conical singularities, and let\/ $\M_\sX$ and\/ $\I_\sXp$ be
as in Theorem \ref{cs5thm1}. Then $\M_\sX$ is a smooth manifold
of dimension~$\dim\I_\sXp$.
\label{cs5cor1}
\end{cor}

This has clear similarities with Theorem \ref{cs2thm3}. Here
is another simple condition for $\M_\sX$ to be a manifold
near $X$, \cite[Def.~6.12]{Joyc10}.

\begin{dfn} Let $(M,J,\om,\Om)$ be an almost Calabi--Yau
$m$-fold and $X$ a compact SL $m$-fold in $M$ with conical
singularities, and let $\I_\sXp,\O_\sXp,U$ and $\Phi$ be as
in Theorem \ref{cs5thm1}. We call $X$ {\it transverse} if
the linear map $\d\Phi\vert_0:\I_\sXp\ra\O_\sXp$ is surjective.
\label{cs5def2}
\end{dfn}

If $X$ is transverse then $\{u\in U:\Phi(u)=0\}$ is a manifold
near 0, so Theorem \ref{cs5thm1} yields~\cite[Cor.~6.13]{Joyc10}:

\begin{cor} Suppose $(M,J,\om,\Om)$ is an almost Calabi--Yau $m$-fold
and\/ $X$ a transverse compact SL\/ $m$-fold in $M$ with conical
singularities, and let\/ $\M_\sX,\I_\sXp$ and\/ $\O_\sXp$ be as
in Theorem \ref{cs5thm1}. Then $\M_\sX$ is near $X$ a smooth manifold
of dimension~$\dim\I_\sXp-\dim\O_\sXp$.
\label{cs5cor2}
\end{cor}

In \cite[\S 7]{Joyc10} we extend all this to {\it families} of
almost Calabi--Yau $m$-folds. Combining Definitions \ref{cs2def9}
and \ref{cs5def1}, we define moduli spaces in families:

\begin{dfn} Let $(M,J,\om,\Om)$ be an almost Calabi--Yau
$m$-fold and $X$ a compact SL $m$-fold in $M$ with conical
singularities at $x_1,\ldots,x_n$. Suppose $\bigl\{(M,J^s,
\om^s,\Om^s):s\in\F\bigr\}$ is a smooth family of deformations
of $(M,J,\om,\Om)$. Define the {\it moduli space $\M_\sX^\sF$
of deformations of\/ $X$ in the family} $\F$ to be the set of
pairs $(s,\hat X)$ such that
\begin{itemize}
\setlength{\parsep}{0pt}
\setlength{\itemsep}{0pt}
\item[(i)] $s\in\F$ and $\hat X$ is a compact SL $m$-fold in
$(M,J^s,\om^s,\Om^s)$ with conical singularities at $\hat x_1,
\ldots,\hat x_n$ with cones $C_1,\ldots,C_n$, for some~$\hat x_i$.
\item[(ii)] There exists a homeomorphism $\hat\iota:X\ra\hat X$
with $\hat\iota(x_i)=\hat x_i$ for $i=1,\ldots,n$ such that
$\hat\iota\vert_{X'}:X'\ra\hat X'$ is a diffeomorphism and
$\hat\iota$ and $\iota$ are isotopic as continuous maps
$X\ra M$, where $\iota:X\ra M$ is the inclusion.
\end{itemize}
Define a {\it projection\/} $\pi^\sF:\M_\sX^\sF\ra\F$ by
$\pi^\sF(s,\hat X)=s$. In \cite[Def.~7.5]{Joyc10} we define a
natural {\it topology\/} on $\M_\sX^\sF$, for which $\pi^\sF$
is continuous.
\label{cs5def3}
\end{dfn}

Here \cite[Th.~7.9]{Joyc10} is the families analogue of
Theorem~\ref{cs5thm1}.

\begin{thm} Suppose $(M,J,\om,\Om)$ is an almost Calabi--Yau
$m$-fold and\/ $X$ a compact SL\/ $m$-fold in $M$ with conical
singularities at\/ $x_1,\ldots,x_n$. Let\/ $\M_\sX,X',
\allowbreak
\I_\sXp,\allowbreak
\O_\sXp,\allowbreak
U,\allowbreak
\Phi$ and $\Xi$ be as in Theorem~\ref{cs5thm1}.

Suppose $\bigl\{(M,J^s,\om^s,\Om^s):s\in\F\bigr\}$ is
a smooth family of deformations of\/ $(M,J,\om,\Om)$, in the
sense of Definition \ref{cs2def8}, such that\/ $\iota_*(\ga)
\cdot[\om^s]=0$ for all\/ $\ga\in H_2(X,\R)$ and\/ $s\in\F$,
where $\iota:X\ra M$ is the inclusion, and\/ $[X]\cdot[\Im\Om^s]=0$
for all\/ $s\in\F$, where $[X]\in H_m(M,\R)$ and\/
$[\Im\Om^s]\in H^m(M,\R)$. Let\/ $\M_\sX^\sF$ and\/
$\pi^\sF:\M_\sX^\sF\ra\F$ be as in Definition~\ref{cs5def3}.

Then there exists an open neighbourhood\/ $U^\sF$ of\/ $(0,0)$
in $\F\t U$, a smooth map $\Phi^\sF:U^\sF\ra\O_\sXp$ with\/
$\Phi^\sF(0,u)\equiv\Phi(u)$, and a map $\Xi^\sF:\{(s,u)\in
U^\sF:\Phi^\sF(s,u)=0\}\ra\M_\sX^\sF$ with\/ $\Xi^\sF(0,u)\equiv
\bigl(0,\Xi(u)\bigr)$ and\/ $\pi^\sF\circ\Xi^\sF(s,u)\equiv s$,
which is a homeomorphism with an open neighbourhood of\/
$(0,X)$ in~$\M_\sX^\sF$.
\label{cs5thm2}
\end{thm}

The conditions $\iota_*(\ga)\cdot[\om^s]=0$ for all $\ga$ and
$[X]\cdot[\Im\Om^s]=0$ are {\it necessary conditions} for the
existence of any SL $m$-fold $\hat X$ with conical singularities
isotopic to $X$ in $(M,J^s,\om^s,\Om^s)$. Here are the families
analogues of Definition \ref{cs5def2} and Corollaries \ref{cs5cor1} and
\ref{cs5cor2}, taken from~\cite[Def.~7.11 \& Cor.s 7.10 \& 7.12]{Joyc10}.

\begin{cor} In the situation of Theorem \ref{cs5thm2},
suppose $X$ has stable singularities. Then $\M_\sX^\sF$
is a smooth manifold of dimension\/ $d+\dim\I_\sXp$ and\/
$\pi^\sF:\M_\sX^\sF\ra\F$ a smooth submersion. For small\/
$s\in\F$ the fibre $(\pi^\sF)^{-1}(s)$ is a nonempty smooth
manifold of dimension $\dim\I_\sXp$,
with\/~$(\pi^\sF)^{-1}(0)=\M_\sX$.
\label{cs5cor3}
\end{cor}

Note the similarity of Corollary \ref{cs5cor3} and
Theorem~\ref{cs2thm4}.

\begin{dfn} In the situation of Definition \ref{cs5thm2},
we call $X$ {\it transverse in} $\F$ if the linear map
$\d\Phi^\sF\vert_{(0,0)}:\R^d\t\I_\sXp\ra\O_\sXp$ is
surjective. If $X$ is transverse in the sense of Definition
\ref{cs5def2} then it is also transverse in~$\F$.
\label{cs5def4}
\end{dfn}

\begin{cor} In the situation of Theorem \ref{cs5thm2}, suppose
$X$ is transverse in $\F$. Then $\M_\sX^\sF$ is near $(0,X)$ a
smooth manifold of dimension\/ $d+\dim\I_\sXp-\dim\O_\sXp$, and\/
$\pi^\sF:\M_\sX^\sF\ra\F$ is a smooth map near~$(0,X)$.
\label{cs5cor4}
\end{cor}

Now there are a number of well-known moduli space problems in
geometry where in general moduli spaces are obstructed and
singular, but after a generic perturbation they become smooth
manifolds. For instance, moduli spaces of instantons on
4-manifolds can be made smooth by choosing a generic metric,
and similar things hold for Seiberg--Witten equations, and
moduli spaces of pseudo-holomorphic curves in symplectic
manifolds.

In \cite[\S 9]{Joyc10} we try (but do not quite succeed) to
replicate this for moduli spaces of SL $m$-folds with conical
singularities, by choosing a {\it generic K\"ahler metric} in a
fixed K\"ahler class. Our first result is~\cite[Th.~9.1]{Joyc10}:

\begin{thm} Let\/ $(M,J,\om,\Om)$ be an almost Calabi--Yau
$m$-fold, $X$ be a compact SL\/ $m$-fold in $M$ with conical
singularities, and\/ $\I_\sXp,\O_\sXp$ be as in Theorem
\ref{cs5thm1}. Then there exists a smooth family of deformations
$\bigl\{(M,J,\om^s,\Om):s\in\F\bigr\}$ of\/ $(M,J,\om,\Om)$
with\/ $[\om^s]=[\om]\in H^2(M,\R)$ for all\/ $s\in\F$, such
that\/ $X$ is transverse in $\F$ in the sense of Definition
\ref{cs5def4}, and\/~$d=\dim\F=\dim\O_\sXp$.
\label{cs5thm3}
\end{thm}

Combining this with Corollary \ref{cs5cor4} we see that
$\M_\sX^\sF$ is a manifold near $(0,X)$ and $\pi^\sF$ a
smooth map near $(0,X)$. It then follows from Sard's Theorem
that for small generic $s\in\F$, the moduli space
$(\pi^\sF)^{-1}(s)$ of deformations of $X$ in $(M,J,\om^s,\Om)$
is a smooth manifold near~$X$.

Thus, given a compact SL $m$-fold $X$ with conical singularities
in $(M,J,\om,\Om)$ we can perturb $\om$ a little bit in its
K\"ahler class to $\om^s$, and the moduli space $\M_\sX^s$ in
$(M,J,\om^s,\Om)$ will be a smooth manifold near $X$. More
generally \cite[Th.~9.3]{Joyc10}, if $W\subseteq\M_\sX$ is
a compact subset then we can perturb $\om$ to $\om^s$ so
$\M_\sX^s$ is a smooth manifold near~$W$.

We would like to conclude that by choosing a sufficiently
generic perturbation $\om^s$ we can make $\M_\sX^s$
smooth everywhere. This is the idea of the following
conjecture,~\cite[Conj.~9.5]{Joyc10}:

\begin{conj} Let\/ $(M,J,\om,\Om)$ be an almost Calabi--Yau
$m$-fold, $X$ a compact SL\/ $m$-fold in $M$ with conical
singularities, and\/ $\I_\sXp,\O_\sXp$ be as in Theorem
\ref{cs5thm1}. Then for a second category subset of K\"ahler
forms $\check\om$ in the K\"ahler class of\/ $\om$, the moduli
space $\check\M_\sX$ of compact SL\/ $m$-folds $\hat X$ with
conical singularities in $(M,J,\check\om,\Om)$ isotopic to $X$
is a manifold of dimension~$\dim\I_\sXp-\dim\O_\sXp$.
\label{cs5conj}
\end{conj}

If we could treat the moduli spaces $\M_\sX$ as compact,
say if we had a good understanding of the compactification
$\oM_\sX$ of $\M_\sX$ in \S\ref{cs8}, then this would follow
from \cite[Th.~9.3]{Joyc10}. However, without knowing
$\M_\sX$ is compact, the condition that $\check\M_\sX$ is
smooth everywhere is in effect the intersection of an
infinite number of genericity conditions on $\check\om$,
and we do not know that this intersection is dense (or even
nonempty) in the K\"ahler class.

Notice that Conjecture \ref{cs5conj} constrains the topology
and cones of SL $m$-folds $X$ with conical singularities that
can occur in a generic almost Calabi--Yau $m$-fold, as we
must have~$\dim\I_\sXp\ge\dim\O_\sXp$.

\section{Asymptotically Conical SL $m$-folds}
\label{cs6}

We now discuss {\it Asymptotically Conical\/} SL $m$-folds $L$
in $\C^m$, \cite[Def.~7.1]{Joyc9}.

\begin{dfn} Let $C$ be a closed SL cone in $\C^m$ with isolated
singularity at 0 for $m>2$, and let $\Si=C\cap{\cal S}^{2m-1}$,
so that $\Si$ is a compact, nonsingular $(m-1)$-manifold, not
necessarily connected. Let $g_\sSi$ be the metric on $\Si$
induced by the metric $g'$ on $\C^m$ in \eq{cs2eq2}, and $r$ the
radius function on $\C^m$. Define $\iota:\Si\t(0,\iy)\ra\C^m$ by
$\iota(\si,r)=r\si$. Then the image of $\iota$ is $C\sm\{0\}$,
and $\iota^*(g')=r^2g_\sSi+\d r^2$ is the cone metric
on~$C\sm\{0\}$.

Let $L$ be a closed, nonsingular SL $m$-fold in $\C^m$. We
call $L$ {\it Asymptotically Conical (AC)} with {\it rate}
$\la<2$ and {\it cone} $C$ if there exists a compact subset
$K\subset L$ and a diffeomorphism $\vp:\Si\t(T,\iy)\ra L\sm K$
for some $T>0$, such that
\e
\bmd{\na^k(\vp-\iota)}=O(r^{\la-1-k})
\quad\text{as $r\ra\iy$ for $k=0,1$.}
\label{cs6eq1}
\e
Here $\na,\md{\,.\,}$ are computed using the cone metric~$\iota^*(g')$.
\label{cs6def1}
\end{dfn}

This is very similar to Definition \ref{cs3def4}, and in fact
there are strong parallels between the theories of SL $m$-folds
with conical singularities and of Asymptotically Conical SL
$m$-folds. We continue to assume $m>2$ throughout.

In \S\ref{cs61}--\S\ref{cs62} we review the results of
\cite[\S 7]{Joyc9} on AC SL $m$-folds. Section \ref{cs63}
covers the {\it deformation theory\/} of AC SL $m$-folds in
$\C^m$ following Marshall \cite{Mars} and Pacini \cite{Paci},
and \S\ref{cs64} discusses {\it examples} of AC SL $m$-folds.

\subsection{Cohomological invariants of AC SL $m$-folds}
\label{cs61}

Let $L$ be an AC SL $m$-fold in $\C^m$ with cone $C$, and set
$\Si=C\cap{\cal S}^{2m-1}$. Using the notation of \S\ref{cs34},
as in \eq{cs3eq13} there is a long exact sequence
\e
\cdots\ra
H^k_{\rm cs}(L,\R)\ra H^k(L,\R)\ra H^k(\Si,\R)\ra
H^{k+1}_{\rm cs}(L,\R)\ra\cdots.
\label{cs6eq2}
\e
Following \cite[Def.~7.2]{Joyc9} we define {\it cohomological
invariants\/} $Y(L),Z(L)$ of~$L$.

\begin{dfn} Let $L$ be an AC SL $m$-fold in $\C^m$ with cone $C$,
and let $\Si=C\cap{\cal S}^{2m-1}$. As $\om',\Im\Om'$ in \eq{cs2eq2}
are closed forms with $\om'\vert_L\equiv\Im\Om'\vert_L\equiv 0$, they
define classes in the relative de Rham cohomology groups $H^k(\C^m;
L,\R)$ for $k=2,m$. But for $k>1$ we have the exact sequence
\begin{equation*}
0=H^{k-1}(\C^m,\R)\ra H^{k-1}(L,\R){\buildrel\cong\over\longra}
H^k(\C^m;L,\R)\ra H^k(\C^m,\R)=0.
\end{equation*}
Let $Y(L)\in H^1(\Si,\R)$ be the image of $[\om']$ in
$H^2(\C^m;L,\R)\cong H^1(L,\R)$ under $H^1(L,\R)\ra H^1(\Si,R)$
in \eq{cs6eq2}, and $Z(L)\in H^{m-1}(\Si,\R)$ be the image of
$[\Im\Om']$ in $H^m(\C^m;L,\R)\cong H^{m-1}(L,\R)$ under
$H^{m-1}(L,\R)\ra H^{m-1}(\Si,R)$ in~\eq{cs6eq2}.
\label{cs6def2}
\end{dfn}

Here are some conditions for $Y(L)$ or $Z(L)$ to be
zero,~\cite[Prop.~7.3]{Joyc9}.

\begin{prop} Let\/ $L$ be an AC SL\/ $m$-fold in $\C^m$ with
cone $C$ and rate $\la$, and let\/ $\Si=C\cap{\cal S}^{2m-1}$.
If\/ $\la<0$ or $b^1(L)=0$ then $Y(L)=0$. If\/ $\la<2-m$ or
$b^0(\Si)=1$ then~$Z(L)=0$.
\label{cs6prop}
\end{prop}

\subsection{Lagrangian Neighbourhood Theorems and regularity}
\label{cs62}

Next we give versions of parts of \S\ref{cs41}--\S\ref{cs42} 
for AC SL $m$-folds rather than SL $m$-folds with conical
singularities. Here are the analogues of Theorems \ref{cs4thm2}
and \ref{cs4thm3}, proved in~\cite[Th.s~7.4 \& 7.5]{Joyc9}.

\begin{thm} Let\/ $C$ be an SL cone in $\C^m$ with isolated
singularity at\/ $0$, and set\/ $\Si=C\cap{\cal S}^{2m-1}$.
Define $\iota:\Si\t(0,\iy)\ra\C^m$ by $\iota(\si,r)=r\si$. Let\/
$\ze$, $U_\sC\subset T^*\bigl(\Si\t(0,\iy)\bigr)$ and\/
$\Phi_\sC:U_\sC\ra\C^m$ be as in Theorem~\ref{cs4thm1}.

Suppose $L$ is an AC SL\/ $m$-fold in $\C^m$ with cone $C$
and rate $\la<2$. Then there exists a compact\/ $K\subset L$
and a diffeomorphism $\vp:\Si\t(T,\iy)\ra L\sm K$ for some
$T>0$ satisfying \eq{cs6eq1}, and a closed\/ $1$-form $\chi$ on
$\Si\t(T,\iy)$ written $\chi(\si,r)=\chi^1(\si,r)+\chi^2(\si,r)\d r$
for $\chi^1(\si,r)\in T_\si^*\Si$ and\/ $\chi^2(\si,r)\in\R$,
satisfying
\e
\begin{gathered}
\bmd{\chi(\si,r)}<\ze r,\quad \vp(\si,r)\equiv
\Phi_\sC\bigl(\si,r,\chi^1(\si,r),\chi^2(\si,r)\bigr)\\
\text{and}\quad\bmd{\na^k\chi}=O(r^{\la-1-k})
\quad\text{as $r\ra\iy$ for $k=0,1,$}
\end{gathered}
\label{cs6eq3}
\e
computing $\na,\md{\,.\,}$ using the cone metric~$\iota^*(g')$.
\label{cs6thm1}
\end{thm}

\begin{thm} Suppose $L$ is an AC SL\/ $m$-fold in $\C^m$ with cone
$C$. Let\/ $\Si,
\allowbreak
\iota,
\allowbreak
\ze,
\allowbreak
U_{\sst C},
\allowbreak
\Phi_{\sst C},K,T,\vp,\chi,\chi^1,\chi^2$ be as in Theorem \ref{cs6thm1}.
Then making $T,K$ larger if necessary, there exists an open tubular
neighbourhood\/ $U_{\sst L}\subset T^*L$ of the zero section $L$ in $T^*L$,
such that under $\d\vp:T^*\bigl(\Si\t(T,\iy)\bigr)\ra T^*L$ we have
\e
(\d\vp)^*(U_{\sst L})=\bigl\{(\si,r,\tau,u)\in
T^*\bigl(\Si\t(T,\iy)\bigr):\bmd{(\tau,u)}<\ze r\bigr\},
\label{cs6eq4}
\e
and there exists an embedding $\Phi_{\sst L}:U_{\sst L}\ra\C^m$ with\/
$\Phi_{\sst L}\vert_L=\id:L\ra L$ and\/ $\Phi_{\sst L}^*(\om')=\hat\om$,
where $\hat\om$ is the canonical symplectic structure on $T^*L$,
such that
\e
\Phi_{\sst L}\circ\d\vp(\si,r,\tau,u)\equiv
\Phi_{\sst C}\bigl(\si,r,\tau+\chi^1(\si,r),u+\chi^2(\si,r)\bigr)
\label{cs6eq5}
\e
for all\/ $(\si,r,\tau,u)\!\in\!T^*\bigl(\Si\!\t\!(T,\iy)\bigr)$ with\/
$\md{(\tau,u)}<\ze r$, computing $\md{\,.\,}$ using~$\iota^*(g')$.
\label{cs6thm2}
\end{thm}

Combining \cite[Prop.~7.6]{Joyc9} and \cite[Th.s 7.7 \& 7.11]{Joyc9}
gives an analogue of Theorem \ref{cs4thm4}, on the {\it regularity}
of $L$ near infinity in $\C^m$. As in \cite[Th.~7.11]{Joyc9}, the
theorem can be strengthened when~$0\le\la<\min\bigl(\D_\Si\cap
(0,\iy)\bigr)$.

\begin{thm} In Theorem \ref{cs6thm1} we have $[\chi]=Y(L)$ in
$H^1\bigl(\Si\t(T,\iy),\R\bigr)\cong H^1(\Si,\R)$, where $Y(L)$
is as in Definition \ref{cs6def2}. Let\/ $\ga$ be the unique
$1$-form on $\Si$ with\/ $\d\ga=\d^*\ga=0$ and\/ $[\ga]=Y(L)\in
H^1(\Si,\R)$, which exists by Hodge theory. Then $\chi=\pi^*(\ga)
+\d E$, where $\pi:\Si\t(T,\iy)\ra\Si$ is the projection
and\/~$E\in C^\iy\bigl(\Si\t(T,\iy)\bigr)$.

If either $\la=\la'$, or $\la,\la'$ lie in the same connected
component of\/ $\R\sm\D_\sSi$, then $L$ is an AC SL\/ $m$-fold
with rate $\la'$ and
\e
\begin{gathered}
\bmd{\na^k(\vp\!-\!\iota)}=O(r^{\la'-1-k}),\;\>
\bmd{\na^k\chi}=O(r^{\la'-1-k}),\;\>
\bmd{\na^{k+1}E}=O(r^{\la'-1-k})\\
\text{for all $k\ge 0$, and}\quad
\md{E}=\begin{cases}O(r^{\la'}), & \la'\ne 0, \\
O\bigl(\md{\log r}\bigr), & \la'=0. \end{cases}
\end{gathered}
\label{cs6eq6}
\e
Here $\na,\md{\,.\,}$ are computed using the cone
metric $\iota^*(g')$ on~$\Si\t(T,\iy)$.
\label{cs6thm3}
\end{thm}

\subsection{Moduli spaces of AC SL $m$-folds}
\label{cs63}

The deformation theory of Asymptotically Conical SL $m$-folds
in $\C^m$ has been studied independently by Pacini \cite{Paci}
and Marshall \cite{Mars}. Pacini's results are earlier, but
Marshall's are more complete.

\begin{dfn} Suppose $L$ is an Asymptotically Conical SL
$m$-fold in $\C^m$ with cone $C$ and rate $\la<2$, as in
Definition \ref{cs6def1}. Define the {\it moduli space
$\M_\sL^\la$ of deformations of\/ $L$ with rate} $\la$
to be the set of AC SL $m$-folds $\hat L$ in $\C^m$ with
cone $C$ and rate $\la$, such that $\hat L$ is diffeomorphic
to $L$ and isotopic to $L$ as an Asymptotically Conical
submanifold of $\C^m$. One can define a natural {\it topology}
on $\M_\sL^\la$, in a similar way to the conical singularities
case of~\cite[Def.~5.6]{Joyc10}.
\label{cs6def3}
\end{dfn}

Note that if $L$ is an AC SL $m$-fold with rate $\la$, then
it is {\it also} an AC SL $m$-fold with rate $\la'$ for any
$\la'\in[\la,2)$. Thus we have defined a 1-{\it parameter
family} of moduli spaces $\M_\sL^{\smash{\la'}}$ for $L$, and
not just one. Since we did not impose any condition on $\la$
in Definition \ref{cs6def1} analogous to \eq{cs3eq11} in the
conical singularities case, it turns out that $\M_\sL^\la$
depends nontrivially on~$\la$.

The following result can be deduced from Marshall
\cite[Th.~6.2.15]{Mars} and \cite[Table~5.1]{Mars}.
(See also Pacini \cite[Th.~2 \& Th.~3]{Paci}.) It implies
conjectures by the author in \cite[Conj.~2.12]{Joyc1}
and~\cite[\S 10.2]{Joyc8}.

\begin{thm} Let\/ $L$ be an Asymptotically Conical SL\/
$m$-fold in $\C^m$ with cone $C$ and rate $\la<2$, and
let\/ $\M_\sL^\la$ be as in Definition \ref{cs6def3}.
Set\/ $\Si=C\cap{\cal S}^{2m-1}$, and let\/ $\D_\sSi,
N_\sSi$ be as in \S\ref{cs31} and\/ $b^k(L),b^k_{\rm
cs}(L)$ as in \S\ref{cs34}. Then
\begin{itemize}
\item[{\rm(a)}] If\/ $\la\in(0,2)\sm\D_\sSi$ then
$\M_\sL^\la$ is a manifold with
\e
\dim\M_\sL^\la=b^1(L)-b^0(L)+N_\sSi(\la).
\label{cs6eq7}
\e
Note that if\/ $0<\la<\min\bigl(\D_\sSi\cap
(0,\iy)\bigr)$ then~$N_\sSi(\la)=b^0(\Si)$.
\item[{\rm(b)}] If\/ $\la\in(2-m,0)$ then $\M_\sL^\la$
is a manifold of dimension~$b^1_{\rm cs}(L)=b^{m-1}(L)$.
\end{itemize}
\label{cs6thm4}
\end{thm}

This is the analogue of Theorems \ref{cs2thm3} and \ref{cs5thm1}
for AC SL $m$-folds. If $\la\in(2-m,2)\sm\D_\sSi$ then the
deformation theory for $L$ with rate $\la$ is {\it unobstructed\/}
and $\M_\sL^\la$ is a {\it smooth manifold\/} with a given
dimension. This is similar to the case of nonsingular compact
SL $m$-folds in Theorem \ref{cs2thm3}, but different to the
conical singularities case in Theorem~\ref{cs5thm1}.

\subsection{Examples}
\label{cs64}

Examples of AC SL $m$-folds $L$ are constructed by Harvey and Lawson
\cite[\S III.3]{HaLa}, the author \cite{Joyc2,Joyc3,Joyc4,Joyc6},
and others. Nearly all the known examples (up to translations) have
minimum rate $\la$ either 0 or $2-m$, which are topologically
significant values by Proposition \ref{cs6prop}. For instance, all
examples in \cite{Joyc3} have $\la=0$, and \cite[Th.~6.4]{Joyc2}
constructs AC SL $m$-folds with $\la=2-m$ in $\C^m$ from any SL
cone $C$ in $\C^m$. The only explicit, nontrivial examples known
to the author with $\la\ne 0,2-m$ are in \cite[Th.~11.6]{Joyc4},
and have~$\la=\frac{3}{2}$.

We shall give three families of examples of AC SL $m$-folds $L$
in $\C^m$ explicitly. The first family is adapted from Harvey
and Lawson~\cite[\S III.3.A]{HaLa}.

\begin{ex} Let $C_{\sst\rm HL}^m$ be the SL cone in $\C^m$ of
Example \ref{cs3ex1}. We shall define a family of AC SL $m$-folds
in $\C^m$ with cone $C_{\sst\rm HL}^m$. Let $a_1,\ldots,a_m\ge 0$
with exactly two of the $a_j$ zero and the rest positive. Write
${\bf a}=(a_1,\ldots,a_m)$, and define
\e
\begin{split}
L_{\sst\rm HL}^{\bf a}=\bigl\{&(z_1,\ldots,z_m)\in\C^m:
i^{m+1}z_1\cdots z_m\in[0,\iy),\\
&\ms{z_1}-a_1=\cdots=\ms{z_m}-a_m\bigr\}.
\end{split}
\label{cs6eq8}
\e
Then $L_{\sst\rm HL}^{\bf a}$ is an AC SL $m$-fold in $\C^m$
diffeomorphic to $T^{m-2}\t\R^2$, with cone $C_{\sst\rm HL}^m$
and rate 0. It is invariant under the $\U(1)^{m-1}$ group
\eq{cs3eq7}. It is surprising that equations of the form
\eq{cs6eq8} should define a nonsingular submanifold of
$\C^m$ {\it without boundary}, but in fact they do.

Now suppose for simplicity that $a_1,\ldots,a_{m-2}>0$ and
$a_{m-1}=a_m=0$. As $\Si_{\sst\rm HL}^m\cong T^{m-1}$ we have
$H^1(\Si_{\sst\rm HL}^m,\R)\cong\R^{m-1}$, and calculation
shows that $Y(L_{\sst\rm HL}^{\bf a})=(\pi a_1,\ldots,\pi
a_{m-2},0)\in\R^{m-1}$ in the natural coordinates. Since
$L_{\sst\rm HL}^{\bf a}\cong T^{m-2}\t\R^2$ we have
$H^1(L_{\sst\rm HL}^{\bf a},\R)=\R^{m-2}$, and
$Y(L_{\sst\rm HL}^{\bf a})$ lies in the image $\R^{m-2}
\subset\R^{m-1}$ of $H^1(L_{\sst\rm HL}^{\bf a},\R)$ in
$H^1(\Si_{\sst\rm HL}^m,\R)$, as in Definition
\ref{cs6def2}. As $b^0(\Si_{\sst\rm HL}^m)=1$, Proposition
\ref{cs6prop} shows that~$Z(L_{\sst\rm HL}^{\bf a})=0$.

Take $C=C_{\sst\rm HL}^m$, $\Si=\Si_{\sst\rm HL}^m$ and
$L=L_{\sst\rm HL}^{\bf a}$ in Theorem \ref{cs6thm4},
and let $0<\la<\min\bigl(\D_\sSi\cap(0,\iy)\bigr)$. Then
$b^1(L)=m-2$, $b^0(L)=1$ and $N_\sSi(\la)=b^0(\Si)=1$, so part
(a) of Theorem \ref{cs6thm4} shows that $\dim\M_\sL^\la=m-2$.
This is consistent with the fact that $L$ depends on $m-2$
real parameters~$a_1,\ldots,a_{m-2}>0$.

The family of all $L_{\sst\rm HL}^{\bf a}$ has $\ha m(m-1)$
connected components, indexed by which two of $a_1,\ldots,a_m$
are zero. Using the theory of \S\ref{cs7}, these can give many
{\it topologically distinct\/} ways to desingularize SL
$m$-folds with conical singularities with these cones.
\label{cs6ex1}
\end{ex}

Our second family, from \cite[Ex.~9.4]{Joyc2}, was chosen as
it is easy to write down.

\begin{ex} Let $m,a_1,\ldots,a_m,k$ and $L^{a_1,\ldots,a_m}_0$
be as in Example \ref{cs3ex2}. For $0\ne c\in\R$ define
\e
\begin{split}
L^{a_1,\ldots,a_m}_c=\bigl\{
\bigl(i{\rm e}^{ia_1\th}x_1&,{\rm e}^{ia_2\th}x_2,\ldots,
{\rm e}^{ia_m\th}x_m\bigr):
\th\in[0,2\pi),\\ 
&x_1,\ldots,x_m\in\R,\qquad 
a_1x_1^2+\cdots+a_mx_m^2=c\bigr\}.
\end{split}
\label{cs6eq9}
\e
Then $L^{a_1,\ldots,a_m}_c$ is an AC SL $m$-fold in $\C^m$ with
rate 0 and cone $L^{a_1,\ldots,a_m}_0$. It is diffeomorphic as an
immersed SL $m$-fold to $({\cal S}^{k-1}\t\R^{m-k}\t{\cal S}^1)/\Z_2$
if $c>0$, and to $(\R^k\t{\cal S}^{m-k-1}\t{\cal S}^1)/\Z_2$ if~$c<0$.
\label{cs6ex2}
\end{ex}

Our third family was first found by Lawlor \cite{Lawl}, made more
explicit by Harvey \cite[p.~139--140]{Harv}, and discussed from a
different point of view by the author in \cite[\S 5.4(b)]{Joyc3}.
Our treatment is based on that of Harvey.

\begin{ex} Let $m>2$ and $a_1,\ldots,a_m>0$, and define
polynomials $p,P$ by
\begin{equation*}
p(x)=(1+a_1x^2)\cdots(1+a_mx^2)-1
\quad\text{and}\quad P(x)=\frac{p(x)}{x^2}.
\end{equation*}
Define real numbers $\phi_1,\ldots,\phi_m$ and $A$ by
\e
\phi_k=a_k\int_{-\iy}^\iy\frac{\d x}{(1+a_kx^2)\sqrt{P(x)}}
\quad\text{and}\quad A=\om_m(a_1\cdots a_m)^{-1/2},
\label{cs6eq10}
\e
where $\om_m$ is the volume of the unit sphere in $\R^m$. Clearly
$\phi_k,A>0$. But writing $\phi_1+\cdots+\phi_m$ as one integral gives
\begin{equation*}
\phi_1+\cdots+\phi_m=\int_0^\iy\frac{p'(x)\d x}{(p(x)+1)\sqrt{p(x)}}
=2\int_0^\iy\frac{\d w}{w^2+1}=\pi,
\end{equation*}
making the substitution $w=\sqrt{p(x)}$. So $\phi_k\in(0,\pi)$
and $\phi_1+\cdots+\phi_m=\pi$. This yields a 1-1 correspondence
between $m$-tuples $(a_1,\ldots,a_m)$ with $a_k>0$, and
$(m\!+\!1)$-tuples $(\phi_1,\ldots,\phi_m,A)$ with $\phi_k\in
(0,\pi)$, $\phi_1+\cdots+\phi_m=\pi$ and~$A>0$.

For $k=1,\ldots,m$ and $y\in\R$, define a function $z_k:\R\ra\C$ by
\begin{equation*}
z_k(y)={\rm e}^{i\psi_k(y)}\sqrt{a_k^{-1}+y^2}, \quad\text{where}\quad
\psi_k(y)=a_k\int_{-\iy}^y\frac{\d x}{(1+a_kx^2)\sqrt{P(x)}}\,.
\end{equation*}
Now write ${\bs\phi}=(\phi_1,\ldots,\phi_n)$, and define 
a submanifold $L^{{\bs\phi},A}$ in $\C^m$ by
\e
L^{{\bs\phi},A}=\bigl\{(z_1(y)x_1,\ldots,z_m(y)x_m):
y\in\R,\; x_k\in\R,\; x_1^2+\cdots+x_m^2=1\bigr\}.
\label{cs6eq11}
\e

Then $L^{{\bs\phi},A}$ is closed, embedded, and diffeomorphic
to ${\cal S}^{m-1}\t\R$, and Harvey \cite[Th.~7.78]{Harv} shows
that $L^{{\bs\phi},A}$ is {\it special Lagrangian}. One can also
show that $L^{{\bs\phi},A}$ is {\it Asymptotically Conical},
with rate $2-m$ and cone the union $\Pi^0\cup\Pi^{\bs\phi}$ of
two special Lagrangian $m$-planes $\Pi^0,\Pi^{\bs\phi}$ in
$\C^m$ given by
\e
\Pi^0=\bigl\{(x_1,\ldots,x_m):x_j\in\R\bigr\},\;\>
\Pi^{\bs\phi}=\bigl\{({\rm e}^{i\phi_1}x_1,\ldots,
{\rm e}^{i\phi_m}x_m):x_j\in\R\bigr\}.
\label{cs6eq12}
\e

As $\la=2-m<0$ we have $Y(L^{{\bs\phi},A})=0$ by Proposition
\ref{cs6prop}. Now $L^{{\bs\phi},A}\cong{\cal S}^{m-1}\t\R$
so that $H^{m-1}(L^{{\bs\phi},A},\R)\cong\R$, and $\Si=(\Pi^0
\cup\Pi^{\bs\phi})\cap{\cal S}^{2m-1}$ is the disjoint union
of two unit $(m\!-\!1)$-spheres ${\cal S}^{m-1}$, so $H^{m-1}
(\Si,\R)\cong\R^2$. The image of $H^{m-1}(L^{{\bs\phi},A},\R)$
in $H^{m-1}(\Si,\R)$ is $\bigl\{(x,-x):x\in\R\bigr\}$ in the
natural coordinates. Calculation shows that $Z(L^{{\bs\phi},A})=
(A,-A)\in H^{m-1}(\Si,\R)$, which is why we defined $A$ this way
in~\eq{cs6eq10}.

Apply Theorem \ref{cs6thm4} with $L=L^{{\bs\phi},A}$ and
$\la\in(2-m,0)$. As $L\cong{\cal S}^{m-1}\t\R$ we have
$b^1_{\rm cs}(L)=1$, so part (b) of Theorem \ref{cs6thm4} shows
that $\dim\M_\sL^\la=1$. This is consistent with the fact that
when $\bs\phi$ is fixed, $L^{{\bs\phi},A}$ depends on one real
parameter $A>0$. Here $\bs\phi$ is fixed in $\M_\sL^\la$ as the
cone $C=\Pi^0\cup\Pi^{\bs\phi}$ of $L$ depends on $\bs\phi$, and
all $\hat L\in\M_\sL^\la$ have the same cone $C$, by definition.
\label{cs6ex3}
\end{ex}

\section{Desingularizing singular SL $m$-folds}
\label{cs7}

We now discuss the work of \cite{Joyc11,Joyc12} on
{\it desingularizing} compact SL $m$-folds with conical
singularities. Here is the basic idea. Let $(M,J,\om,\Om)$
be an almost Calabi--Yau $m$-fold, and $X$ a compact SL
$m$-fold in $M$ with conical singularities $x_1,\ldots,x_n$
and cones $C_1,\ldots,C_n$. Suppose $L_1,\ldots,L_n$ are AC SL
$m$-folds in $\C^m$ with the same cones $C_1,\ldots,C_n$ as~$X$.

If $t>0$ then $tL_i=\{t\,{\bf x}:{\bf x}\in L_i\}$ is also
an AC SL $m$-fold with cone $C_i$. We construct a 1-parameter
family of compact, nonsingular {\it Lagrangian} $m$-folds $N^t$
in $(M,\om)$ for $t\in(0,\de)$ by gluing $tL_i$ into $X$ at
$x_i$, using a partition of unity.

When $t$ is small, $N^t$ is {\it close to special Lagrangian}
(its phase is nearly constant), but also {\it close to singular}
(it has large curvature and small injectivity radius). We prove
using analysis that for small $t\in(0,\de)$ we can deform $N^t$
to a {\it special\/} Lagrangian $m$-fold $\smash{\ti N^t}$ in
$M$, using a small Hamiltonian deformation. The proof involves
a delicate balancing act, showing that the advantage of being
close to special Lagrangian outweighs the disadvantage of
being nearly singular.

Doing this in full generality is rather complex. There are two
kinds of {\it obstructions} to the existence of $\smash{\ti N^t}$.
Firstly, if $Y(L_i)\ne 0$ then $N^t$ may not exist as a
{\it Lagrangian} $m$-fold. Secondly, if $X'$ is not connected
then we may not be able to deform $N^t$ to a {\it special\/}
Lagrangian $m$-fold $\smash{\ti N^t}$ because of problems with
small eigenvalues of the Laplacian $\De$ on $N^t$. In each case,
$\smash{\ti N^t}$ exists for small $t$ if the $Y(L_i)$ or $Z(L_i)$
satisfy an equation.

We also extend the results to desingularization in {\it families}
of almost Calabi--Yau $m$-folds $(M,J^s,\om^s,\Om^s)$. The
cohomology classes $[\om^s]$ and $[\Im\Om^s]$ contribute to the
obstruction equations in $Y(L_i)$ and $Z(L_i)$ for the existence
of $\smash{\ti N^t}$. Thus, a singular SL $m$-fold $X$ which has
no desingularizations $\smash{\ti N^t}$ in $(M,J,\om,\Om)$ can
still admit desingularizations $\smash{\ti N^{s,t}}$ in
$(M,J^s,\om^s,\Om^s)$ for small~$s\ne 0$.

We begin in \S\ref{cs71} by explaining desingularization
in the simplest case, in one almost Calabi--Yau $m$-fold
$(M,J,\om,\Om)$ when $Y(L_i)=0$ and $X'$ is connected.
Section \ref{cs72} extends this to $X'$ not connected, and
\S\ref{cs73} to $Y(L_i)\ne 0$, introducing the two
kinds of obstructions. Section \ref{cs74} discusses
desingularization in families~$(M,J^s,\om^s,\Om^s)$.

\subsection{Desingularization in the simplest case}
\label{cs71}

Our simplest desingularization result is~\cite[Th.~6.13]{Joyc11}.

\begin{thm} Suppose $(M,J,\om,\Om)$ is an almost Calabi--Yau
$m$-fold and\/ $X$ a compact SL\/ $m$-fold in $M$ with conical
singularities at\/ $x_1,\ldots,x_n$ and cones $C_1,\ldots,C_n$.
Let\/ $L_1,\ldots,L_n$ be Asymptotically Conical SL\/ $m$-folds
in $\C^m$ with cones $C_1,\ldots,C_n$ and rates $\la_1,\ldots,
\la_n$. Suppose $\la_i<0$ for $i=1,\ldots,n$, and\/ $X'=X\sm
\{x_1,\ldots,x_n\}$ is connected.

Then there exists $\ep>0$ and a smooth family $\bigl\{
\smash{\ti N^t}:t\in(0,\ep]\bigr\}$ of compact, nonsingular
SL\/ $m$-folds in $(M,J,\om,\Om)$, such that\/ $\smash{\ti N^t}$
is constructed by gluing $tL_i$ into $X$ at\/ $x_i$ for
$i=1,\ldots,n$. In the sense of currents, $\smash{\ti N^t}\ra
X$ as~$t\ra 0$.
\label{cs7thm1}
\end{thm}

Here is a sketch of the proof, divided into seven steps.
\begin{list}{}{
\setlength{\itemsep}{1pt}
\setlength{\parsep}{1pt}
\setlength{\labelwidth}{35pt}
\setlength{\leftmargin}{35pt}}
\item[Step 1.] Apply Theorem \ref{cs4thm1} to $C_i$ for
$i=1,\ldots,n$, and Theorem \ref{cs4thm3} to $X$, and
Theorem \ref{cs6thm2} to $L_i$ for $i=1,\ldots,n$. This
gives {\it Lagrangian neighbourhoods} $U_\sCi,\Phi_\sCi$
for $C_i$, and $U_\sXp,\Phi_\sXp$ for $X'$, and $U_\sLi,
\Phi_\sLi$ for~$L_i$.

Moreover $U_\sXp,\Phi_\sXp$ and $U_\sCi,\Phi_\sCi$ are
related via $\Up_i$ and an exact 1-form $\eta_i$ on
$\Si_i\t(0,R')$ from Theorem \ref{cs4thm2}, and $U_\sLi,
\Phi_\sLi$ and $U_\sCi,\Phi_\sCi$ are related via an exact
1-form $\chi_i$ on $\Si_i\t(T,\iy)$ from Theorem~\ref{cs6thm1}.

\item[Step 2.] Let $t>0$ be small. We define a nonsingular
Lagrangian $m$-fold $N^t$ in $(M,\om)$, roughly as follows.
Choose $\tau\in(0,1)$ satisfying certain conditions. At
distance at least $2t^\tau$ from $x_1,\ldots,x_n$ we define
$N^t$ to be $X'$. At distance up to $t^\tau$ from $x_i$ we
define $N^t$ to be~$\Up_i(tL_i\cap B_R)$.

Between distances $t^\tau$ and $2t^\tau$ from $x_i$ we define
$N^t$ to be a Lagrangian annulus $\Si_i\t(t^\tau,2t^\tau)$
interpolating between $X'$ and $\Up_i(tL_i\cap B_R)$, using
the Lagrangian neighbourhoods $U_\sCi,\Phi_\sCi$, $U_\sXp,
\Phi_\sXp$ and $U_\sLi,\Phi_\sLi$. This is equivalent to
choosing a closed 1-form $\xi_i^t(\si,r)$ on $\Si_i\t[t^\tau,
2t^\tau]$ which interpolates between $t^2\chi_i(\si,t^{-1}r)$
at $r=t^\tau$ and $\eta_i(\si,r)$ at~$r=2t^\tau$.
\item[Step 3.] Let ${\rm e}^{i\th^t}$ be the phase function of
$N^t$, so that $N^t$ is special Lagrangian if $\sin\th^t\equiv 0$.
We bound various norms of $\psi^m\sin\th^t$ in terms of powers
of $t$. These bounds imply that $N^t$ is {\it close to special
Lagrangian} when $t$ is small. We also estimate other geometrical
quantities, like the curvature and injectivity radius of $N^t$,
in terms of powers of~$t$.
\item[Step 4.] We glue together the Lagrangian neighbourhoods
$U_\sCi,\Phi_\sCi$, $U_\sXp,\Phi_\sXp$ and $U_\sLi,\Phi_\sLi$ to
define a Lagrangian neighbourhood $U_\sNt,\Phi_\sNt$ for~$N^t$.
\item[Step 5.] Let $f\in C^\iy(N^t)$. Then $\d f$ is a 1-form on
$N^t$, and the graph $\Ga(\d f)$ is a submanifold of $T^*N^t$. If
$f$ is small in $C^1$ then $\Ga(\d f)\subset U_\sNt\subset T^*N^t$,
and then $\smash{\ti N^t}=\Phi_\sNt\bigl(\Ga(\d f)\bigr)$ is a
nonsingular Lagrangian $m$-fold in $(M,\om)$. Every small
Hamiltonian deformation of $N^t$ can be written in this way.

We show that $\smash{\ti N^t}$ is {\it special\/} Lagrangian
if and only if
\e
\d^*(\psi^m\cos\th^t\d f)(x)=\psi^m\sin\th^t+
Q^t\bigl(x,\d f(x),\na^2f(x)\bigr)
\label{cs7eq1}
\e
for all $x\in N^t$, as in \eq{cs4eq7} and \eq{cs5eq3}, where $Q^t$
is smooth and $Q^t(x,y,z)=O\bigl(t^{-2}\ms{y}+\ms{z}\bigr)$ for
small~$y,z$.
\item[Step 6.] Working in the Sobolev space $L^{2m}_3(N^t)$,
we show that the operator
\end{list}
\e
\begin{gathered}
P^t:\bigl\{u\in L^{2m}_3(N^t):\ts\int_{N^t}u\,\d V^t=0\bigr\}
\ra\bigl\{v\in L^{2m}_1(N^t):\ts\int_{N^t}v\,\d V^t=0\bigr\}\\
\text{given by}\quad P^t(u)=\d^*(\psi^m\cos\th^t\d u)
\end{gathered}
\label{cs7eq2}
\e
\begin{list}{}{
\setlength{\itemsep}{1pt}
\setlength{\parsep}{1pt}
\setlength{\labelwidth}{35pt}
\setlength{\leftmargin}{35pt}}
\item[]has an inverse $(P^t)^{-1}$ which is (in a rather weak
sense) bounded independently of $t$. The restriction to $u$
with $\int_{N^t}u\,\d V^t=0$ is necessary as $P^t(1)=0$, so
$P^t$ is not invertible on spaces including~1.
\item[Step 7.] We inductively construct a sequence
$(f_k)_{k=0}^\iy$ in $L^{2m}_3(N^t)$ with $f_0=0$,
$\int_{N^t}f_k\,\d V^t\!=\!0$ and $f_k\!=\!(P^t)^{-1}\bigl(
\psi^m\sin\th^t\!+\!Q^t(x,\d f_{k-1},\na^2\d f_{k-1})\bigr)$, so
\e
\d^*(\psi^m\cos\th^t\d f_k)\equiv\psi^m\sin\th^t+
Q^t\bigl(x,\d f_{k-1}(x),\na^2f_{k-1}(x)\bigr).
\label{cs7eq3}
\e

Using the bounds on $\psi^m\sin\th^t$ from Step 3 and on
$(P^t)^{-1}$ from Step 6 we show that $(f_k)_{k=0}^\iy$
exists and converges in $L^{2m}_3(N^t)$ for small $t$.
The limit $f$ satisfies \eq{cs7eq1}, and is smooth by
elliptic regularity. Then $\smash{\ti N^t}=\Phi_\sNt
\bigl(\Ga(\d f)\bigr)$ is the SL $m$-fold we seek.
\end{list}

The condition $\la_i<0$ in Theorem \ref{cs6prop} is there
for two reasons. Firstly, it forces $Y(L_i)=0$ by Proposition
\ref{cs6prop}, and therefore $\chi_i$ is an {\it exact\/}
1-form on $\Si_i\t(T,\iy)$, since $[\chi_i]=Y(L_i)\in
H^1(\Si_i,\R)$ by Theorem \ref{cs6thm3}. This exactness
makes it possible to define the closed 1-form $\xi_i^t$
in Step~2.

Secondly, we need $\la_i<0$ so that the contributions
to $\psi^m\sin\th^t$ from tapering $\chi_i$ off to zero
on the annulus $\Si_i\t[t^\tau,2t^\tau]$ are small
enough for the method to work. If $\la_i>0$ then
$\lnm{\psi^m\sin\th^t}{2}$ is too large, and we cannot
prove that the sequence $(f_k)_{k=0}^\iy$ in Step 7 converges.

\subsection{Desingularization when $X'$ is not connected}
\label{cs72}

In \cite[Th.~7.10]{Joyc11} we extend Theorem \ref{cs7thm1}
to the case when $X'$ is not connected.

\begin{thm} Suppose $(M,J,\om,\Om)$ is an almost Calabi--Yau
$m$-fold and\/ $X$ a compact SL\/ $m$-fold in $M$ with conical
singularities at\/ $x_1,\ldots,x_n$ and cones $C_1,\ldots,C_n$.
Define $\psi:M\ra(0,\iy)$ as in \eq{cs2eq4}. Let\/ $L_1,\ldots,L_n$
be Asymptotically Conical SL\/ $m$-folds in $\C^m$ with cones
$C_1,\ldots,C_n$ and rates $\la_1,\ldots,\la_n$. Suppose $\la_i<0$
for $i=1,\ldots,n$. Write $X'=X\sm\{x_1,\ldots,x_n\}$
and\/~$\Si_i=C_i\cap{\cal S}^{2m-1}$.

Set\/ $q=b^0(X')$, and let\/ $X_1',\ldots,X_q'$ be the connected
components of\/ $X'$. For $i=1,\ldots,n$ let\/ $l_i=b^0(\Si_i)$,
and let\/ $\Si_i^1,\ldots,\Si_i^{\smash{l_i}}$ be the connected
components of\/ $\Si_i$. Define $k(i,j)=1,\ldots,q$ by $\Up_i
\circ\vp_i\bigl(\Si_i^j\t(0,R')\bigr)\subset X'_{\smash{k(i,j)}}$
for $i=1,\ldots,n$ and $j=1,\ldots,l_i$. Suppose that
\e
\sum_{\substack{1\le i\le n, \; 1\le j\le l_i: \\
k(i,j)=k}}\psi(x_i)^mZ(L_i)\cdot[\Si_i^j\,]=0
\quad\text{for all\/ $k=1,\ldots,q$.}
\label{cs7eq4}
\e

Suppose also that the compact\/ $m$-manifold\/ $N$ obtained by
gluing $L_i$ into $X'$ at\/ $x_i$ for $i=1,\ldots,n$ is connected.
A sufficient condition for this to hold is that\/ $X$ and\/ $L_i$
for $i=1,\ldots,n$ are connected.

Then there exists $\ep>0$ and a smooth family $\smash{\bigl\{
\ti N^t:t\in(0,\ep]\bigr\}}$ of compact, nonsingular SL\/
$m$-folds in $(M,J,\om,\Om)$ diffeomorphic to $N$, such
that\/ $\smash{\ti N^t}$ is constructed by gluing $tL_i$
into $X$ at\/ $x_i$ for $i=1,\ldots,n$. In the sense of
currents in Geometric Measure Theory, $\smash{\ti N^t}\ra
X$ as~$t\ra 0$.
\label{cs7thm2}
\end{thm}

The new issues when $X'$ is not connected occur in Steps
6 and 7 of \S\ref{cs71}. Suppose $b^0(X')=q>1$, so that
$X'$ has $q$ connected components $X_1',\ldots,X_q'$. Then
the operator $P^t$ of \eq{cs7eq2} turns out to have $q-1$
{\it small positive eigenvalues} $\la_1,\ldots,\la_{q-1}$
of size $O(t^{m-2})$. The corresponding eigenfunctions
$w_1,\ldots,w_{q-1}$ are approximately constant on the
parts of $N^t$ coming from each $X_k'$, and change
rapidly on the `small necks' in between.

As $(P^t)^{-1}w_k=\la_k^{-1}w_k$ and $\la_k=O(t^{m-2})$ we see
that $(P^t)^{-1}$ is $O(t^{2-m})$ on $\an{w_1,\ldots,w_{q-1}}$,
and so cannot be bounded independently of $t$. To repair the
proof, roughly speaking we set $W^t=\an{1,w_1,\ldots,w_{q-1}}$,
and let $(W^t)^\perp$ be the orthogonal subspace to $W^t$ in
$L^2(N^t)$. Then $P^t$ maps
\e
P^t:L^{2m}_3(N^t)\cap(W^t)^\perp\ra L^{2m}_1(N^t)\cap(W^t)^\perp
\label{cs7eq5}
\e
and has an inverse $(P^t)^{-1}$ bounded independently of $t$
on these spaces, in a weak sense. (Actually, we do something
more complicated than this, in which $W^t$ is an approximation
to $\an{1,w_1,\ldots,w_{q-1}}$ defined explicitly in terms of
bounded harmonic functions on~$L_1,\ldots,L_n$.)

In Step 7, the sequence $(f_k)_{k=0}^\iy$ is constructed as
before. The bound on the inverse of \eq{cs7eq5} can be used to
inductively bound the components of $f_k$ in $(W^t)^\perp$.
But we still need to bound the components $\pi_\sWt(f_k)$
of $f_k$ in $W^t$. Since $f_0=0$ and $Q(x,0,0)=0$, equation
\eq{cs7eq3} gives $P^tf_1=\psi^m\sin\th^t$, so that $f_1=
(P^t)^{-1}(\psi^m\sin\th^t)$. It turns out that we need
$\pi_\sWt(f_1)=o(t^2)$ for $f_k$ to remain small as~$k\ra\iy$.

As $(P^t)^{-1}=O(t^{2-m})$ on $\an{w_1,\ldots,w_{q-1}}$,
this holds if $\pi_\sWt(\psi^m\sin\th^t)=o(t^m)$.
Calculation shows that the dominant term in $\pi_\sWt
(\psi^m\sin\th^t)$ is $O(t^m)$, and proportional to
the left hand side of \eq{cs7eq4}. Therefore
$\pi_\sWt(\psi^m\sin\th^t)=o(t^m)$ if and only if
\eq{cs7eq4} holds, and this is the condition for the
sequence $(f_k)_{k=0}^\iy$ to remain bounded and
converge to a small solution of~\eq{cs7eq1}.

If $X'$ is connected, so that $q=1$, then $k(i,j)\equiv 1$
and \eq{cs7eq4} becomes
\begin{equation*}
\sum_{i=1}^n\psi(x_i)^mZ(L_i)\cdot\sum_{j=1}^{l_i}[\Si_i^j\,]=0.
\end{equation*}
But $\sum_{j=1}^{l_i}[\Si_i^j\,]=[\Si_i]$, and $Z(L_i)\cdot[\Si_i]=0$
as $Z(L_i)$ is the image of a class in $H^{m-1}(L_i,\R)$, and $\Si_i$
is the boundary of $L_i$. Therefore \eq{cs7eq4} holds automatically
when $X'$ is connected, and Theorem \ref{cs7thm2} reduces to Theorem
\ref{cs7thm1} in this case.

\subsection{Desingularization when $Y(L_i)\ne 0$}
\label{cs73}

In \cite[Th.~6.13]{Joyc12} we extend Theorem \ref{cs7thm1} to
the case $\la_i\le 0$, allowing~$Y(L_i)\ne 0$.

\begin{thm} Let\/ $(M,J,\om,\Om)$ be an almost Calabi--Yau $m$-fold
for $2\!<\!m\!<\nobreak\!6$, and\/ $X$ a compact SL\/ $m$-fold
in $M$ with conical singularities at\/ $x_1,\ldots,x_n$ and cones
$C_1,\ldots,C_n$. Let\/ $L_1,\ldots,L_n$ be Asymptotically
Conical SL\/ $m$-folds in $\C^m$ with cones $C_1,\ldots,C_n$
and rates $\la_1,\ldots,\la_n$. Suppose that\/ $\la_i\le 0$
for\/ $i=1,\ldots,n$, that\/ $X'=X\sm\{x_1,\ldots,x_n\}$ is
connected, and that there exists $\varrho\in H^1(X',\R)$
such that\/ $\bigl(Y(L_1),\ldots,Y(L_n)\bigr)$ is the image
of\/ $\varrho$ under the map $H^1(X',\R)\ra\bigoplus_{i=1}^n
H^1(\Si_i,\R)$ in \eq{cs3eq13}, where~$\Si_i=C_i\cap{\cal S}^{2m-1}$.

Then there exists $\ep>0$ and a smooth family $\bigl\{
\smash{\ti N^t}:t\in(0,\ep]\bigr\}$ of compact, nonsingular SL\/
$m$-folds in $(M,J,\om,\Om)$, such that\/ $\smash{\ti N^t}$ is
constructed by gluing $tL_i$ into $X$ at\/ $x_i$ for
$i=1,\ldots,n$. In the sense of currents, $\smash{\ti N^t}\ra
X$ as~$t\ra 0$.
\label{cs7thm3}
\end{thm}

There is also \cite[Th.~6.12]{Joyc12} an analogue of Theorem
\ref{cs7thm2}, combining the modifications of Theorems
\ref{cs7thm2} and \ref{cs7thm3}, but for brevity we will
not give it.

The new issues when $Y(L_i)\ne 0$ come mostly in Step 2
of \S\ref{cs71}. As $[\chi_i]=Y(L_i)\in H^1(\Si_i,\R)$ by
Theorem \ref{cs6thm3}, if $Y(L_i)\ne 0$ then $\chi_i$ is
no longer an {\it exact\/} form. Therefore, in Step 2 we
cannot choose a closed 1-form $\xi_i^t$ on $\Si_i\t[t^\tau,
2t^\tau]$ interpolating between $t^2\chi_i$ at $r=t^\tau$
and $\eta_i$ at $r=2t^\tau$, since $t^2\chi_i$ and $\eta_i$
have different cohomology classes.

Thus we cannot choose $N^t$ to coincide with $X$ away from
$x_i$, and work locally near $x_i$, as we did in \S\ref{cs71}.
Instead, we define $N^t$ away from $x_i$ to be $\Phi_\sXp\bigl(
\Ga(t^2\al)\bigr)$, where $\al$ is a 1-form on $X'$ satisfying
$\d\al=\d^*(\psi^m\al)=0$, and $\md{\al}=O(r^{-1})$ near $x_i$.
We show using analysis on manifolds with ends that there is a
unique such 1-form $\al$ with $[\al]=\varrho$ for each~$\varrho
\in H^1(X',\R)$.

To glue $\Up_i(tL_i\cap B_R)$ and $\Phi_\sXp\bigl(\Ga(t^2\al)
\bigr)$ together as Lagrangian $m$-folds the cohomology classes
of $t^2\chi_i$ and $t^2\al$ must agree in $H^1(\Si_i,\R)$. This
holds if the image of $\varrho$ under the map $H^1(X',\R)\ra
\bigoplus_{i=1}^nH^1(\Si_i,\R)$ is $\bigl(Y(L_1),\ldots,Y(L_n)
\bigr)$, as in the theorem. Thus, the existence of $\varrho$
with this property is a {\it necessary condition} for the
existence of $N^t$ as a {\it Lagrangian} $m$-fold in~$(M,\om)$.

After constructing $N^t$ we need to estimate norms of
$\psi^m\sin\th^t$ in Step 3. The condition $\d^*(\psi^m\al)=0$
means that linear terms in $t^2\al$ contribute 0 to
$\psi^m\sin\th^t$. However, quadratic terms in $t^2\al$
contribute $O(t^4)$ to $\psi^m\sin\th^t$ on most of $N^t$,
so all norms of $\psi^m\sin\th^t$ are at least~$O(t^4)$.

Now to show that $(f_k)_{k=0}^\iy$ converges in Step 7
we need $\lnm{\psi^m\sin\th^t}{2}=o(t^{1+m/2})$ for
small $t$. As $\lnm{\psi^m\sin\th^t}{2}$ has $O(t^4)$
contributions, this is possible only if $m<6$. Therefore
we have to restrict to complex dimension $m<6$ when
$Y(L_i)\ne 0$ for this method of proof to work.

\subsection{Desingularization in families $(M,J^s,\om^s,\Om^s)$}
\label{cs74}

Next we explain the work of \cite[\S 7--\S 8]{Joyc12} on
desingularization in {\it families} of almost Calabi--Yau
$m$-folds $(M,J^s,\om^s,\Om^s)$. The analogue of Theorem
\ref{cs7thm1} is \cite[Th.~7.15]{Joyc12}, but for brevity
we will not give it. Here \cite[Th.~7.14]{Joyc12} is the
families analogue of Theorem~\ref{cs7thm2}.

\begin{thm} Suppose $(M,J,\om,\Om)$ is an almost Calabi--Yau
$m$-fold and\/ $X$ a compact SL\/ $m$-fold in $M$ with conical
singularities at\/ $x_1,\ldots,x_n$ and cones $C_1,\ldots,C_n$.
Define $\psi:M\ra(0,\iy)$ as in \eq{cs2eq4}. Let\/ $L_1,\ldots,L_n$
be Asymptotically Conical SL\/ $m$-folds in $\C^m$ with cones
$C_1,\ldots,C_n$ and rates $\la_1,\ldots,\la_n$. Suppose $\la_i<0$
for $i=1,\ldots,n$. Write $X'=X\sm\{x_1,\ldots,x_n\}$
and\/~$\Si_i=C_i\cap{\cal S}^{2m-1}$.

Set\/ $q=b^0(X')$, and let\/ $X_1',\ldots,X_q'$ be the connected
components of\/ $X'$. For $i=1,\ldots,n$ let\/ $l_i=b^0(\Si_i)$,
and let\/ $\Si_i^1,\ldots,\Si_i^{\smash{l_i}}$ be the connected
components of\/ $\Si_i$. Define $k(i,j)=1,\ldots,q$ by $\Up_i
\circ\vp_i\bigl(\Si_i^j\t(0,R')\bigr)\subset X'_{\smash{k(i,j)}}$
for $i=1,\ldots,n$ and\/ $j=1,\ldots,l_i$. Suppose the compact\/
$m$-manifold\/ $N$ obtained by gluing $L_i$ into $X'$ at\/ $x_i$
for $i=1,\ldots,n$ is connected. A sufficient condition for this
to hold is that\/ $X$ and\/ $L_i$ for $i=1,\ldots,n$ are connected.

Suppose $\bigl\{(M,J^s,\om^s,\Om^s):s\in\F\bigr\}$ is a smooth
family of deformations of\/ $(M,J,\om,\Om)$, with base space
$\F\subset\R^d$. Let\/ $\iota_*:H_2(X,\R)\ra H_2(M,\R)$ be
the natural inclusion. Suppose that
\e
[\om^s]\cdot\iota_*(\ga)=0
\quad\text{for all\/ $s\in\F$ and\/ $\ga\in H_2(X,\R)$.}
\label{cs7eq6}
\e
Define $\G\subseteq\F\t(0,1)$ to be the subset of\/
$(s,t)\in\F\t(0,1)$ with
\e
\bigl[\Im\Om^s\bigr]\cdot\bigl[\,\smash{\ov{\!X'_k\!}}\,]=t^m
\!\!\!\!\!\!\!\!\sum_{\substack{1\le i\le n, \; 1\le j\le l_i: \\
k(i,j)=k}}\!\!\!\!\!\!\!
\psi(x_i)^mZ(L_i)\cdot[\Si_i^j\,]
\quad\text{for $k=1,\ldots,q$.}
\label{cs7eq7}
\e

Then there exist\/ $\ep\in(0,1)$ and\/ $\ka>1$ and a smooth family
\e
\bigl\{\smash{\ti N^{s,t}}:(s,t)\in\G,\quad t\in(0,\ep],
\quad \md{s}\le t^{\ka+m/2}\bigr\},
\label{cs7eq8}
\e
such that\/ $\smash{\ti N^{s,t}}$ is a compact, nonsingular
SL\/ $m$-fold in $(M,J^s,\om^s,\Om^s)$ diffeomorphic to $N$,
which is constructed by gluing $tL_i$ into $X$ at\/ $x_i$ for
$i=1,\ldots,n$. In the sense of currents in Geometric Measure
Theory, $\smash{\ti N^{s,t}}\ra X$ as~$s,t\ra 0$.
\label{cs7thm4}
\end{thm}

To prove it we modify Steps 1--7 of \S\ref{cs71} in the following
ways. In Step 1 we generalize Theorem \ref{cs4thm3} to give smooth
families of maps $\Up_i^s:B_R\ra M$ and $\Phi_\sXp^s:U_\sXp\ra M$
for small $s\in\F$ with $(\Up_i^s)^*(\om^s)=\om'$, $(\Phi_\sXp^s)^*
(\om^s)=\hat\om$ and $\Up_i^0=\Up_i$, $\Phi_\sXp^0=\Phi_\sXp$. Using
these, in Step 2 we define a smooth family of Lagrangian $m$-folds
$N^{s,t}$ in $(M,\om^s)$ for small $s\in\F$ and $t\in(0,\de)$. In
the rest of the proof we make everything depend on $s\in\F$, and
deform $N^{s,t}$ to an SL $m$-fold $\smash{\ti N^{s,t}}$ in
$(M,J^s,\om^s,\Om^s)$ for small $s\in\F$ and~$t\in(0,\de)$.

To allow $X'$ not connected, as in \S\ref{cs72}, we introduce
a vector subspace $W^{s,t}\subset C^\iy(N^{s,t})$, and we need
$\pi_\sWst(\psi^m\sin\th^{s,t})=o(t^m)$. The dominant terms in
$\pi_\sWst(\psi^m\sin\th^{s,t})$ are of two kinds: $O(t^m)$
terms involving the $Z(L_i)$, as in \S\ref{cs72}, and also terms
in $\bigl[\Im\Om^s\bigr]\cdot\bigl[\,\smash{\ov{\!X'_k\!}}\,]$.
Equation \eq{cs7eq7} requires these two terms to cancel, so that
$\pi_\sWst(\psi^m\sin\th^{s,t})=o(t^m)$, and the rest of the
proof works.

Here \cite[Th.~8.10]{Joyc12} is the families analogue of
Theorem~\ref{cs7thm3}.

\begin{thm} Let\/ $(M,J,\om,\Om)$ be an almost Calabi--Yau
$m$-fold for $2\!<\!m\!<\nobreak\!6$, and\/ $X$ a compact SL\/
$m$-fold in $M$ with conical singularities at\/ $x_1,\ldots,x_n$
and cones $C_1,\ldots,C_n$. Let\/ $L_1,\ldots,L_n$ be
Asymptotically Conical SL\/ $m$-folds in $\C^m$ with cones
$C_1,\ldots,C_n$ and rates $\la_1,\ldots,\la_n$. Suppose 
$\la_i\le 0$ for $i=1,\ldots,n$, and\/ $X'=X\sm\{x_1,\ldots,
x_n\}$ is connected.

Suppose $\bigl\{(M,J^s,\om^s,\Om^s):s\in\F\bigr\}$ is a smooth
family of deformations of\/ $(M,J,\om,\Om)$, with base space
$\F\subset\R^d$, satisfying
\e
[\Im\Om^s]\cdot[X]=0
\quad\text{for all\/ $s\in\F$, where $[X]\in H_m(M,\R)$.}
\label{cs7eq9}
\e
Define $\varpi\in H^2_{\rm cs}(X',\R)$ to be the image
of\/ $\bigl(Y(L_1),\ldots,Y(L_n)\bigr)$ under the map
$\bigoplus_{i=1}^nH^1(\Si_i,\R)\ra H^2_{\rm cs}(X',\R)$
in \eq{cs3eq13}. Define $\G\subseteq\F\t(0,1)$ to be
\e
\G=\bigl\{(s,t)\in\F\t(0,1):
[\om^s]\cdot\iota_*(\ga)=t^2\varpi\cdot\ga
\;\>\text{for all\/ $\ga\in H_2(X,\R)$}\bigr\},
\label{cs7eq10}
\e
where $\iota_*:H_2(X,\R)\ra H_2(M,\R)$ is the natural inclusion.

Then there exist\/ $\ep\in(0,1)$, $\ka>1$ and\/ $\vartheta\in(0,2)$
and a smooth family
\e
\bigl\{\smash{\ti N^{s,t}}:(s,t)\in\G,\quad t\in(0,\ep],
\quad \md{s}\le t^\vartheta\bigr\},
\label{cs7eq11}
\e
such that\/ $\smash{\ti N^{s,t}}$ is a compact, nonsingular
SL\/ $m$-fold in $(M,J^s,\om^s,\Om^s)$, which is constructed
by gluing $tL_i$ into $X$ at\/ $x_i$ for $i=1,\ldots,n$. In
the sense of currents in Geometric Measure Theory,
$\smash{\ti N^{s,t}}\ra X$ as~$s,t\ra 0$.
\label{cs7thm5}
\end{thm}

In \S\ref{cs73} we saw that when $Y(L_i)\ne 0$ there is a
topological obstruction to defining $N^t$ as a Lagrangian
$m$-fold in $(M,\om)$, so that $N^t$ exists only if the
$Y(L_i)$ satisfy an equation. In this case there is also
an obstruction to defining $N^{s,t}$ as a Lagrangian
$m$-fold in $(M,\om^s)$, but now the condition in
\eq{cs7eq10} for $N^{s,t}$ to exist involves both
$Y(L_i)$, which determine $\varpi$, and~$[\om^s]$.

Here is how to understand the relation between the
conditions for $N^t$ to exist in Theorem \ref{cs7thm3},
and for $N^{s,t}$ to exist in Theorem \ref{cs7thm5}. As
\eq{cs3eq13} is exact, $\bigl(Y(L_1),\ldots,Y(L_n)\bigr)$
is the image of $\varrho\in H^1(X',\R)$ if and only if
the image $\varpi$ of $\bigl(Y(L_1),\ldots,Y(L_n)\bigr)$
in $H^2_{\rm cs}(X',\R)$ is zero.

Now $\om^0=\om$ and $[\om]\cdot\iota_*(\ga)=0$ for all
$\ga\in H_2(X,\R)$, as $X'$ is Lagrangian in $(M,\om)$.
Thus when $s=0$, equation \eq{cs7eq10} reduces to
$t^2\varpi\cdot\ga=0$ for all $\ga$. But $H^2_{\rm cs}
(X',\R)\cong H_2(X,\R)^*$ by \eq{cs3eq15}. Thus when
$s=0$ equation \eq{cs7eq10} is equivalent to $\varpi=0$,
which is equivalent to the existence of $\varrho$ in
Theorem~\ref{cs7thm3}.

Note that as the conditions \eq{cs7eq7} and \eq{cs7eq10} for
the existence of $\smash{\ti N^{s,t}}$ involve both $s$ and
$t$, it can happen that an SL $m$-fold $X$ with conical
singularities admits no desingularizations $\smash{\ti N^t}$
in $(M,J,\om,\Om)$, but does admit desingularizations
$\smash{\ti N^{s,t}}$ in $(M,J^s,\om^s,\Om^s)$ for small
$s\ne 0$. Thus we can overcome obstructions to the
existence of desingularizations by varying the underlying
almost Calabi--Yau $m$-fold~$(M,J,\om,\Om)$.

We also prove a theorem \cite[Th.~8.9]{Joyc12} combining
Theorems \ref{cs7thm4} and \ref{cs7thm5}, desingularizing
in families when $Y(L_i)\ne 0$ and $X'$ is not connected.
However, for technical reasons it is not as strong as the
author would like, in that we must assume both sides of
\eq{cs7eq7} are zero rather than just that \eq{cs7eq7} holds.

\section{Discussion: how moduli spaces fit together}
\label{cs8}

We now consider the {\it boundary\/} $\pd\M_\sN$ of a
moduli space $\M_\sN$ of SL $m$-folds.

\begin{dfn} Let $(M,J,\om,\Om)$ be an almost Calabi--Yau
$m$-fold, $N$ a compact, nonsingular SL $m$-fold in $M$,
and $\M_\sN$ the moduli space of deformations of $N$ in $M$.
Then $\M_\sN$ is a smooth manifold of dimension $b^1(N)$, by
Theorem \ref{cs2thm3}. In general $\M_\sN$ will be a 
{\it noncompact\/} manifold, but we can construct a natural
{\it compactification} $\oM_\sN$ as follows.

Regard $\M_\sN$ as a moduli space of special Lagrangian
{\it integral currents} in the sense of Geometric Measure
Theory, as discussed in \cite[\S 6]{Joyc9}. An introduction
to Geometric Measure Theory can be found in Morgan \cite{Morg}.
Let $\oM_\sN$ be the closure of $\M_\sN$ in the space of
integral currents. As elements of $\M_\sN$ have uniformly
bounded volume, $\oM_\sN$ is {\it compact\/}
by~\cite[5.5]{Morg}.

Define the {\it boundary} $\pd\M_\sN$ to be $\oM_\sN
\sm\M_\sN$. Then elements of $\pd\M_\sN$ are {\it singular special
Lagrangian integral currents}. Essentially, they are singular SL
$m$-folds $X$ in $M$ which are limits of nonsingular $\hat N\in\M_\sN$
in an appropriate sense. By a result of Almgren, the singular set
of each $X\in\pd\M_\sN$ has Hausdorff dimension at most~$m-2$.
\label{cs8def1}
\end{dfn}

In good cases, say if $(M,J,\om,\Om)$ is suitably generic,
it seems reasonable that $\pd\M_\sN$ should be divided into
a number of {\it strata}, each of which is a moduli space of
singular SL $m$-folds with singularities of a particular type,
and is itself a manifold with singularities. In particular, some
or all of these strata could be moduli spaces $\M_\sX$ of SL
$m$-folds with isolated conical singularities, as in~\S\ref{cs5}. 

In this case, using \S\ref{cs7} for each $\hat X\in\M_\sX$ we
can try to construct desingularizations $\smash{\ti N^t}$ in
$\M_\sN$ by gluing in AC SL $m$-folds $\hat L_i$ at the
singular points $\hat x_i$ of $\hat X$ for $i=1,\ldots,n$. In
good cases, say when the cones $C_i$ of $\hat X$ are {\it stable},
every $\hat N\in\M_\sN$ close to $\M_\sX$ might be constructed
uniquely from some $\hat X,\hat L_1,\ldots,\hat L_n$, and so we
could identify an open neighbourhood of $\M_\sX$ in $\M_\sN$ with
a submanifold of the product of moduli spaces~$\M_\sX\t
\M_{\smash{\sst L_1}}^0\t\cdots\t\M_{\smash{\sst L_n}}^0$.

The goal of this section is to work towards such a description
of $\M_\sN$ near a boundary stratum $\M_\sX$ which is a moduli
space of SL $m$-folds with conical singularities. Our treatment
will be informal or conjectural in places, and is far from
giving a complete picture of~$\pd\M_\sN$.

\subsection{Topological calculations and dimension counting}
\label{cs81}

We shall consider the following situation.

\begin{dfn} Let $(M,J,\om,\Om)$ be an almost Calabi--Yau
$m$-fold for $2<m<6$. Here the assumption $m<6$ is only so that
we can apply Theorem \ref{cs7thm3} and \cite[Th.~6.12]{Joyc12},
and all of the topological calculations below actually hold
when $m>2$. Define $\psi:M\ra(0,\iy)$ as in~\eq{cs2eq4}.

Let $X$ be a compact SL $m$-fold in $M$ with conical singularities
at $x_1,\ldots,x_n$ and cones $C_1,\ldots,C_n$. Write $\Si_i=C_i
\cap{\cal S}^{2m-1}$ and $X'=X\sm\{x_1,\ldots,x_n\}$. Let $\M_\sX$
be the moduli space of deformations of $X$ in $M$, as in Definition
\ref{cs5def1}. Write $\hat X$ for a general element of $\M_\sX$. Let
$\I_\sXp$ be the image of $H^1_{\rm cs}(X',\R)$ in $H^1(X',\R)$, as
in Theorem~\ref{cs5thm1}.

Let $L_1,\ldots,L_n$ be Asymptotically Conical SL $m$-folds in
$\C^m$ with cones $C_1,\ldots,C_n$ and rate 0. Let $\M_\sLi^0$
be the moduli space of deformations of $L_i$ with rate 0, as in
Definition \ref{cs6def3}. Write $\hat L_i$ for a general element
of~$\M_\sLi^0$.

Let $q=b^0(X')$ and $X_1',\ldots,X_q'$ be the connected components
of $X'$. For $i=1,\ldots,n$ let $l_i=b^0(\Si_i)$, and let
$\Si_i^1,\ldots,\Si_i^{\smash{l_i}}$ be the connected components
of $\Si_i$. Define $k(i,j)=1,\ldots,q$ by $\Up_i\circ\vp_i\bigl(
\Si_i^j\t(0,R')\bigr)\subset X'_{\smash{k(i,j)}}$, as usual.

Define $\cY_i\subset H^1(\Si_i,\R)$ and $\cZ_i\subset H^{m-1}(\Si_i,
\R)$ to be the images of the map $H^k(L_i,\R)\ra H^k(\Si_i,\R)$ of
\eq{cs6eq2} for $k=1,m-1$. Define maps $\pi_\sYi:\M_\sLi^0\ra\cY_i$
and $\pi_\sZi:\M_\sLi^0\ra\cZ_i$ by $\pi_\sYi(\hat L_i)=Y(\hat L_i)$
and $\pi_\sZi(\hat L_i)=Z(\hat L_i)$. These are well-defined as
$Y(L_i),Z(L_i)$ are images of classes in $H^k(L_i,\R)$ by Definition
\ref{cs6def2}. Write general elements of $\cY_i$ as $\ga_i$, and of
$\cZ_i$ as~$\de_i$.

Let the vector subspace $\cY$ in $\cY_1\t\cdots\t\cY_n$ be the
intersection of $\cY_1\t\cdots\t\cY_n$ with the image of the map
$H^1(X',\R)\ra\bigoplus_{i=1}^nH^1(\Si_i,\R)$ in \eq{cs3eq13}.
Let the vector subspace $\cZ$ in $\cZ_1\t\cdots\t\cZ_n$ be the set
of all $(\de_1,\ldots,\de_n)$ for which
\e
\sum_{\substack{1\le i\le n, \; 1\le j\le l_i: \\
k(i,j)=k}}\psi(x_i)^m\de_i\cdot[\Si_i^j\,]=0
\quad\text{for all $k=1,\ldots,q$.}
\label{cs8eq2}
\e

Suppose $\smash{\bigl(Y(L_1),\ldots,Y(L_n)\bigr)}\in\cY$. This is
equivalent to the existence of $\varrho$ in Theorem \ref{cs7thm3}.
Suppose $\smash{\bigl(Z(L_1),\ldots,Z(L_n)\bigr)}\in\cZ$. This is
equivalent to equation \eq{cs7eq4} of Theorem \ref{cs7thm2}. Let
$N$ be the compact $m$-manifold obtained by gluing $L_i$ into $X'$
at $x_i$ for $i=1,\ldots,n$, as in Theorem \ref{cs7thm2}. Suppose
$N$ is connected.

Let $\smash{\ti N^t}$ for $t\in(0,\ep]$ be the desingularizations
of $X$ constructed in Theorem \ref{cs7thm2} when $Y(L_i)=0$ (as
then $L_i$ is actually AC with rate $\la_i<0$ by
\cite[Th.~7.11(b)]{Joyc9}), and in Theorem \ref{cs7thm3} when
$X'$ is connected, and in \cite[Th.~6.12]{Joyc12} in the general
case. Then each $\smash{\ti N^t}$ is a compact SL $m$-fold in
$M$ diffeomorphic to $N$. Let $\M_\sN$ be the moduli space of
deformations of $\smash{\ti N^t}$, which is independent of $t$.
Then $\M_\sN$ is a smooth manifold of dimension $b^1(N)$, by
Theorem~\ref{cs2thm3}.
\label{cs8def2}
\end{dfn}

The next four results compute the dimensions of various spaces.

\begin{lem} In the situation above we have~$\dim\I_\sXp=
b^1_{\rm cs}(X')+q-\sum_{i=1}^nl_i$.
\label{cs8lem1}
\end{lem}

\begin{proof} From \eq{cs3eq13} we see that $\I_\sXp$ fits
into an exact sequence
\begin{equation*}
0\ra H^0(X',\R)\ra\bigoplus_{i=1}^nH^0(\Si_i,\R)\ra
H^1_{\rm cs}(X',\R)\ra\I_\sXp\ra 0.
\end{equation*}
The lemma follows by alternating sum of dimensions.
\end{proof}

\begin{prop} In the situation above, $\M_\sLi^0$ is
a smooth manifold with
\e
\begin{gathered}
\dim\cY_i=b^1(L_i)-b^0(L_i)+l_i-b^1_{\rm cs}(L_i),\quad
\dim\cZ_i=l_i-b^0(L_i)\\
\text{and}\quad
\dim\M_\sLi^0=b^1(L_i)-b^0(L_i)+l_i.
\end{gathered}
\label{cs8eq3}
\e
Also the projection $\pi_\sYi\t\pi_\sZi:\M_\sLi^0\ra\cY_i\t\cZ_i$
is a smooth submersion.
\label{cs8prop1}
\end{prop}

\begin{proof} From \eq{cs6eq2} we see that $\cY_i$ fits into an
exact sequence
\begin{equation*}
0\ra H^0(L_i,\R)\ra H^0(\Si_i,\R)\ra
H^1_{\rm cs}(L_i,\R)\ra H^1(L_i,\R)\ra\cY_i\ra 0.
\end{equation*}
Taking alternating sums of dimensions gives $\dim\cY_i$ in
\eq{cs8eq3}. Similarly, $\cZ_i$ fits into an exact sequence
\begin{equation*}
0\ra\cZ_i\ra H^{m-1}(\Si_i,\R)\ra
H^m_{\rm cs}(L_i,\R)\ra H^m(L_i,\R)\ra 0.
\end{equation*}
But Poincar\'e duality gives $b^{m-1}(\Si_i)=b^0(\Si_i)=l_i$,
$b^m_{\rm cs}(L_i)=b^0(L_i)$ and $b^m(L_i)=b^0_{\rm cs}(L_i)=0$,
so we deduce $\dim\cZ_i$ in~\eq{cs8eq3}.

Suppose $0<\la_i<\min\smash{\bigl(\D_\sSii\cap(0,\iy)\bigr)}$. Then
\cite[Th.~7.11(b)]{Joyc9} shows that any AC SL $m$-fold $\hat L_i$
with cone $C_i$ and rate $\la_i$ is also Asymptotically Conical
with rate 0. Hence $\M_\sLi^{\la_i}=\M_\sLi^0$, in the notation
of Definition \ref{cs6def3}. Part (a) of Theorem \ref{cs6thm4}
then shows that $\M_\sLi^0$ is smooth with dimension given
in~\eq{cs8eq3}.

We can deduce that $\pi_\sYi\t\pi_\sZi$ is a smooth submersion from
the proof of Theorem \ref{cs6thm4} in Marshall \cite[\S 6]{Mars}.
Smoothness holds for fairly general reasons. To show that
$\pi_\sYi\t\pi_\sZi$ is a {\it submersion} we need to verify
that the natural projection $T_\sLi\M_\sLi^0\ra\cY_i\t\cZ_i$
is {\it surjective}, and this follows from the determination
of $T_\sLi\M_\sLi^0$ in~\cite[\S 5.2]{Mars}.
\end{proof}

\begin{prop} $\dim\cZ=1\!-\!q\!+\!\sum_{i=1}^n\dim\cZ_i
=1\!-\!q\!+\!\sum_{i=1}^nl_i\!-\!\sum_{i=1}^nb^0(L_i)$.
\label{cs8prop2}
\end{prop}

\begin{proof} Let $\de_i\in\cZ_i$ for $i=1,\ldots,n$. Then
$\de_i\cdot[\Si_i]=0$, since $\de_i$ is the image of a
class in $H^1(L_i)$ and $\Si_i$ is a boundary in $L_i$.
As $[\Si_i]=\sum_{j=1}^{l_i}[\Si_i^j]$, summing the left
hand side of \eq{cs8eq2} over $k=1,\ldots,q$ yields
$\sum_{i=1}^n\psi(x_i)^m\de_i\cdot[\Si_i]$, which is zero.
Thus, for any $(\de_1,\ldots,\de_n)\in\cZ_1\t\cdots\t\cZ_n$,
the sum of \eq{cs8eq2} over $k=1,\ldots,q$ holds automatically.
That is, the $q$ equations \eq{cs8eq2} on $\de_1,\ldots,\de_n$
are {\it dependent}, and represent at most $q-1$ independent
restrictions on~$\de_1,\ldots,\de_n$.

We claim that \eq{cs8eq2} is {\it exactly} $q-1$ independent
restrictions on $\de_1,\ldots,\de_n$. Then $\dim\cZ=1-q+
\sum_{i=1}^n\dim\cZ_i$, and the proposition follows from
\eq{cs8eq3}. To see this, note that $X'$ has $q$ connected
components $X_1',\ldots,X_q'$, which are joined into one
connected $N$ by gluing in $L_1,\ldots,L_n$. Define a
{\it link\/} to be a triple $(X_j',X_k',L_i^l)$, where
$1\le j<k\le q$ and $L_i^l$ is a connected component of some
$L_i$ which is glued into both $X_j'$ and $X_k'$ at~$x_i$.

Then we can choose a {\it minimal set\/} of $q-1$ links which
join $X_1',\ldots,X_q'$ into one component. It is not difficult
to show that from $\smash{(X_j',X_k',L_i^l)}$ we can construct
$\de_i\in\cZ_i$ as the image of a class in $H^{m-1}(L_i^l,\R)$,
such that the $q-1$ classes $\de_i$ obtained from the minimal
set of $q-1$ links give linearly independent left hand sides of
\eq{cs8eq2}, thought of as vectors in $\R^q$. Hence \eq{cs8eq2}
is at least $q-1$ independent restrictions on $\de_1,\ldots,\de_n$,
and the proof is complete.
\end{proof}

\begin{prop} $b^1(N)=\dim\cY+1+b^1_{\rm cs}(X')+
\sum_{i=1}^nb^1_{\rm cs}(L_i)-\sum_{i=1}^nl_i$.
\label{cs8prop3}
\end{prop}

\begin{proof} Regard $X'$ as the interior of a compact
manifold $\bar X'$ with boundary $\coprod_{i=1}^m\Si_i$,
and $L_i$ as the interior of a compact manifold $\bar L_i$
with boundary $\Si_i$. Then $N$ is constructed by gluing
$\bar X'$ and $\bar L_1,\ldots,\bar L_n$ together along
$\Si_1,\ldots,\Si_n$.

Thus the disjoint union $\coprod_{i=1}^n\Si_i$ is a subset
of $N$, with $N\sm\coprod_{i=1}^n\Si_i$ diffeomorphic to the
disjoint union of $X'$ and $L_1,\ldots,L_n$. The pair $(N;
\coprod_{i=1}^n\Si_i)$ gives an exact sequence in cohomology:
\e
\begin{split}
\cdots&\ra H^{k-1}(N,\R)\ra\bigoplus_{i=1}^n
H^{k-1}(\Si_i,\R)\ra H^k_{\rm cs}(X',\R)
\op\bigoplus_{i=1}^nH^k_{\rm cs}(L_i,\R)\\
&\ra H^k(N,\R)\ra\cdots,
\end{split}
\label{cs8eq4}
\e
since $H^k(N;\coprod_{i=1}^n\Si_i,\R)\cong H^k_{\rm cs}(X',\R)
\op\bigoplus_{i=1}^nH^k_{\rm cs}(L_i,\R)$ by excision.

Now from the definitions of $\cY_i,\cY$ and exactness of
\eq{cs3eq13} and \eq{cs6eq2} we find that the kernel of 
$\bigoplus_{i=1}^nH^1(\Si_i,\R)\ra H^2_{\rm cs}(X',\R)
\op\bigoplus_{i=1}^nH^2_{\rm cs}(L_i,\R)$ in \eq{cs8eq4} is
$\cY$. Thus as $b^0_{\rm cs}(X')=b^0_{\rm cs}(L_i)=0$ we have
an exact sequence
\begin{align*}
0&\ra H^0(N,\R)\ra\bigoplus_{i=1}^n
H^0(\Si_i,\R)\ra H^1_{\rm cs}(X',\R)
\op\bigoplus_{i=1}^nH^1_{\rm cs}(L_i,\R)\\
&\ra H^1(N,\R)\ra\cY\ra 0.
\end{align*}
The proposition follows by alternating sums of dimensions,
as~$b^0(N)=1$.
\end{proof}

\subsection{Describing the moduli space $\M_\sN$ near its boundary}
\label{cs82}

We continue to use the notation of \S\ref{cs81}. From Lemma
\ref{cs8lem1} and Propositions \ref{cs8prop1}--\ref{cs8prop3}
we deduce the following theorem. Smoothness of $\F_\sLon^\sX$
and the first line of \eq{cs8eq6} follow as $\pi_\sYi\t\pi_\sZi$
is a smooth submersion by Proposition \ref{cs8prop1}, and the
rest of \eq{cs8eq6} from the dimension formulae above.

\begin{thm} In the situation of Definition \ref{cs8def2}, define
a family $\F_\sLon^\sX$ of\/ $n$-tuples of AC SL\/ $m$-folds by
\e
\begin{split}
\F_\sLon^\sX=\Bigl\{&
(\hat L_1,\ldots,\hat L_n)\in\M_{\sst L_1}^0\t\cdots\t\M_{\sst L_n}^0:\\
&\bigl(Y(\hat L_1),\ldots,Y(\hat L_n)\bigr)\in\cY,\quad
\bigl(Z(\hat L_1),\ldots,Z(\hat L_n)\bigr)\in\cZ\Bigr\}.
\end{split}
\label{cs8eq5}
\e
Then $\F_\sLon^\sX$ is a smooth manifold with
\e
\begin{split}
\dim\F_\sLon^\sX&=\dim\cY+\dim\cZ+
\ts\sum_{i=1}^n(\dim\M_\sLi^0-\dim\cY_i-\dim\cZ_i)\\
&=\dim\cY+1-q+\ts\sum_{i=1}^nb^1_{\rm cs}(L_i)\\
&=b^1(N)-\dim\I_\sXp.
\end{split}
\label{cs8eq6}
\e
\label{cs8thm1}
\end{thm}

The significance of the theorem is that $\F_\sLon^\sX$ is the family
of $n$-tuples of AC SL $m$-folds $\smash{\hat L_1,\ldots,\hat L_n}$
which can be used to desingularize $X$ using the results of
\S\ref{cs7}. Now $\M_\sN$ is smooth with $\dim\M_\sN=b^1(N)$ by
Theorem \ref{cs2thm3}. If the cones $C_i$ are {\it stable} then
Corollary \ref{cs5cor1} shows that $\M_\sX$ is smooth with
$\dim\M_\sX=\dim\I_\sXp$. So we see from Theorem \ref{cs8thm1} that:

\begin{cor} Suppose the SL cones $C_1,\ldots,C_n$ are stable,
in the sense of Definition \ref{cs3def3}. Then the moduli
spaces $\M_\sX,\M_\sN$ and\/ $\F_\sLon^\sX$ are smooth
manifolds with\/~$\dim\M_\sX+\dim\F_\sLon^\sX=\dim\M_\sN$.
\label{cs8cor}
\end{cor}

We claim that in the stable singularities case, $\M_\sN$ is
roughly speaking locally diffeomorphic to $\M_\sX\t\F_\sLon^\sX$
near $\M_\sX\subset\pd\M_\sN$. That is, each $\hat N$ in this
region of $\M_\sN$ can be constructed from some unique $\hat
X\in\M_\sX$ and $(\hat L_1,\ldots,\hat L_n)\in\F_\sLon^\sX$
by gluing $\hat L_i$ into $\hat X$ at $\hat x_i$. This is the
reason for the formula $\dim\M_\sX+\dim\F_\sLon^\sX=\dim\M_\sN$.
To explain why, we make the following definition:

\begin{dfn} Suppose $(M,J,\om,\Om)$ is Calabi--Yau, so that
$\psi\equiv 1$. Choose $\hat X\in\M_\sX$ with singular points
$\hat x_1,\ldots,\hat x_n$ and $(\hat L_1,\ldots,\hat L_n)\in
\F_\sLon^\sX$. Then the definition of $\F_\sLon^\sX$ implies
that $\hat X,\hat L_i$ satisfy the hypotheses of Theorem
\ref{cs7thm2} if $X(\hat L_i)=0$, or Theorem \ref{cs7thm3}
if $q=1$, or \cite[Th.~6.12]{Joyc12} in the general case.
Thus, these theorems give $\ep>0$ such that for $t\in(0,\ep]$
there exists a compact SL $m$-fold $\smash{\ti N^t}$ in $M$
constructed by gluing $\hat L_i$ into $\hat X$ at~$\hat x_i$.

Observe that $(t\hat L_1,\ldots,t\hat L_n)\in\F_\sLon^\sX$ for
$t>0$. Define a subset $U$ in $\M_\sX\t\F_\sLon^\sX$ and a map
$\Psi:U\ra\M_\sN$ by $(\hat X,(t\hat L_1,\ldots,t\hat L_n))\in U$
if $t\in(0,\ep]$, where $\ep>0$ depends on $\hat X,\hat L_i$ as
above, and then~$\Psi(\hat X,(t\hat L_1,\ldots,t\hat L_n))=
\smash{\ti N^t}$.

To make $\Psi$ well-defined we have to ensure that $\smash{\ti N^t}$
is independent of choices made in its construction. Actually the
choice of $\varrho\in H^1(X',\R)$ in Theorem \ref{cs7thm3} does
affect $\smash{\ti N^t}$. Now $\varrho$ is unique up to addition
of the kernel of $H^1(X',\R)\ra\bigoplus_{i=1}^nH^1(\Si_i,\R)$ in
\eq{cs3eq13}. Choose a vector subspace of $H^1(X',\R)$ transverse
to this kernel, and restrict $\varrho$ to lie in this subspace.

This gives a way to choose $\varrho$ uniquely. Once this is done
the $N^t$ are independent of choices up to a small Hamiltonian
isotopy, and $\smash{\ti N^t}$ is the unique SL $m$-fold in this
Hamiltonian isotopy class close to $N^t$, so $\smash{\ti N^t}$
is independent of the remaining choices.
\label{cs8def3}
\end{dfn}

Here is why we assumed $M$ is Calabi--Yau above. If $M$ is
only {\it almost\/} Calabi--Yau, then $\psi$ need not be
constant. But $\cZ$ depends on $\psi(x_i)$ by \eq{cs8eq2}.
Thus, if we vary $X$ to $\hat X\in\M_\sX$ then we should
define $\F_\sLon^{\smash{\sst\hat X}}$ using $\hat\cZ$
defined with $\psi(\hat x_i)$ in \eq{cs8eq2} instead of
$\cZ$. So the family $\F_\sLon^\sX$ should vary with
$\hat X\in\M_\sX$ rather than being constant, but only
in a rather trivial way.

We claim that when the $C_i$ are stable, the map $\Psi$ is
a {\it local diffeomorphism} from the interior $U^\circ$
of $U$ to its image in $\M_\sN$. One can justify this as
follows. By \cite[\S 9.4]{Joyc8} we can define {\it natural
coordinates} on $\M_\sN$, local diffeomorphisms $\M_\sN\ra
H^1(N,\R)$ defined uniquely up to translations in $H^1(N,\R)$.
In the same way \cite[\S 6.5]{Joyc10} defines local
diffeomorphisms $\M_\sX\ra\I_\sXp$ uniquely up to translations
in $\I_\sXp$, and a similar thing applies for $\M_\sLi^0$, so
that we can construct natural coordinate systems on~$\F_\sLon^\sX$.

Using the topological calculations of \S\ref{cs81}, one can
show that for $C_i$ stable the natural coordinate systems
on $\M_\sN$ can be identified with products of the natural
coordinate systems on $\M_\sX$ and $\F_\sLon^\sX$, and $\Psi$
is just the product map in these coordinates. Thus $\Psi$ is
a local diffeomorphism on~$U^\circ$.

When the $C_i$ are not stable, things are more complicated.
Then $\M_\sX$ may be singular. If $X$ is {\it transverse},
which we expect for {\it generic} $(M,J,\om,\Om)$ by Conjecture
\ref{cs5conj}, then $\M_\sX$ is smooth of dimension $\dim\I_\sXp
-\dim\O_\sXp$ near $X$. In this case Theorem \ref{cs8thm1}
implies that $\dim\M_\sX+\dim\F_\sLon^\sX=\dim\M_\sN-\dim\O_\sXp$.
We expect $\Psi$ to be a {\it smooth immersion} wherever $\M_\sX$
is smooth, with image of codimension $\dim\O_\sXp$ in~$\M_\sN$.

Thus, when the $C_i$ are not stable and $X$ is transverse, the
desingularization results of \S\ref{cs7} do not yield the whole
of $\M_\sN$ locally, but only a subset of codimension $\dim\O_\sXp
=\sum_{i=1}^n\sind(C_i)$. Where do these extra degrees of freedom
in $\M_\sN$ come from? A rough answer is that the $\sum_{i=1}^n
\sind(C_i)$ reduction in $\dim\M_\sX$ reappears as an extra
$\sind(C_i)$ degrees of freedom to deform each $L_i$ as an AC
SL $m$-fold, but with rate $\la_i<2$ rather than rate~0.

Choose $\la_i$ with $\max(\D_\sSii\cap[0,2))<\la_i<2$. Then
by part (a) of Theorem \ref{cs6thm4}, the moduli space
$\M_\sLi^{\la_i}$ of deformations of $L_i$ with rate $\la_i$
is smooth with
\e
\dim\M_\sLi^{\la_i}\!=\!b^1(L_i)\!-\!b^0(L_i)\!+\!N_\sSii(\la_i)
\!=\!\dim\M_\sLi^0\!+\!N_\sSii(2)\!-\!m_\sSii(2)\!-\!b^0(\Si_i),
\label{cs8eq7}
\e
using the notation of Definition \ref{cs3def2}, and the
fact that $N_\sSii$ is monotone increasing and upper
semicontinuous, and increases by $m_\sSii(2)$ at~2.

Suppose $C_i$ is {\it rigid}, as in Definition \ref{cs3def3}.
Then \eq{cs3eq5} and \eq{cs8eq7} give
\e
\dim\M_\sLi^{\la_i}=\dim\M_\sLi^0+\sind(C_i)+2m.
\label{cs8eq8}
\e
Thus, deforming $L_i$ with rate $\la_i$ rather than 0 gives
an extra $\sind(C_i)+2m$ degrees of freedom. Here the $2m$
comes from {\it translations} in $\C^m$, since the AC SL
$m$-folds of rate $\la\ge 1$ are closed under translations,
and the $\sind(C_i)$ from new, nontrivial deformations of~$L_i$.

So when $X$ is transverse and the $C_i$ are rigid, we can
understand the difference in dimension $\sum_{i=1}^n\sind(C_i)$
between $\M_\sX\t\F_\sLon^\sX$ and $\M_\sN$ as coming from an
extra $\sind(C_i)$ nontrivial deformations of $L_i$ an an AC
SL $m$-fold with rate $\la_i$ rather than 0. If the $C_i$ are
not rigid, we should take into account also extra infinitesimal
deformations of $C_i$ as an SL cone.

\subsection{The index of singularities of SL $m$-folds}
\label{cs83}

We can now make more rigorous some speculations by the author
in \cite[\S 10.3]{Joyc8}. Suppose, as above, that $\M_\sN$ is
a moduli space of compact, nonsingular SL $m$-folds in
$(M,J,\om,\Om)$, and that $\M_\sX$ is a moduli space of
singular SL $m$-folds in $\pd\M_\sN$ with singularities
of a particular type, and~$X\in\M_\sX$.

Define the {\it index} of the singularities of $X$ to be
$\ind(X)=\dim\M_\sN-\dim\M_\sX$, provided $\M_\sX$ is smooth
near $X$ so $\dim\M_\sX$ is well-defined. Note that $\ind(X)$
depends not just on $X$ and its singularities, but also on $N$
through $\dim\M_\sN=b^1(N)$. Thus there could be topologically
distinct desingularizations $N_1,N_2,\ldots$ yielding different
values of~$\ind(X)$.

We can also work in {\it families} $\F$ of almost Calabi--Yau
$m$-folds $(M,J^s,\om^s,\Om^s)$. Defining $\M_\sN^\sF$ as in
Definition \ref{cs2def9} and $\M_\sX^\sF$ as in Definition
\ref{cs5def3}, the {\it index of\/ $X$ in} $\F$ is $\ind^\sF(X)
=\dim\M_\sN^\sF-\dim\M_\sX^\sF$. Note that $\ind(X)\le\dim\M_\sN
=b^1(N)$, as $\dim\M_\sX\ge 0$, and~$\ind^\sF(X)\le\dim\M_\sN^\sF$.

Combining Corollary \ref{cs5cor2} and Theorem \ref{cs8thm1},
we can compute $\ind(X)$ when $X$ is {\it transverse} with
conical singularities. In the families case, if $X$ is
{\it transverse in} $\F$ then a similar proof shows that
$\ind^\sF(X)$ is given by the same formula~\eq{cs8eq9}.

\begin{thm} Let\/ $X$ be a compact, transverse SL $m$-fold in
$(M,J,\om,\Om)$ with conical singularities at\/ $x_1,\ldots,x_n$
and cones $C_1,\ldots,C_n$. Construct desingularizations $N$
of\/ $X$ by gluing AC SL\/ $m$-folds $L_1,\ldots,L_n$ in at\/
$x_1,\ldots,x_n$, as in \S\ref{cs7}. Let\/ $q,\cY$ be as in
Definition \ref{cs8def2}. Then
\e
\ind(X)=\dim\cY+1-q+\ts\sum_{i=1}^nb^1_{\rm cs}(L_i)+
\sum_{i=1}^n\sind(C_i).
\label{cs8eq9}
\e
\label{cs8thm2}
\end{thm}

When $n=1$ this proves \cite[Conj.~2.13]{Joyc1}, in the
transverse case.

Suppose $C_i$ is not {\it rigid}, for instance if $\Si_i$ is not
connected. Then $C_i$ may lie in a smooth, connected moduli space
${\cal C}_i$ of SL cones in $\C^m$, upon which $\SU(m)$ does not
act transitively. In this case, as in \cite[\S 8.3]{Joyc10}, it is
better to define $\M_\sX$ to be the moduli space of $\hat X$ with
cones $\hat C_i\in{\cal C}_i$, rather than with fixed cones $C_i$.
Under suitable transversality assumptions, this increases the
dimension of $\M_\sX$ by the codimension of the $\SU(m)$ orbit of
$C_i$ in ${\cal C}_i$, so this codimension should be subtracted
from the r.h.s.\ of~\eq{cs8eq9}.

Here is why the index is an important idea. As $\ind(X)$ is
the codimension of $\M_\sX$ in $\M_\sN$, the {\it largest
pieces} of $\pd\M_\sN$ are the $\M_\sX$ with {\it smallest
index}. So we argue that singularities with {\it small index}
are the most generic, and the most interesting. In good cases,
we expect $\oM_\sN$ to be a {\it compact manifold with
singular boundary}. Thus the largest pieces of $\pd\M_\sN$ should
have codimension 1 in $\M_\sN$, and hence index 1. If $\pd\M_\sN$
has no index 1 strata, then $\oM_\sN$ is roughly a
compact, singular manifold {\it without boundary}.

As tools, $\ind(X)$ or $\ind^\sF(X)$ probably work best together
with a {\it genericity assumption} on $(M,J,\om,\Om)$ or $\F$.
Suppose that $\om$ is {\it generic in its K\"ahler class}. Then
as in Conjecture \ref{cs5conj}, we expect SL $m$-folds $X$ in
$(M,J,\om,\Om)$ with conical singularities to be {\it transverse},
so we can compute $\ind(X)$ using~\eq{cs8eq9}.

Since $\ind(X)\le b^1(N)$, this places {\it strong restrictions}
on the kinds of singularities that can occur in $\pd\M_\sN$. In
the families case, when $\F$ is generic in a suitable sense we
expect SL $m$-folds $X$ in $(M,J^s,\om^s,\Om^s)$ with conical
singularities to be {\it transverse in} $\F$. Then we can
compute $\ind^\sF(X)$ using \eq{cs8eq9}. Since $\ind^\sF(X)
\le\dim\M_\sN^\sF$, this again places strong restrictions on
the kinds of singularities that can occur in~$\pd\M_\sN^\sF$.

For some problems we only need to know about singularities
with index up to a certain value. For example, in \cite{Joyc1}
the author proposed to define an invariant of almost Calabi--Yau
3-folds by a weighted count of SL homology 3-spheres in a given
homology class. To understand how this invariant transforms as
we deform $(M,J,\om,\Om)$, we restrict to a generic 1-dimensional
family $\F$, and then we meet only singularities with index~1.

The index may also be useful in the {\it SYZ Conjecture} \cite{SYZ}.
This explains Mirror Symmetry of (almost) Calabi--Yau 3-folds
$M,\check M$ in terms of {\it fibrations} by SL 3-tori $N$, $\check N$.
The corresponding moduli spaces $\M_\sN$, $\M_{\smash{\sst\check N}}$
have dimension 3, so that $\pd\M_\sN$, $\pd\M_{\smash{\sst\check N}}$
can only contain singularities of index 1, 2 or~3.

\section{Applications to connected sums}
\label{cs9}

We shall now apply the results of \S\ref{cs7} to the case where
the SL $m$-fold $X$ with conical singularities is actually a
{\it nonsingular, immersed SL\/ $m$-fold}, the singular points
$x_i$ are {\it self-intersection points} of $X$ satisfying an
{\it angle criterion}, and the AC SL $m$-folds $L_i$ are chosen
from the $L^{{\bs\phi},A}$ of Example \ref{cs6ex3} due to Lawlor
\cite{Lawl}. The desingularizations $N$ are {\it multiple
connected sums} of $X$ with itself.

For the connected sum $X_1\# X_2$ of two SL 3-folds at one
point in a Calabi--Yau 3-fold, our results were conjectured
by the author in \cite[\S 5--\S 6]{Joyc1}. Butscher \cite{Buts}
proves existence of SL connected sums $X_1\# X_2$ at one point
by gluing in Lawlor necks $L^{{\bs\phi},A}$, where $X_1,X_2$
are compact SL $m$-folds in $\C^m$ with boundary.

Closer to our results, Lee \cite{Lee} considers a compact,
connected, immersed SL $m$-fold $X$ in a Calabi--Yau $m$-fold
$M$, whose self-intersection points $x_i$ satisfy an angle
criterion. She glues in $L^{{\bs\phi},A}$ at $x_i$ to get a
family of compact, {\it embedded\/} SL $m$-folds in $M$. Her
result is re-proved in Theorem \ref{cs9thm1} below.

\subsection{Transverse intersections of SL planes $\Pi^+,\Pi^-$
in $\C^m$}
\label{cs91}

Let $\Pi^+,\Pi^-$ be two SL $m$-planes $\R^m$ in $\C^m$. We call
$\Pi^\pm$ {\it transverse} if they intersect transversely, that
is, if $\Pi^+\cap\Pi^-=\{0\}$. Note that this has nothing to do
with the use of `transverse' in \S\ref{cs5}. We now classify
transverse pairs of SL $m$-planes up to the action of $\SU(m)$.
Related results may be found in \cite[\S 1]{Buts} and
\cite[\S 2]{Lee}, but we believe our approach is clearer.

\begin{prop} Let\/ $\Pi^+,\Pi^-$ be transverse SL\/ $m$-planes $\R^m$
in $\C^m$. Then there exists $B\in\SU(m)$ and\/ $\phi_1,\ldots,\phi_m
\in(0,\pi)$ such that\/ $B(\Pi^+)=\Pi^0$ and\/ $B(\Pi^-)=\Pi^{\bs\phi}$,
where ${\bs\phi}=(\phi_1,\ldots,\phi_m)$ and
\e
\Pi^0=\bigl\{(x_1,\ldots,x_m):x_j\in\R\bigr\},\;\>
\Pi^{\bs\phi}=\bigl\{({\rm e}^{i\phi_1}x_1,\ldots,
{\rm e}^{i\phi_m}x_m):x_j\in\R\bigr\},
\label{cs9eq1}
\e
as in \eq{cs6eq12}. Moreover $\phi_1,\ldots,\phi_m$ are
unique up to order, so that we can make them unique by
assuming that\/ $\phi_1\le\phi_2\le\cdots\le\phi_m$, and\/
$\phi_1+\cdots+\phi_m=k\pi$ for some~$k=1,\ldots,m-1$.
\label{cs9prop1}
\end{prop}

To prove this, first find $B'\in\SU(m)$ with $B'(\Pi^+)=\Pi^0$,
and then find $B''\in\SO(m)$ which `diagonalizes' $B'(\Pi^-)$
to get $\Pi^{\bs\phi}$, and set $B=B''B'$. The process is like
diagonalizing a real quadatic form on $\R^m$ with an $\SO(m)$
matrix. We use the proposition to divide transverse pairs
$\Pi^\pm$ into {\it types}.

\begin{dfn} Let $\Pi^+,\Pi^-$ be transverse SL $m$-planes $\R^m$
in $\C^m$. Then Proposition \ref{cs9prop1} gives $B\in\SU(m)$
and unique $0<\phi_1\le\phi_2\le\cdots\le\phi_m<\pi$ with
$\phi_1+\cdots+\phi_m=k\pi$ for some $k=1,\ldots,m-1$, such
that $B(\Pi^+)=\Pi^0$ and $B(\Pi^-)=\Pi^{\bs\phi}$. We say
that the intersection point 0 of $\Pi^+$, $\Pi^-$ is of
{\it type}~$k$.

Note that this depends on the {\it order\/} of $\Pi^+,\Pi^-$.
Exchanging $\Pi^\pm$ replaces $\phi_j$ by $\pi-\phi_{m+1-j}$
for $j=1,\ldots,m$, and therefore $\Pi^+,\Pi^-$ intersect
with type $k$ if and only if $\Pi^-,\Pi^+$ intersect with
type~$m-k$.
\label{cs9def1}
\end{dfn}

When $k=1$ in Proposition \ref{cs9prop1}, Example \ref{cs6ex3}
gives AC SL $m$-folds $L^{{\bs\phi},A}$ for $A>0$ with cone $\Pi^0
\cup\Pi^{\bs\phi}$. Thus $L^{\pm,A}=B^{-1}L^{{\bs\phi},A}$ is an
AC SL $m$-fold with cone $\Pi^+\cup\Pi^-$. When $k=m-1$ we can
exchange $\Pi^\pm$ to get $k=1$, and do the same thing. So
combining Example \ref{cs6ex3} and Proposition \ref{cs9prop1} gives:

\begin{prop} Let\/ $\Pi^+,\Pi^-$ be transverse SL\/ $m$-planes
$\R^m$ in $\C^m$, for $m>2$. Regard $C=\Pi_+\cup\Pi_-$ as an SL
cone in $\C^m$ with isolated singularity at\/ $0$. Then $\Si=C
\cap{\cal S}^{2m-1}$ is a disjoint union $\Si^+\cup\Si^-$, where
$\Si^\pm$ are the unit spheres ${\cal S}^{m-1}$ in $\Pi^\pm$. Then
\begin{itemize}
\setlength{\itemsep}{0.1pt}
\setlength{\parsep}{0.1pt}
\item[{\rm(a)}] Suppose $\Pi^+,\Pi^-$ intersect with type $1$.
Then there is a $1$-parameter family of AC SL\/ $m$-folds
$L^{\pm,A}$ in $\C^m$ asymptotic to $C$ with rate $2-m$ for
$A>0$, with\/ $Z(L^{\pm,A})\cdot[\Si^+]=A$ and\/~$Z(L^{\pm,A})
\cdot[\Si^-]=-A$.
\item[{\rm(b)}] Suppose $\Pi^+,\Pi^-$ intersect with type $m-1$.
Then there is a $1$-parameter family of AC SL\/ $m$-folds
$L^{\pm,A}$ in $\C^m$ asymptotic to $C$ with rate $2-m$ for
$A>0$, with\/ $Z(L^{\pm,A})\cdot[\Si^+]=-A$ and\/~$Z(L^{\pm,A})
\cdot[\Si^-]=A$.
\end{itemize}
The $L^{\pm,A}$ are images of the $L^{{\bf\phi},A}$ of Example
\ref{cs6ex3} under\/ $\SU(m)$ rotations.
\label{cs9prop2}
\end{prop}

We will call the $L^{\pm,A}$ {\it Lawlor necks}. They are
diffeomorphic to~${\cal S}^{m-1}\t\R$.

\subsection{Desingularizing immersed SL $m$-folds}
\label{cs92}

Here is some notation for self-intersection points of
immersed SL $m$-folds.

\begin{dfn} Let $(M,J,\om,\Om)$ be an almost Calabi--Yau
$m$-fold for $m>2$, and define $\psi:M\ra(0,\iy)$ as in
\eq{cs2eq4}. Let $X$ be a compact, nonsingular immersed
SL $m$-fold in $M$. That is, $X$ is a compact $m$-manifold
(not necessarily connected) and $\iota:X\ra M$ an immersion,
with special Lagrangian image.

Call $x\in M$ a {\it self-intersection point\/} of $X$
if $\iota^*(x)$ is at least two points in $X$. Call such
an $x$ {\it transverse} if $\iota^*(x)$ is exactly two
points $x^+,x^-$ in $X$, and $\iota_*(T_{x^+}X)\cap
\iota_*(T_{x^-}X)=\{0\}$ in~$T_xM$.

Let $x$ be a transverse self-intersection point of $X$, and
$x^\pm$ as above. Choose an isomorphism $\up:\C^m\ra T_xM$
with $\up^*(\om)=\om'$ and $\up^*(\Om)=\psi(x)^m\Om'$, where
$\om',\Om'$ are as in \eq{cs2eq2}. Define $\Pi^\pm=\up^*\bigl(
\iota_*(T_{x^\pm}X)\bigr)$. Then $\Pi^+,\Pi^-$ are transverse
SL planes $\R^m$ in~$\C^m$.

Define the {\it type} $k=1,\ldots,m-1$ of $x$ to be the type of
$\Pi^+,\Pi^-$, in the sense of Definition \ref{cs9def1}. This is
independent of the choice of $\up$, but it {\it does} depend on
the order of $x^+,x^-$, and exchanging $x^\pm$ replaces $k$ by~$m-k$.
\label{cs9def2}
\end{dfn}

We now apply Theorem \ref{cs7thm1} to desingularize {\it connected},
immersed SL $m$-folds $X$ in $M$. Our result is equivalent to Lee's
main result \cite[Th.s~3 \& 4]{Lee}, save that Lee also allows $m=2$,
and considers only the Calabi--Yau case.

\begin{thm} Let\/ $(M,J,\om,\Om)$ be an almost Calabi--Yau $m$-fold
for $m>2$, and\/ $X$ a compact, connected, immersed SL\/ $m$-fold in
$M$. Suppose $x_1,\ldots,x_n$ are transverse self-intersection
points of\/ $X$ with type $1$ or $m-1$.

Then there exists $\ep>0$ and a smooth family $\bigl\{\smash{\ti
N^t}:t\in(0,\ep]\bigr\}$ of compact, immersed SL\/ $m$-folds in
$(M,J,\om,\Om)$, such that\/ $\smash{\ti N^t}$ is constructed by
gluing a Lawlor neck $L^{\pm,t^mA_i}$ into $X$ at\/ $x_i$ for
$i=1,\ldots,n$, and so is a multiple connected sum of\/ $X$ with
itself. In the sense of currents, $\smash{\ti N^t}\ra X$ as $t\ra 0$.
If\/ $x_1,\ldots,x_n$ are the only self-intersection points of\/
$X$ then $\smash{\ti N^t}$ is embedded.
\label{cs9thm1}
\end{thm}

\begin{proof} For each $i=1,\ldots,n$ let $\iota^*(x_i)$ be
$x_i^+,x_i^-$, and define $\up_i:\C^m\ra T_{x_i}M$ and
$\Pi_i^\pm$ as in Definition \ref{cs9def2}. Then $\Pi_i^+,
\Pi_i^-$ are transverse SL planes with type 1 or $m-1$, by
assumption. Therefore Proposition \ref{cs9prop2} gives a
1-parameter family $L_i^{\smash{\pm,A}}$ for $A>0$ of AC SL
$m$-folds in $\C^m$ asymptotic to~$\Pi_i^+\cup\Pi_i^-$.

Choose some $A_1,\ldots,A_n>0$, for instance $A_i\equiv 1$,
and let $L_i=L_i^{\smash{\pm,A_i}}$ for $i=1,\ldots,n$.
Apply Theorem \ref{cs7thm1} to $X$ with conical singular
points $x_i$, cones $C_i=\Pi_i^+\cup\Pi_i^-$ and AC SL
$m$-folds $L_i$ for $i=1,\ldots,n$. As $X$ is connected
$X'=X\sm\{x_1^\pm,\ldots,x_n^\pm\}$ is connected, and
$L_i$ has rate $\la_i=2-m<0$, so the hypotheses hold.

Thus Theorem \ref{cs7thm1} gives $\ep>0$ and the family
$\smash{\ti N^t}$, and most of the theorem follows. If
$X$ has other self-intersection points than $x_1,\ldots,x_n$
then the $\smash{\ti N^t}$ are immersed. But if $x_1,\ldots,x_n$
are the only self-intersection points then $X'$ and the $L_i$ are
embedded, so the $\smash{\ti N^t}$ are embedded for small~$t$.
\end{proof}

When $m=3$ the only possible types are 1 and $m-1$, giving:

\begin{cor} Let\/ $X$ be a compact, connected, immersed SL\/
$3$-fold with transverse self-intersection points in an almost
Calabi--Yau $3$-fold $M$. Then $X$ is a limit of embedded SL\/
$3$-folds.
\label{cs9cor1}
\end{cor}

Next we apply Theorem \ref{cs7thm2} to desingularize non-connected
$X$. If a self-intersection point $x_i$ has type $m-1$ we can swap
$x_i^+,x_i^-$ to get type 1, so for simplicity we suppose the $x_i$
all have type~1.

\begin{thm} Let\/ $(M,J,\om,\Om)$ be an almost Calabi--Yau $m$-fold
for $m>2$, and\/ $X$ a compact, immersed SL\/ $m$-fold in $M$. Define
$\psi:M\ra(0,\iy)$ as in \eq{cs2eq4}. Suppose $x_1,\ldots,x_n\in M$
are transverse self-intersection points of\/ $X$ with type $1$, and
let\/ $x_i^\pm\in X$ be as in Definition \ref{cs9def2}. Set\/
$q=b^0(X)$, and let\/ $X_1,\ldots,X_q$ be the connected components
of\/ $X$. Suppose $A_1,\ldots,A_n>0$ satisfy
\e
\sum_{\substack{i=1,\ldots,n:\\ x_i^+\in X_k}}\psi(x_i)^mA_i=
\sum_{\substack{i=1,\ldots,n:\\ x_i^-\in X_k}}\psi(x_i)^mA_i
\quad\text{for all\/ $k=1,\ldots,q$.}
\label{cs9eq2}
\e
Let\/ $N$ be the oriented multiple connected sum of\/ $X$ with
itself at the pairs of points $x_i^+,x_i^-$ for $i=1,\ldots,n$.
Suppose $N$ is connected.

Then there exists $\ep>0$ and a smooth family $\bigl\{\smash{\ti
N^t}:t\in(0,\ep]\bigr\}$ of compact, immersed SL\/ $m$-folds in
$(M,J,\om,\Om)$ diffeomorphic to $N$, such that\/ $\smash{\ti N^t}$
is constructed by gluing a Lawlor neck $L^{\pm,t^mA_i}$ into $X$ at\/
$x_i$ for $i=1,\ldots,n$. In the sense of currents, $\smash{\ti N^t}
\ra X$ as $t\ra 0$. If\/ $x_1,\ldots,x_n$ are the only
self-intersection points of\/ $X$ then $\smash{\ti N^t}$ is embedded.
\label{cs9thm2}
\end{thm}

\begin{proof} Use the notation of the proof of Theorem \ref{cs9thm1}.
The assumption that $N$ is connected above is one of the hypotheses
of Theorem \ref{cs7thm2}. Part (a) of Proposition \ref{cs9prop2}
gives $Z(L_i)\cdot[\Si_i^\pm]=\pm A_i$, and using this we find that
\eq{cs7eq4} is equivalent to \eq{cs9eq2}. The other hypotheses of
Theorem \ref{cs7thm2} are established in the proof of Theorem
\ref{cs9thm1}. Thus Theorem \ref{cs7thm2} applies, and the rest
of the proof follows Theorem~\ref{cs9thm1}.
\end{proof}

To decide whether desingularizations $\smash{\ti N^t}$ of $X$
exist, we need to know when \eq{cs9eq2} admits solutions $A_i>0$.
Since $q-1$ of equations \eq{cs9eq2} are independent as in the
proof of Proposition \ref{cs8prop2}, there can only exist
{\it nonzero} solutions $A_i$ if~$n\ge q$.

Here is a {\it graphical method\/} for deciding. Draw $q$
vertices, numbered $1,\ldots,q$. For $i=1,\ldots,n$ draw a
directed edge from vertex $j$ to vertex $k$, where $x_i^+\in X_j$
and $x_i^-\in X_k$. Then \eq{cs9eq2} admits solutions $A_i>0$ if
and only if whenever we divide the $q$ vertices into two disjoint
subsets $B$ and $C$, there is at least one directed edge going
from $B$ to $C$, and at least one going from $C$ to~$B$.

As in \S\ref{cs8} we can compare the dimensions of moduli spaces
$\M_\sX,\M_\sN$ in Theorem \ref{cs9thm2}. Since $m>2$ it is easy
to show that each connected sum {\it either} reduces $b^0$ by 1
and fixes $b^1$, {\it or} fixes $b^0$ and increases $b^1$ by 1.
As $b^0(X)=q$, $b^0(N)=1$, and $\dim\M_\sX=b^1(X)$ by Theorem
\ref{cs2thm3} in the immersed case, we see that
\e
\dim\M_\sN=b^1(N)=n+1-q+b^1(X)=n+1-q+\dim\M_\sX.
\label{cs9eq3}
\e

If $n=q$ then $\dim\M_\sN=1+\dim\M_\sX$, and $\oM_\sN$
is a {\it manifold with boundary} $\M_\sX$ near $X$, and the
singularities of $X$ have {\it index} $\ind(X)=1$ in the sense of
\S\ref{cs83}. Note that \eq{cs8eq9} does not give the right answer
for $\ind(X)$ in this case, as $C_i=\Pi_i^+\cup\Pi_i^-$ is not
{\it rigid}, and to get the right definition of $\M_\sX$ we have
to allow singular cones $\hat C_i$ not just some fixed $\Pi_i^+
\cup\Pi_i^-$, but any union $\hat\Pi_i^+\cup\hat\Pi_i^-$ of
transverse SL $m$-planes $\hat\Pi_i^\pm$ in $\C^m$ with type~1.

\subsection{Desingularizing immersed SL $m$-folds in families}
\label{cs93}

Next we extend Theorem \ref{cs9thm2} to {\it families} of
almost Calabi--Yau $m$-folds. Following the proof of Theorem
\ref{cs9thm2} and using Theorem \ref{cs7thm4} we prove:

\begin{thm} Let\/ $(M,J,\om,\Om)$ be an almost Calabi--Yau $m$-fold
for $m>2$, and\/ $X$ a compact, immersed SL\/ $m$-fold in $M$ with
immersion $\iota$. Define $\psi:M\ra(0,\iy)$ as in \eq{cs2eq4}.
Suppose $x_1,\ldots,x_n\in M$ are transverse self-intersection
points of\/ $X$ with type $1$, and let\/ $x_i^\pm\in X$ be as in
Definition \ref{cs9def2}. Set\/ $q=b^0(X)$, and let\/ $X_1,\ldots,X_q$
be the connected components of\/ $X$. Let\/ $N$ be the oriented
multiple connected sum of\/ $X$ with itself at the pairs of points
$x_i^+,x_i^-$ for $i=1,\ldots,n$. Suppose $N$ is connected.

Suppose $\bigl\{(M,J^s,\om^s,\Om^s):s\in\F\bigr\}$ is a smooth
family of deformations of\/ $(M,J,\om,\Om)$ with\/ $\iota^*\bigl(
[\om^s]\bigr)=0$ in $H^2(X,\R)$ for all\/ $s\in\F$. Let\/
$A_1,\ldots,A_n>0$. Define $\G\subseteq\F\t(0,1)$ to be the
subset of\/ $(s,t)\in\F\t(0,1)$ with
\e
[\Im\Om^s]\cdot[X_k]=
t^m\sum_{\substack{i=1,\ldots,n:\\ x_i^+\in X_k}}\psi(x_i)^mA_i-
t^m\sum_{\substack{i=1,\ldots,n:\\ x_i^-\in X_k}}\psi(x_i)^mA_i
\label{cs9eq4}
\e
for all\/ $k=1,\ldots,q$. Then there exist\/ $\ep\in(0,1)$, $\ka>1$
and a smooth family
\e
\bigl\{\smash{\ti N^{s,t}}:(s,t)\in\G,\quad t\in(0,\ep],
\quad \md{s}\le t^{\ka+m/2}\bigr\},
\label{cs9eq5}
\e
such that\/ $\smash{\ti N^{s,t}}$ is a compact, nonsingular SL\/
$m$-fold in $(M,J^s,\om^s,\Om^s)$ diffeomorphic to $N$, constructed
by gluing a Lawlor neck $L^{\pm,t^mA_i}$ into $X$ at\/ $x_i$ for
$i=1,\ldots,n$. In the sense of currents, $\smash{\ti N^{s,t}}\ra X$
as $s,t\ra 0$. If\/ $x_1,\ldots,x_n$ are the only self-intersection
points of\/ $X$ then $\smash{\ti N^{s,t}}$ is embedded.
\label{cs9thm3}
\end{thm}

Thus, the main condition for the existence of desingularizations
$\smash{\ti N^{s,t}}$ of $X$ in $(M,J^s,\om^s,\Om^s)$ is that
there should exist solutions $A_1,\ldots,A_n>0$ to \eq{cs9eq4}.
Note that the sum of \eq{cs9eq4} over $k=1,\ldots,q$ gives
$[\Im\Om^s]\cdot[X]=0$, which is clearly a necessary condition
for $\smash{\ti N^{s,t}}$ to exist with~$[\smash{\ti N^{s,t}}]=[X]$.

In \cite{Joyc1} the author proposed to define an invariant
$I_3:H_3(M,\Z)\ra\Q$ of almost Calabi--Yau 3-folds $(M,J,\om,\Om)$
by counting {\it SL homology $3$-spheres} in a given homology
class with a topological weight. Theorem \ref{cs9thm3} will be
important for this programme, because it determines the {\it
transformation rules} $I_3$ satisfies as we deform $(M,J,\om,\Om)$
so that $J$ passes through certain real hypersurfaces in the
complex structure moduli space.

To explain this we need the idea of {\it SL\/ $m$-folds with
phase} ${\rm e}^{i\th}$, as in~\cite{Joyc1,Joyc2}.

\begin{dfn} Let $(M,J,\om,\Om)$ be an almost Calabi--Yau $m$-fold, and
$N$ an oriented real $m$-dimensional submanifold of $M$. Fix $\th\in\R$.
We call $N$ a {\it special Lagrangian submanifold}, or {\it SL $m$-fold}
for short, with {\it phase} ${\rm e}^{i\th}$ if
\e
\om\vert_N\equiv 0 \quad\text{and}\quad 
(\sin\th\,\Re\Om-\cos\th\,\Im\Om)\vert_N\equiv 0,
\label{cs9eq6}
\e
and $\cos\th\,\Re\Om+\sin\th\,\Im\Om$ is a positive $m$-form on 
the oriented $m$-fold~$N$.

If $N$ is {\it compact\/} it easily follows that $[\Om]\cdot[N]=
R{\rm e}^{i\th}$, where $[\Om]\in H^m(M,\R)$ and $[N]\in H_m(M,\Z)$,
and $R=\int_N\cos\th\,\Re\Om+\sin\th\,\Im\Om>0$. Thus the homology
class $[N]$ determines the phase ${\rm e}^{i\th}$ of~$N$.
\label{cs9def3}
\end{dfn}

The definition of SL $m$-fold used in the rest of the paper,
Definition \ref{cs2def6}, is of SL $m$-fold with phase 1.
If $N$ has phase ${\rm e}^{i\th}$ in $(M,J,\om,\Om)$ then
it has phase 1 in $(M,J,\om,{\rm e}^{-i\th}\Om)$, so if we
are dealing with SL $m$-folds in only one homology class then
we can rescale $\Om$ to make the phase 1. But when we consider
several SL $m$-folds $N_1,N_2,\ldots$ we cannot always take
them to have phase~1.

Using this notation, we rewrite Theorem \ref{cs9thm3} when $n=1$
and $q=2$, so that we take the {\it connected sum} $X_1\# X_2$ of
SL $m$-folds $X_1,X_2$ at one point~$x$.

\begin{thm} Let\/ $(M,J,\om,\Om)$ be an almost Calabi--Yau $m$-fold
for $m>2$, and\/ $X_1,X_2$ be compact, connected SL\/ $m$-folds in
$M$ with the same phase ${\rm e}^{i\th}$, which intersect transversely
at\/ $x\in M$ with type $1$. Suppose $\bigl\{(M,J^s,\om^s,\Om^s):
s\in\F\bigr\}$ is a smooth family of deformations of\/ $(M,J,\om,\Om)$
with\/ $\iota^*\bigl([\om^s]\bigr)=0$ in $H^2(X_k,\R)$ for all\/
$k=1,2$ and\/ $s\in\F$. Write
\e
[\Om^s]\cdot[X_k]=
R_k^s{\rm e}^{i\th_k^s}
\quad\text{for $k=1,2$ and\/}\quad
[\Om^s]\cdot\bigl([X_1]+[X_2]\bigr)=R^s{\rm e}^{i\th^s},
\label{cs9eq7}
\e
where $R_k^s,R^s>0$ and\/ $\th_k^s,\th^s\in\R$ depend continuously
on $s$ with\/ $\th_k^0=\th^0=\th$. Make $\F$ smaller if necessary
so that\/ $R_k^s,\th_k^s,R^s,\th^s$ are well-defined. Define
\e
\G=\bigl\{(s,t)\in\F\t(0,1):
R_1^s\sin(\th_1^s-\th^s)=t^m\psi(x)^m\bigr\}.
\label{cs9eq8}
\e

Then there exist\/ $\ep\in(0,1)$, $\ka>1$ and a smooth family
\e
\bigl\{\smash{\ti N^{s,t}}:(s,t)\in\G,\quad t\in(0,\ep],
\quad \md{s}\le t^{\ka+m/2}\bigr\},
\label{cs9eq9}
\e
such that\/ $\smash{\ti N^{s,t}}$ is a compact SL\/ $m$-fold with
phase ${\rm e}^{i\th^s}$ in $(M,J^s,\om^s,\Om^s)$ diffeomorphic to
$X_1\# X_2$, constructed by gluing a Lawlor neck $L^{\pm,t^m}$
into $X_1\cup X_2$ at\/ $x$. In the sense of currents,
$\smash{\ti N^{s,t}}\ra X$ as $s,t\ra 0$. If\/ $X_1,X_2$ are
embedded and\/ $x$ is their only intersection point then\/
$\smash{\ti N^{s,t}}$ is embedded.
\label{cs9thm4}
\end{thm}

This follows from Theorem \ref{cs9thm3} with $X=X_1\cup X_2$, $n=1$,
$x_1=x$ and $A_1=1$, replacing $(M,J,\om,\Om),(M,J^s,\om^s,\Om^s)$
by $(M,J,\om,{\rm e}^{-i\th}\Om),(M,J^s,\om^s,{\rm e}^{-i\th^s}\Om^s)$
so that $X_1,X_2$ and $\smash{\ti N^{s,t}}$ have phase 1. Equation
\eq{cs9eq4} becomes
\begin{equation*}
[{\rm e}^{-i\th^s}\Im\Om^s]\cdot[X_1]=t^m\psi(x)^m
\quad\text{and}\quad
[{\rm e}^{-i\th^s}\Im\Om^s]\cdot[X_2]=-t^m\psi(x)^m.
\end{equation*}
By \eq{cs9eq7} these are equivalent to
\e
R_1^s\sin(\th_1^s-\th^s)=t^m\psi(x)^m
\quad\text{and}\quad
R_2^s\sin(\th_2^s-\th^s)=-t^m\psi(x)^m.
\label{cs9eq10}
\e
But $R_1^s\sin(\th_1^s-\th^s)=-R_2^s\sin(\th_2^s-\th^s)$
as $R_1^s{\rm e}^{i\th_1^s}+R_2^s{\rm e}^{i\th_2^s}=
R^s{\rm e}^{i\th^s}$ by \eq{cs9eq7}. Thus both equations
of \eq{cs9eq10} are equivalent, so we use only the
first in~\eq{cs9eq8}.

We can interpret Theorem \ref{cs9thm4} like this. From \eq{cs9eq7}
we see that $\th^s$ always lies between $\th_1^s$ and $\th_2^s$
for small $s\in\F$. Thus making $\F$ smaller if necessary we can
divide $\F$ into three regions:
\e
\begin{gathered}
\F^+=\{s\in\F:\th_1^s>\th^s>\th_2^s\},\quad
\F^-=\{s\in\F:\th_1^s<\th^s<\th_2^s\},\\
\text{and}\quad
\F^0=\{s\in\F:\th_1^s=\th^s=\th_2^s\}.
\end{gathered}
\label{cs9eq11}
\e
If $[X_1],[X_2]$ are linearly dependent in $H_m(M,\R)$ then
$\th_1^s\equiv\th_2^s\equiv\th^s$, giving $\G=\emptyset$ in
\eq{cs9eq8}, and the theorem is trivial.

So suppose $[X_1],[X_2]$ are linearly independent. Then for $\F$
suitably generic $\F^0$ will be a {\it smooth real hypersurface}
in $\F$, which divides $\F\sm\F^0$ into two open regions $\F^\pm$.
Call $\F^+$ the {\it positive side} and $\F^-$ the {\it negative
side} of $\F^0$. Now $\th_1^s-\th^s$ is small close to $\F^0$ in
$\F$, and so $\sin(\th_1^s-\th^s)$ has the same sign as
$\th_1^s-\th^s$. Therefore $R_1^s\sin(\th_1^s-\th^s)=t^m\psi(x)^m$
in \eq{cs9eq8} admits a unique solution $t>0$ for small $s$ if
and only if $\th_1^s>\th^s$, that is, if and only if~$s\in\F^+$.

We thus have the following picture, described when $m=3$ in
Conjecture 6.5 of \cite{Joyc1}, which we have now proved. By
Theorem \ref{cs2thm4} we can extend $X_1,X_2$ to families of
SL $m$-folds $X_k^s$ with phase ${\rm e}^{i\th_k^s}$ in
$(M,J^s,\om^s,\Om^s)$ for $k=1,2$ and small $s\in\F$, such
that $X_1^s,X_2^s$ intersect transversely with type 1 at
$x^s\in M$ close to $x$. On the hypersurface $\F^0$ in $\F$
the phases of $X_1^s,X_2^s$ are equal.

On the positive side $\F^+$ of $\F^0$ there exist SL $m$-fold
connected sums $X_1^s\# X_2^s$ with phase ${\rm e}^{i\th^s}$.
On the negative side $\F^-$ there are no such $X_1^s\# X_2^s$.
Thus, as we cross hypersurfaces $\F^0$ in $\F$ where the phases
of two SL $m$-folds $X_1,X_2$ become equal, we create or destroy
new SL $m$-folds $X_1\# X_2$ by connected sum at points $x$ where
$X_1,X_2$ intersect transversely with type 1 or~$m-1$.

The conjectured invariant $I_3:H_3(M,\Z)\ra\Q$ of \cite{Joyc1}
should change in a predictable fashion as we cross hypersurfaces
$\F^0$, owing to the creation and destruction of SL homology
3-spheres. Theorem \ref{cs9thm3} also gives criteria for the
existence of multiple connected sums $X_1\# X_2\# \cdots\# X_q$
of SL $m$-folds. Using this I can derive a complete set of
transformation rules for $I_3$, and also extend the programme
to all $m\ge 3$. I hope to write a paper about this soon.

\section{Stable $T^2$-cone singularities in SL 3-folds}
\label{cs10}

We now study SL 3-folds with conical singularities modelled
on the {\it stable $T^2$-cone} $C=C_{\sst\rm HL}^3$ of Example
\ref{cs3ex1}, given by
\ea
C&=\bigl\{(z_1,z_2,z_3)\in\C^m:z_1z_2z_3\in[0,\iy),
\quad \md{z_1}=\md{z_2}=\md{z_3}\bigr\}.
\label{cs10eq1}\\
\intertext{Example \ref{cs6ex1} gives three families of AC SL
3-folds $L_{\sst\rm HL}^{\bf a}$ with rate 0 and cone $C$, which
we will write as $L_j^a$ for $j=1,2,3$ and $a>0$, given by}
L_1^a&=\bigl\{(z_1,z_2,z_3)\in\C^3:z_1z_2z_3\in[0,\iy),\quad
\ms{z_1}-a=\ms{z_2}=\ms{z_3}\bigr\},
\label{cs10eq2}\\
L_2^a&=\bigl\{(z_1,z_2,z_3)\in\C^3:z_1z_2z_3\in[0,\iy),\quad
\ms{z_1}=\ms{z_2}-a=\ms{z_3}\bigr\},
\label{cs10eq3}\\
L_3^a&=\bigl\{(z_1,z_2,z_3)\in\C^3:z_1z_2z_3\in[0,\iy),\quad
\ms{z_1}=\ms{z_2}=\ms{z_3}-a\bigr\}.
\label{cs10eq4}
\ea
Then $L_j^a$ is diffeomorphic to~${\cal S}^1\t\R^2$.

Identify $\Si=C\cap{\cal S}^5$ with $T^2=\U(1)^2$ by the map
\e
({\rm e}^{i\th_1},{\rm e}^{i\th_2})\longmapsto
\ts\bigl(\frac{1}{\sqrt{3}}\,{\rm e}^{i\th_1},
\frac{1}{\sqrt{3}}\,{\rm e}^{i\th_2},
\frac{1}{\sqrt{3}}\,{\rm e}^{-i\th_1-i\th_2}\bigr).
\label{cs10eq5}
\e
This identifies $H^1(\Si,\R)\cong H^1(T^2,\R)=\R^2$. Under
this identification, as in Example \ref{cs6ex1} we have
$Z(L_j^a)=0$ for all $j,a$ and 
\e
Y(L_1^a)=(\pi a,0),\quad Y(L_2^a)=(0,\pi a)\quad\text{and}\quad
Y(L_3^a)=(-\pi a,-\pi a).
\label{cs10eq6}
\e
For all $j,a$ define a {\it holomorphic disc} $D_j^a$ with
$\pd D_j^a\subset L_j^a$ and area $\pi a$ by
\e
D_j^a=\bigl\{(z_1,z_2,z_3)\in\C^3:\ms{z_j}\le a,\quad
\text{$z_k=0$ for $j\ne k$}\bigr\}.
\label{cs10eq7}
\e

The cone $C$ is interesting as it has three {\it topologically
distinct\/} families of AC SL 3-folds asymptotic to it, giving
three different ways to desingularize singularities of SL
3-folds with cone $C$. It is also significant as it is the
{\it only} nontrivial example of a {\it stable} SL cone in
$\C^m$ known to the author. Mark Haskins has a proof that $C$
is the only stable $T^2$-cone in $\C^3$ up to $\SU(3)$
isomorphisms (personal communication). But SL $m$-folds with
stable conical singularities have particularly good properties,
as in Corollaries \ref{cs5cor1} and~\ref{cs8cor}.

The author discussed singular SL 3-folds with cone $C$ in
\cite[\S 3--\S 4]{Joyc1} and \cite{Joyc5}. We can now prove
some of the conjectures in these papers, and give nontrivial
applications of Theorems \ref{cs7thm3} and \ref{cs7thm5}.
For simplicity we consider only SL 3-folds with one or two
singular points, but the results of \S\ref{cs7} apply to
arbitrarily many singularities~$x_1,\ldots,x_n$.

\subsection{SL 3-folds with one $T^2$-cone singularity}
\label{cs101}

We shall need a lemma on 3-manifolds with boundary.

\begin{lem} Let\/ $N$ be a compact, oriented\/ $3$-manifold
with boundary $\Si$. Then the natural map $H^1(N,\R)\ra
H^1(\Si,\R)$ has image of dimension~$\ha\,b^1(\Si)$.
\label{cs10lem1}
\end{lem}

\begin{proof} Let $N^\circ=N\sm\Si$. The pair $(N,\Si)$ has an exact
sequence in cohomology
\e
\cdots\ra H^1(N,\R)=H^1(N^\circ,\R)\ra H^1(\Si,\R)\ra
H^2_{\rm cs}(N^\circ,\R)\ra\cdots.
\label{cs10eq8}
\e
But $H^1(N,\R)\cong H^2_{\rm cs}(N^\circ,\R)^*$ and
$H^1(\Si,\R)\cong H^1(\Si,\R)^*$ by Poincar\'e duality
for $N$ and $\Si$. These isomorphisms identify the map
$H^1(N,\R)\ra H^1(\Si,\R)$ with the dual of the map
$H^1(\Si,\R)\ra H^2_{\rm cs}(N^\circ,\R)$ in \eq{cs10eq8}.
Hence the image of $H^1(N,\R)\ra H^1(\Si,\R)$ has the same
dimension as the cokernel of $H^1(\Si,\R)\ra H^2_{\rm cs}
(N^\circ,\R)$, and the lemma follows by exactness in~\eq{cs10eq8}.
\end{proof}

Consider the following situation, as in~\cite[\S 4]{Joyc1}.

\begin{dfn} Let $(M,J,\om,\Om)$ be an almost Calabi--Yau
3-fold, and $X$ a compact, connected SL 3-fold with exactly
one conical singularity at $x$, with cone $C$ in \eq{cs10eq1}.
Then $X'=X\sm\{x\}$ is connected. Let $\Si=C\cap{\cal S}^5$,
and identify $H^1(\Si,\R)\cong\R^2$ as above. Since $X'$ is
the interior of a compact, oriented 3-manifold $\bar X'$
with boundary $\Si\cong T^2$, Lemma \ref{cs10lem1} shows
that the map $H^1(X',\R)\ra H^1(\Si,\R)$ of \eq{cs3eq13}
has image~$\R$.

Similarly, the map $H^1(X',\Q)\ra H^1(\Si,\Q)$ has image
$\Q\subset\Q^2$. Thus there exist coprime integers $k_1,k_2$,
such that $(k_2,-k_1)\in\Z^2\subset\Q^2\subset\R^2$ generates
the images of $H^1(X',\R)\ra H^1(\Si,\R)$ and $H^1(X',\Q)\ra
H^1(\Si,\Q)$. Define $k_3=-k_1-k_2$. By exactness in \eq{cs3eq13},
the map $\R^2=H^1(\Si,\R)\ra H^2_{\rm cs}(X',\R)$ has kernel
$\ban{(k_2,-k_1)}$. Therefore this map is given by
\e
(y_1,y_2)\mapsto(k_1y_1+k_2y_2)\chi
\quad\text{for some nonzero $\chi\in H^2_{\rm cs}(X',\R)$.}
\label{cs10eq9}
\e
Then $k_1,k_2,k_3$ and $\chi$ are unique up to an overall
change of sign.
\label{cs10def1}
\end{dfn}

The integers $k_1,k_2,k_3$ were introduced in \cite[Def.~4.3]{Joyc1}.
We now carry out the topological calculations of \S\ref{cs8} for
desingularizing $X$ by gluing in $L_j^a$ at~$x$.

\begin{prop} In the situation of Definition \ref{cs10def1},
fix $j=1,2$ or $3$ and let\/ $N_j$ be the compact, nonsingular
$3$-manifold obtained by gluing $L_j^a$ into $X'$ at\/ $x$.
Use the notation of\/ \S\ref{cs8}, with\/ $n=1$ and\/
$L_1=L_j^a$. Then $\cZ_1=\cZ=\{0\}$ and
\begin{gather}
\dim\I_\sXp=b^1_{\rm cs}(X'),\quad
\cY_1=\begin{cases} \ban{(1,0)}, & j=1, \\ \ban{(0,1)}, & j=2, \\
\ban{(1,1)}, & j=3,\end{cases}\quad
\cY=\begin{cases}\cY_1, & k_j=0, \\ \{0\}, & k_j\ne 0, \end{cases}
\label{cs10eq10}\\
\F^\sX_{\smash{\sst L_1}}=\begin{cases}\{L_j^a:a>0\}, & k_j=0, \\
\emptyset, & k_j\ne 0, \end{cases}\quad
b^1(N_j)=\begin{cases} b^1_{\rm cs}(X')+1, & k_j=0,\\
b^1_{\rm cs}(X'), & k_j\ne 0.\end{cases}
\label{cs10eq11}
\end{gather}
\label{cs10prop1}
\end{prop}

\begin{proof} In \S\ref{cs81} as $L_1\cong{\cal S}^1\t\R^2$ and
$\Si_1\cong T^2$ we have $l_1=1$, $b^0(L_1)=1$ and $b^1_{\rm cs}
(L_1)=0$, and $q=1$ as $X'$ is connected. Thus \eq{cs8eq3} gives
$\dim\cZ_1=0$, and $\cZ\subseteq\cZ_1$, so $\cZ_1=\cZ=\{0\}$,
and Lemma \ref{cs8lem1} gives~$\dim\I_\sXp=b^1_{\rm cs}(X')$.

Now $\cY_1\subset H^1(\Si,\R)=\R^2$ is the image of
$H^1(L_{\smash{j}}^a,\R)\ra H^1(\Si,\R)$ by Definition
\ref{cs8def2}, and calculation shows it is as in \eq{cs10eq10}.
But the image of $H^1(X',\R)\ra H^1(\Si,\R)$ is $\ban{(k_2,-k_1)}$
by Definition \ref{cs10def1}, and $\cY$ is the intersection of
this image with $\cY_1$. As $k_3=-k_1-k_2$, we see that $\cY$ is
as given in~\eq{cs10eq10}.

Equation \eq{cs8eq3} gives $\dim\M_{\smash{\sst L_1}}^0=1$,
so $\M_{\smash{\sst L_1}}^0=\{L_{\smash{j}}^a:a>0\}$. But by
\eq{cs8eq5}, as $\cZ_1=\{0\}$ we see that $\F^\sX_{\smash{\sst
L_1}}$ is the subset of $\hat L_1\in\M_{\smash{\sst L_1}}^0$
with $Y(\hat L_1)\in\cY$, so \eq{cs10eq10} gives the first
equation of \eq{cs10eq11}. The second equation of
\eq{cs10eq11} follows from Proposition \ref{cs8prop3},
since $\dim\cY$ is 1 if $k_j=0$ and 0 if $k_j\ne 0$
by~\eq{cs10eq10}.
\end{proof}

Now from \S\ref{cs82}, $\F^\sX_{\smash{\sst L_1}}$ is the
family of AC SL 3-folds $L_1$ which satisfy the hypotheses
of the desingularization results of \S\ref{cs7}. Thus if $k_j=0$,
applying Theorem \ref{cs7thm3} with $L_1=L_{\smash{j}}^1$ gives:

\begin{thm} Suppose $(M,J,\om,\Om)$ is an almost Calabi--Yau
$3$-fold, and\/ $X$ a compact, connected SL\/ $3$-fold with
exactly one conical singularity at\/ $x$, with cone $C$ in
\eq{cs10eq1}. Let\/ $k_1,k_2,k_3$ be as in Definition
\ref{cs10def1}, and suppose $k_j=0$ for $j=1,2$ or $3$. Then
there exists $\ep>0$ and a smooth family $\bigl\{\smash{\ti N_j^t}:
t\in(0,\ep]\bigr\}$ of compact SL\/ $3$-folds in $(M,J,\om,\Om)$
constructed by gluing $L_{\smash{j}}^{\smash{t^2}}$ into $X$ at\/
$x$. In the sense of currents, $\smash{\ti N_j^t}\ra X$ as~$t\ra 0$.
\label{cs10thm1}
\end{thm}

If $k_j\ne 0$ then for topological reasons there exist
no {\it Lagrangian} 3-folds $N_j^t$ constructed by gluing
$tL_{\smash{j}}^1$ into $X$ at $x$, and hence no SL 3-folds
$\smash{\ti N_j^t}$. As the $k_j$ are not all zero and
$k_1+k_2+k_3=0$ there can be at most one $j$ with~$k_j=0$.

In the situation of Theorem \ref{cs10thm1}, by Corollary
\ref{cs5cor1} and \eq{cs10eq10} the moduli space $\M_\sX$
of deformations of $X$ is a smooth manifold of dimension
$b^1_{\rm cs}(X')$, and by Theorem \ref{cs2thm3} and
\eq{cs10eq11} the moduli space $\M_{\smash{\sst N_j}}$
of deformations of $\smash{\ti N_j^t}$ is a smooth manifold
of dimension~$b^1_{\rm cs}(X')+1$.

Hence the singularities of $X$ have {\it index one} in the
sense of \S\ref{cs83}, and $\oM_{\smash{\sst N_j}}$
is near $X$ a {\it nonsingular manifold with boundary} $\M_\sX$.
So in this case we have a very good understanding of the boundary
$\pd\M_{\smash{\sst N_j}}$ of $\M_{\smash{\sst N_j}}$, as
discussed in~\S\ref{cs8}.

Following \cite[\S 3.3]{Joyc1} we can explain why $\smash{\ti N_j^t}$
becomes singular as $t\ra 0$, using the $D_{\smash{j}}^a$ of
\eq{cs10eq7}. As $\smash{\ti N_j^t}$ is made by gluing
$L_{\smash{j}}^{\smash{t^2}}$ into $X$ at $x$, and there
is a holomorphic disc $D_{\smash{j}}^{\smash{t^2}}$ with
area $\pi t^2$ and $\pd D_{\smash{j}}^{\smash{t^2}}\subset
L_{\smash{j}}^{\smash{t^2}}$, we expect there to exist a
holomorphic disc $\smash{\ti D^t}$ with area $\pi t^2$ and
$\pd\smash{\ti D^t}\subset\smash{\ti N_j^t}$, for small~$t$.

As $t\ra 0$ the area of $\smash{\ti D^t}$ goes to zero, and
$\smash{\ti D^t}$ collapses to a point. Its boundary ${\cal S}^1$
in $\smash{\ti N_j^t}$ also collapses to a point, giving the singular
SL 3-fold $X$. The author expects that singularities with cone $C$
are the generic kind of singularity of SL 3-folds occurring when
areas of holomorphic discs become zero.

\subsection{SL 3-folds with one $T^2$-cone singularity in families}
\label{cs102}

Next we apply Theorem \ref{cs7thm5} to desingularize SL 3-folds
with one singularity $x$ with cone $C$ in {\it families} of
almost Calabi--Yau 3-folds~$(M,J^s,\om^s,\Om^s)$.

\begin{thm} Suppose $(M,J,\om,\Om)$ is an almost Calabi--Yau
$3$-fold, and\/ $X$ a compact, connected SL\/ $3$-fold with
exactly one conical singularity at\/ $x$, with cone $C$ in
\eq{cs10eq1}. Let\/ $k_1,k_2,k_3,\chi$ be as in Definition
\ref{cs10def1}. Suppose $\bigl\{(M,J^s,\om^s,\Om^s):s\in\F\bigr\}$
is a smooth family of deformations of\/ $(M,J,\om,\Om)$ with
\e
[\Im\Om^s]\cdot[X]=0
\quad\text{for all\/ $s\in\F$, where $[X]\in H_3(M,\R)$.}
\label{cs10eq12}
\e
Let\/ $\iota_*:H_2(X,\R)\ra H_2(M,\R)$ be the inclusion,
fix $j=1,2$ or $3$, and define
\e
\G_j\!=\!\bigl\{(s,t)\!\in\!\F\!\t\!(0,1):
[\om^s]\cdot\iota_*(\ga)\!=\!\pi t^2k_j(\chi\cdot\ga)
\;\>\text{for all\/ $\ga\in H_2(X,\R)$}\bigr\}.
\label{cs10eq13}
\e

Then there exist\/ $\ep\in(0,1)$, $\ka>1$ and\/ $\vartheta\in(0,2)$
and a smooth family
\e
\bigl\{\smash{\ti N^{s,t}_{\smash{j}}}:(s,t)\in\G_j,\quad
t\in(0,\ep], \quad \md{s}\le t^\vartheta\bigr\},
\label{cs10eq14}
\e
such that\/ $\smash{\ti N^{s,t}_{\smash{j}}}$ is a compact
SL\/ $3$-fold in $(M,J^s,\om^s,\Om^s)$ constructed by gluing
$L_{\smash{j}}^{\smash{t^2}}$ into $X$ at\/ $x$. In the sense of
currents, $\smash{\ti N^{s,t}_{\smash{j}}}\ra X$ as~$s,t\ra 0$.
\label{cs10thm2}
\end{thm}

\begin{proof} We apply Theorem \ref{cs7thm5}, with $L_1=L_{\smash{j}}^1$.
As $X$ is connected, $X'$ is connected. Equation \eq{cs10eq12}
gives \eq{cs7eq9}. The values of $Y(L_{\smash{j}}^a)$ in $H^1(\Si,\R)
\cong\R^2$ are given in \eq{cs10eq6}, and the map $H^1(\Si,\R)
\ra H^2_{\rm cs}(X',\R)$ is given in \eq{cs10eq9}. Combining
these two shows that the image of $Y(L_{\smash{j}}^a)$ in $H^2_{\rm cs}
(X',\R)$ is $\pi ak_j\,\chi$. Thus putting $a=1$ as $L_1=L_{\smash{j}}^1$,
we have $\varpi=\pi k_j\,\chi$ in Theorem \ref{cs7thm5}, and so
$\G_j$ in \eq{cs10eq13} agrees with $\G$ in \eq{cs7eq10}. The
result then follows from Theorem~\ref{cs7thm5}.
\end{proof}

We now specialize to the case that $b^1_{\rm cs}(X')=0$.
Then $b^2(X')=0$ by \eq{cs3eq14}, so \eq{cs3eq13} gives
an exact sequence
\begin{equation*}
0\ra H^1(X',\R)\ra H^1(\Si_i,\R)\cong\R^2\ra H^2_{\rm cs}
(X',\R)\ra 0.
\end{equation*}
As $b^1(X')=b^2_{\rm cs}(X')$ by \eq{cs3eq14}, this
gives~$b^1(X')=b^2_{\rm cs}(X')=1$.

Now $H_2(X,\R)\cong H^2_{\rm cs}(X',\R)^*\cong\R$ by
\eq{cs3eq15}, and $\chi\ne 0$ by \eq{cs10eq9}. Thus there
exists a unique $\ga_0\in H_2(X,\R)$ with $\chi\cdot\ga_0=1$.
The condition $[\om^s]\cdot\iota_*(\ga)=\pi t^2k_j(\chi\cdot\ga)$
for all $\ga\in H_2(X,\R)$ in \eq{cs10eq13} then becomes the
single real equation~$[\om^s]\cdot\iota_*(\ga_0)=\pi t^2k_j$.

As in \eq{cs9eq11}, let us divide the family $\F$ into three regions:
\e
\begin{gathered}
\F^+=\bigl\{s\in\F:[\om^s]\cdot\iota_*(\ga_0)>0\bigr\},\quad
\F^-=\bigl\{s\in\F:[\om^s]\cdot\iota_*(\ga_0)<0\bigr\},\\
\text{and}\quad
\F^0=\bigl\{s\in\F:[\om^s]\cdot\iota_*(\ga_0)=0\bigr\}.
\end{gathered}
\label{cs10eq15}
\e
If $\iota_*(\ga_0)\ne 0$ and $\F$ is sufficiently generic,
then $\F^0$ will be a {\it smooth real hypersurface} in $\F$,
which divides $\F\sm\F^0$ into two open regions $\F^\pm$. Call
$\F^+$ the {\it positive side} and $\F^-$ the {\it negative
side} of~$\F^0$.

We investigate the existence of deformations $X^s$ and
desingularizations $\smash{\ti N^{s,t}_{\smash{j}}}$ of $X$ in
$(M,J^s,\om^s,\Om^s)$, for small $s$ in each of these regions.
\begin{itemize}
\setlength{\itemsep}{1pt}
\setlength{\parsep}{1pt}
\item[(a)] From \S\ref{cs5}, a {\it necessary condition}
for the existence of any SL 3-fold $X^s$ with a conical
singularity isotopic to $X$ in $(M,J^s,\om^s,\Om^s)$ is that
$\iota_*(\ga)\cdot[\om^s]=0$ for all $\ga\in H_2(X,\R)$.
Thus, such $X^s$ can exist only if~$s\in\F^0$.

As $\dim\I_\sXp=b^1_{\rm cs}(X')=0$, Corollary \ref{cs5cor3} then
shows that for {\it small\/} $s\in\F^0$, there is a {\it unique}
deformation $X^s$ of $X$ close to $X$ in~$(M,J^s,\om^s,\Om^s)$.
\item[(b)] If $s\in\F^0$ then $[\om^s]\cdot\iota_*(\ga_0)=\pi t^2k_j$
has solutions $t>0$ if and only if $k_j=0$, and then any $t>0$ is a
solution. For small $s\in\F^0$ we may apply Theorem \ref{cs10thm1}
to the unique $X^s$ above to get a 1-parameter family of SL 3-folds
$\smash{\ti N^{s,t}_{\smash{j}}}$ in $(M,J^s,\om^s,\Om^s)$ for
small $t>0$, with~$b^1(\smash{\ti N^{s,t}_{\smash{j}}})=b^1(N_j)=1$.
\item[(c)] If $s\in\F^+$ then $[\om^s]\cdot\iota_*(\ga_0)=\pi t^2k_j$
has a unique solution $t>0$ if and only if $k_j>0$. Theorem
\ref{cs10thm2} then shows that a desingularization $\smash{\ti
N^{s,t}_{\smash{j}}}$ exists provided $t\le\ep$ and $\md{s}\le
t^\vartheta$, which is unique as~$b^1(\smash{\ti N^{s,t}_{
\smash{j}}})=b^1(N_j)=0$.

By applying Theorem \ref{cs10thm2} not just to $X$ but to $X^s$
for small $s\in\F^0$, one can show that such $\smash{\ti N^{s,t}_{
\smash{j}}}$ actually exist for all small~$s\in\F^+$.
\item[(d)] If $s\in\F^-$ then $[\om^s]\cdot\iota_*(\ga_0)=\pi t^2k_j$
has a unique solution $t>0$ if and only if $k_j<0$. As for $s\in
\F^+$, we find that a unique desingularization $\smash{\ti
N^{s,t}_{\smash{j}}}$ then exists for small~$s\in\F^-$.
\end{itemize}

As $k_1,k_2,k_3$ are not all zero with $k_1+k_2+k_3=0$, there
is at least one $k_j<0$, and at least one $k_j>0$. Suppose 
$k_1<0$ and $k_2,k_3>0$, for example. Imagine $s\in\F$ moving
along a curve near 0 starting in $\F^-$, crossing $\F^0$ and
ending in $\F^+$. Then initially there is one SL homology
3-sphere $\smash{\ti N^{s,t}_1}$ in $(M,J^s,\om^s,\Om^s)$.
As $s$ crosses $\F^0$ this SL 3-fold collapses to a singular
SL 3-fold $X^s$, with a conical singularity with cone $C$.

As $s$ moves into $\F^+$ it is desingularized in {\it two
topologically distinct ways} to give two SL homology
3-spheres $\smash{\ti N^{s,t}_2},\smash{\ti N^{s,t}_3}$.
We have found a process by which {\it one SL homology
$3$-sphere can turn into two SL homology $3$-spheres}
as we deform the underlying almost Calabi--Yau 3-fold
$(M,J^s,\om^s,\Om^s)$. This was described in
\cite[\S 4.2]{Joyc1}, and we have now
proved~\cite[Conj.~4.4]{Joyc1}.

In \cite[Prop.~4.5]{Joyc1} we compute $H_1(X,\Z)$,
$H_1(N_j,\Z)$ and show:

\begin{prop} In the situation of Definition \ref{cs10def1},
let\/ $N_j$ be the compact, nonsingular $3$-manifold obtained
by gluing $L_{\smash{j}}^a$ into $X'$ at\/ $x$, for $j=1,2,3$.
Suppose $b^1_{\rm cs}(X')=0$. Then $H_1(X,\Z)$ is finite.
If\/ $k_j\ne 0$ then $H_1(N_j,\Z)$ is also finite
with\/~$\md{H_1(N_j,\Z)}=\md{k_j}\cdot\md{H_1(X,\Z)}$.
\label{cs10prop2}
\end{prop}

In the situation above with $k_1<0$ and $k_2,k_3>0$, as
$k_1+k_2+k_3=0$ we have $\md{k_1}=\md{k_2}+\md{k_3}$,
and therefore $\md{H_1(N_1,\Z)}=\md{H_1(N_2,\Z)}+
\md{H_1(N_3,\Z)}$ by Proposition \ref{cs10prop2}.
Now $\md{H_1(N_j,\Z)}$ is the {\it number of flat\/
$\U(1)$-connections} on $N_j$. Thus when $\smash{\ti
N^{s,t}_1}$ turns into $\smash{\ti N^{s,t}_2}$ and
$\smash{\ti N^{s,t}_3}$, the number of SL homology
3-spheres with flat $\U(1)$-connections {\it does not change}.

This is physically significant since 3-{\it branes} in String
Theory correspond in a classical limit to SL 3-folds with flat
$\U(1)$-connections, as in \cite{SYZ} for instance. The proposal
of \cite{Joyc1} is to count SL homology 3-spheres with flat
$\U(1)$-connections in a given homology class. We have shown
that this number is {\it conserved\/} under a nontrivial kind
of transition in the family of SL homology 3-spheres.

\subsection{SL 3-folds with two $T^2$-cone singularities}
\label{cs103}

Next we study SL 3-folds with {\it two} singularities with cone~$C$.

\begin{dfn} Let $(M,J,\om,\Om)$ be an almost Calabi--Yau
3-fold, and $X$ a compact, connected SL 3-fold with exactly
two conical singularities at $x_1,x_2$, both with cone $C$
in \eq{cs10eq1}. Then $X'=X\sm\{x_1,x_2\}$ is connected. Write
$\Si_1,\Si_2$ for the two copies of $\Si=C\cap{\cal S}^5$ at
$x_1,x_2$, and identify $H^1(\Si_i,\R)\cong\R^2$ as above.
Write elements of $H^1(\Si_1,\R)\op H^1(\Si_2,\R)=\R^2\op\R^2$
as~$\bigl((u,v),(y,z)\bigr)$.

Since $X'$ is the interior of a compact, oriented 3-manifold
$\bar X'$ with boundary $\Si_1\cup\Si_2$, the map $H^1(X',\R)
\ra H^1(\Si_1,\R)\op H^1(\Si_2,\R)$ of \eq{cs3eq13} has image
$\R^2$ by Lemma \ref{cs10lem1}. Choose a basis $\bigl((u_1,v_1),
(y_1,z_1)\bigr)$, $\bigl((u_2,v_2),(y_2,z_2)\bigr)$ for this image.
As it is also a basis over $\R$ for the image of $H^1(X',\Q)\ra
H^1(\Si_1,\Q)\op H^1(\Si_2,\Q)$ we can take~$u_1,\ldots,z_2\in\Q$.

Let $\al_1,\al_2$ be closed 1-forms on $\bar X'$ such that the
images of $[\al_1],[\al_2]$ in $H^1(\Si_1,\R)\op H^1(\Si_2,\R)$
are this basis. Then $\al_1\w\al_2$ is a closed 2-form on $\bar X'$,
an oriented 3-manifold with boundary $\Si_1\cup\Si_2$, so by
Stokes' Theorem we have $\int_{\Si_1\cup\Si_2}\al_1\w\al_2=0$.
This gives the consistency condition
\e
u_1v_2-u_2v_1+y_1z_2-y_2z_1=0.
\label{cs10eq16}
\e
\label{cs10def2}
\end{dfn}

Applying Theorem \ref{cs7thm3} gives a necessary and sufficient
criterion for when we can desingularize $X$ by gluing in AC SL
3-folds $L_{\smash{j_1}}^{\smash{a_1}},L_{\smash{j_2}}^{\smash{a_2}}$
at~$x_1,x_2$.

\begin{thm} Suppose $(M,J,\om,\Om)$ is an almost Calabi--Yau
$3$-fold, and\/ $X$ a compact, connected SL\/ $3$-fold with
exactly two conical singularities at\/ $x_1,x_2$, both with
cone $C$ in \eq{cs10eq1}. Let\/ $u_1,\ldots,z_2$ be as in
Definition \ref{cs10def2}. Choose $j_1,j_2=1,2,3$ and\/
$a_1,a_2>0$, and set\/ $L_i=L_{\smash{j_i}}^{\smash{a_i}}$
for $i=1,2$. Then there exists $\ep>0$ and a smooth family
$\bigl\{\smash{\ti N^t}:t\in(0,\ep]\bigr\}$ of compact SL\/
$3$-folds in $(M,J,\om,\Om)$ constructed by gluing $tL_i$
into into $X$ at\/ $x_i$ if and only if
\e
\bigl(Y(L_1),Y(L_2)\bigr)\in\ban{\bigl((u_1,v_1),(y_1,z_1)\bigr),
\bigl((u_2,v_2),(y_2,z_2)\bigr)}\subset\R^2\op\R^2,
\label{cs10eq17}
\e
where $Y(L_i)$ are given in \eq{cs10eq6}. In the sense of
currents, $\smash{\ti N^t}\ra X$ as~$t\ra 0$.
\label{cs10thm3}
\end{thm}

Here $\varrho$ exists in Theorem \ref{cs7thm3} if and
only if $\bigl(Y(L_1),Y(L_2)\bigr)$ lies in the image of
$H^1(X',\R)\ra H^1(\Si_1,\R)\op H^1(\Si_2,\R)$, that is,
if and only if \eq{cs10eq17} holds. We are interested in
how many {\it topologically distinct\/} ways of
desingularizing $X$ there are, and in the {\it index} of
the singularities of~$X$.

Let us use the notation of \S\ref{cs8}, so that $\I_\sXp$
is the image of $H^1_{\rm cs}(X',\R)$ in $H^1(X',\R)$, and
$\M_\sX$ the moduli space of deformations of $X$ in
$(M,J,\om,\Om)$ as in \S\ref{cs5}, and $N$ the compact
3-manifold obtained by gluing $L_i=L_{\smash{j_i}}^{\smash{a_i}}$
into $X$ at $x_i$ for $i=1,2$, and $\M_\sN$ the moduli space
of deformations of $\smash{\ti N^t}$ in~$(M,J,\om,\Om)$.

Corollary \ref{cs5cor1} and Theorem \ref{cs2thm3} show that
$\M_\sX$ and $\M_\sN$ are smooth with
\e
\dim\M_\sX=\dim\I_\sXp
\quad\text{and}\quad
\dim\M_\sN=b^1(N).
\label{cs10eq18}
\e
As $q=l_1=l_2=1$, Lemma \ref{cs8lem1} gives $\dim\I_\sXp=
b^1_{\rm cs}(X')-1$. Equation \eq{cs8eq3} shows that $\cZ_i=\cZ=
\{0\}$ and $\dim\cY_i=1$, so that $\cY_i=\an{Y(L_i)}$. Therefore
\e
\cY=\ban{\bigl(Y(L_1),0\bigr),\bigl(0,Y(L_2)\bigr)}\cap
\ban{\bigl((u_1,v_1),(y_1,z_1)\bigr),\bigl((u_2,v_2),(y_2,z_2)\bigr)}.
\label{cs10eq19}
\e
Proposition \ref{cs8prop3} then gives
\e
b^1(N)=\dim\cY+b^1_{\rm cs}(X')-1=\dim\cY+\dim\I_\sXp,
\label{cs10eq20}
\e
so that $\dim\M_\sN\!=\!\dim\M_\sX\!+\!\dim\cY$, and
$\ind(X)\!=\!\dim\cY$ in the sense of~\S\ref{cs83}.

For generic $u_1,\ldots,z_2$ there will be no choices of
$j_1,j_2,a_1,a_2$ for which \eq{cs10eq17} holds, and so
no way to desingularize $X$ in $(M,J,\om,\Om)$. Here are
four examples where other things happen.

\begin{ex} Suppose $u_1,\ldots,z_2$ are given by
\e
\bigl((u_1,v_1),(y_1,z_1)\bigr)\!=\!\bigl((1,0),(0,0)\bigr),\;\>
\bigl((u_2,v_2),(y_2,z_2)\bigr)\!=\!\bigl((0,0),(1,0)\bigr),
\label{cs10eq21}
\e
so that \eq{cs10eq16} holds. Then we can desingularize
$X$ using $L_1=L_{\smash{1}}^{\smash{a_1}}$ and $L_2=
L_{\smash{1}}^{\smash{a_2}}$ for any $a_1,a_2>0$, with
$\ind(X)=\dim\cY=2$. We must have $j_1=j_2=1$, so there
is only one topological possibility.
\label{cs10ex1}
\end{ex}

\begin{ex} Let $\bigl((u_1,v_1),(y_1,z_1)\bigr)=\bigl((1,0),(r,0)
\bigr)$ for $r>0$ in $\Q$, and $\bigl((u_2,v_2),(y_2,z_2)\bigr)$
be generic with $v_2\!+\!rz_2\!=\!0$. The only way to desingularize
$X$ is with $L_1=L_{\smash{1}}^{\smash{a}}$ and $L_2=L_{\smash{1}
}^{\smash{ra}}$ for $a>0$, with~$\ind(X)=\dim\cY=1$.
\label{cs10ex2}
\end{ex}

\begin{ex} Let $r>0$ be in $\Q$ with $r\ne 1$, and let
\e
\bigl((u_1,v_1),(y_1,z_1)\bigr)\!=\!\bigl((1,0),(0,r)\bigr),\;\>
\bigl((u_2,v_2),(y_2,z_2)\bigr)\!=\!\bigl((0,r),(1,0)\bigr).
\label{cs10eq22}
\e
Then there are exactly {\it two topologically distinct\/} ways
to desingularize~$X$:
\begin{itemize}
\setlength{\itemsep}{0pt}
\setlength{\parsep}{0pt}
\item[(a)] $j_1=1$, $j_2=2$, $a_1=a>0$, $a_2=ra>0$,
$L_1=L_{\smash{1}}^{\smash{a}}$ and
$L_2=L_{\smash{2}}^{\smash{ra}}$,
\item[(b)] $j_1=2$, $j_2=1$, $a_1=ra>0$, $a_2=a>0$,
$L_1=L_{\smash{2}}^{\smash{ra}}$ and
$L_2=L_{\smash{1}}^{\smash{a}}$.
\end{itemize}
Both have $\ind(X)=\dim\cY=1$.
\label{cs10ex3}
\end{ex}

\begin{ex} Suppose $u_1,\ldots,z_2$ are given by
\e
\bigl((u_1,v_1),(y_1,z_1)\bigr)\!=\!\bigl((1,0),(0,1)\bigr),\;\>
\bigl((u_2,v_2),(y_2,z_2)\bigr)\!=\!\bigl((0,1),(1,0)\bigr),
\label{cs10eq23}
\e
which is the case $r=1$ in Example \ref{cs10ex3}. Then there are
exactly {\it three topologically distinct\/} ways to desingularize
$X$, each with~$\ind(X)=\dim\cY=1$:
\begin{itemize}
\setlength{\itemsep}{0pt}
\setlength{\parsep}{0pt}
\item[(a)] $j_1=1$, $j_2=2$, $a_1=a_2=a>0$,
$L_1=L_{\smash{1}}^{\smash{a}}$ and
$L_2=L_{\smash{2}}^{\smash{a}}$,
\item[(b)] $j_1=2$, $j_2=1$, $a_1=a_2=a>0$,
$L_1=L_{\smash{2}}^{\smash{a}}$ and
$L_2=L_{\smash{1}}^{\smash{a}}$,
\item[(c)] $j_1=j_2=3$, $a_1=a_2=a>0$, and
$L_1=L_2=L_{\smash{3}}^{\smash{a}}$.
\end{itemize}
\label{cs10ex4}
\end{ex}

We can explain this using the {\it holomorphic discs} $D_j^a$
of \eq{cs10eq7}. For each desingularization $\smash{\ti N^t}$
when $t$ is small, we expect there to exist unique holomorphic
discs $\smash{\ti D_1^t,\ti D_2^t}$ in $(M,J)$, where
$\smash{\ti D_i^t}$ is near $x_i$, has area $\pi a_i$ and
boundary $\pd\smash{\ti D_i^t}\subset\smash{\ti N^t}$
for~$i=1,2$.

In Example \ref{cs10ex1} the homology classes
$[\pd\smash{\ti D_1^t}],[\pd\smash{\ti D_2^t}]\in
H_1(\smash{\ti N^t},\R)$ are {\it linearly independent}.
Therefore by deforming $\smash{\ti N^t}$ as an SL 3-fold
we can vary the areas of $\smash{\ti D_1^t,\ti D_2^t}$
independently. These two areas give two real parameters
for the desingularization of $X$, which is why~$\ind(X)=2$.

However, in Examples \ref{cs10ex2}--\ref{cs10ex4} the
homology classes $[\pd\smash{\ti D_1^t}],[\pd\smash{\ti D_2^t}]$
are {\it proportional\/} in $H_1(\smash{\ti N^t},\R)$. This
forces $\area(\smash{\ti D_2^t})=c\cdot\area(\smash{\ti D_1^t})$
to hold under deformations of $\smash{\ti N^t}$, where $c=r$ in
Example \ref{cs10ex2} and part (a) of Example \ref{cs10ex3},
$c=r^{-1}$ in part (b) of Example \ref{cs10ex3}, and $c=1$
in Example \ref{cs10ex4}. In particular, the areas of
$\smash{\ti D_1^t,\ti D_2^t}$ can only become zero
{\it simultaneously}.

Therefore the two singularities $x_1,x_2$ in $X$ are not
independent, but {\it coupled together}. For an SL $m$-fold
$X$ with conical singularities $x_1,\ldots,x_n$ one might
na\"\i vely expect each $x_i$ to be desingularized separately,
and $\ind(X)$ to be the sum of contributions from each $x_i$.
But in Examples \ref{cs10ex2}--\ref{cs10ex4} we see that
$x_1,x_2$ can only be desingularized together, and $\ind(X)=1$
is not a sum of separate contributions from~$x_1,x_2$.

In Examples \ref{cs10ex2}--\ref{cs10ex4} the moduli spaces
$\M_\sX$, $\M_\sN$ are smooth with $\dim\M_\sN=\dim\M_\sX+1$.
Therefore $\oM_\sN$ is near $X$ a {\it nonsingular
manifold with boundary} $\M_\sX$. So we have a good
understanding of the boundary $\pd\M_\sN$ of $\M_\sN$, as
in \S\ref{cs8}. It is also interesting that in Examples
\ref{cs10ex3} and \ref{cs10ex4} we have {\it two or three
different moduli spaces $\M_\sN$ with common boundary}~$\M_\sX$.

Finally we discuss the ideas of \cite{Joyc5} on {\it SL
fibrations} of (almost) Calabi--Yau 3-folds, as required by
the {\it SYZ Conjecture} \cite{SYZ}. Let $(M,J,\om,\Om)$
be an almost Calabi--Yau 3-fold and $f:M\ra B$ an {\it SL
fibration}. That is, $B$ is a compact 3-manifold, $f$ is
continuous and piecewise smooth, and for some $\De\subset B$
with $B\sm\De$ open and dense the fibres $X_b=f^{-1}(b)$ are SL
3-tori for $b\in B\sm\De$, and singular SL 3-folds for~$b\in\De$.

The idea of \cite[\S 7--\S 8]{Joyc5} is that for $\om$ generic in
its K\"ahler class, $\De$ should be of codimension 1 in $B$, and
for $b\in\De$ generic $X_b$ should have 2 (or $2n$) singularities
with cone $C$, as in Definition \ref{cs10def2}. Then $\De$ is
locally a hypersurface which divides $B\sm\De$ into two pieces.
These are {\it two different moduli spaces} of SL 3-tori with
{\it common boundary} $\De$, as in Examples \ref{cs10ex3} and
\ref{cs10ex4}. Our results provide a partial proof of
speculations in~\cite[\S 8.2(a)]{Joyc5}.

\end{document}